%% file: Non_local_Arxiv.tex
\documentclass{amsart}
\usepackage{graphicx}

\usepackage{appendix}
\usepackage{amssymb}
\usepackage{amsmath}

\usepackage{wrapfig}
\usepackage{subfig}
\usepackage{enumerate}
\usepackage{fancybox}
\usepackage{framed}
\usepackage{xcolor}
\usepackage{mathrsfs}
\usepackage{esint}
\usepackage{float}
\usepackage{color}
\usepackage[top=2.5cm,bottom=2.5cm,left=2.5cm,right=2.5cm]{geometry}
\usepackage[bordercolor=white,backgroundcolor=gray!30,linecolor=black,colorinlistoftodos]{todonotes}
\usepackage{hyperref}

 \numberwithin{equation}{section}

\mathchardef\mhyphen="2D

\input{macros}
\DeclareFontFamily{OT1}{pzc}{}
\DeclareFontShape{OT1}{pzc}{m}{it}%
              {<-> s * [0.900] pzcmi7t}{}
\DeclareMathAlphabet{\mathpzc}{OT1}{pzc}%
                                 {m}{it}

\title{Analysis and Application of a non-local Hessian} 

\author{Jan Lellmann, Konstantinos Papafitsoros, Carola Sch\"onlieb and Daniel Spector}
\address[Jan Lellmann, Konstantinos Papafitsoros, Carola Sch\"onlieb]{Department of Applied Mathematics and Theoretical Physics,  University of Cambridge, United Kingdom}
\address[Daniel Spector]{Technion, Israel Institute of Technology, Haifa, Israel and Department of Applied Mathematics, National Chiao Tung University, Hsinchu, Taiwan}

\email{J.Lellmann@damtp.cam.ac.uk}
\email{k.papafitsoros@maths.cam.ac.uk}
\email{C.B.Schoenlieb@damtp.cam.ac.uk}
 \email{dspector@tx.technion.ac.il}

\begin{document}

\maketitle

\begin{abstract}
In this work we introduce a formulation for a non-local Hessian that combines the ideas of \emph{higher-order} and \emph{non-local} regularization for image restoration, extending the idea of non-local gradients to higher-order derivatives. By carefully choosing the weights, the model allows to improve on the current state of the art higher-order method, Total Generalized Variation, with respect to overall quality and particularly close to jumps in the data. In the spirit of recent work by Brezis et al., our formulation also has analytic implications:\ for a suitable choice of weights, it can be shown to converge to classical second-order regularizers, and in fact allows a novel characterization of higher-order Sobolev and BV spaces.

\end{abstract}

\section{Introduction}

The total variation model of image restoration due to Rudin, Osher and Fatemi \cite{rudin1992nonlinear} is now classical, the problem of being given a noisy image $g \in L^2(\Omega)$ on an open set $\Omega \subseteq \mathbb{R}^2$ and selecting a restored image via minimization of the energy
\begin{align*}
E(u):= \int_\Omega (u-g)^2\;dx+\alpha\TV(u).
\end{align*}
Here, $\alpha>0$ is a regularization parameter at our disposal and $\TV(u):=|Du|(\Omega)$ is the total variation of the measure that represents the distributional derivative of a function $u \in \BV(\Omega)$, the space of functions of bounded variation.  Among the known defects of the model is the staircasing effect, stemming from the use of the $\TV$ term in regularization.  It is natural, then, to investigate the replacement of the total variation with another regularizer, for instance a higher-order term (see 
\cite{Sch98a,LLT03,LT06,HS06,CEP07, mineJMIV} for the bounded Hessian framework, \cite{ChambolleLions,TGV,setzer2011infimal} for infimal convolution and generalizations, and \cite{Lellmann2013a} for anisotropic variants) or a non-local one (see, for example, the work of Buades, Coll and Morel \cite{buades2005review}, Kinderman, Osher and Jones \cite{kindermanosherjones2005}, and Gilboa and Osher \cite{gilboa2008nonlocal}).  In this work, we introduce and analyze a regularizer that is both higher-order and non-local, and utilize it in a model for image restoration.  The heuristic interpretation of our model is based on a non-local Hessian, and our use of this nomenclature is justified through a rigorous localization analysis.  Moreover, we give some numerical experiments to demonstrate that using our non-local Hessian, and by a suitable choice of weights, it is possible to derive specialized models that compete with current state of the art higher-order methods such as the Total Generalized Variation \cite{TGV}.

\subsection*{Background on non-local gradients} 
In the first order setting, our analysis of non-local gradients and their associated energies find its origins in the 2001 paper of Bourgain, Brezis and Mironescu \cite{bourgain2001another}. In their paper, they had introduced energies of the form
\begin{equation}\label{BBM_functional}
F_{n}u:=\int_{\Omega}\left(\int_{\Omega}\frac{|u(x)-u(y)|^{pq}}{|x-y|^{pq}}\rho_{n}(x-y)dx\right)^\frac{1}{q}dy,
\end{equation}
where $\Omega$ is a smooth bounded domain in $\mathbb{R}^{N}$ and $1\le p<\infty$, and in the special case $q=1$.  Here, the functions $\rho_{n}$
are radial mollifiers that are assumed to satisfy the following three properties for all $n\in \mathbb{N}$:
\begin{align}
\rho_{n}(x)&\ge 0,\label{rho_property_1}\\
\int_{\mathbb{R}^{N}}\rho_{n}(x)dx&=1, \label{rho_property_2}\\
\lim_{n\to\infty} \int_{|x|>\gamma}\rho_{n}(x)dx&=0,\quad \forall \gamma>0. \label{rho_property_3}
\end{align} 
The work \cite{bourgain2001another} connects the finiteness of the limit as $n \to \infty$ of the functional \eqref{BBM_functional} with the inclusion of a function $u \in L^p(\Omega)$ in the Sobolev space $W^{1,p}(\Omega)$ if $p>1$ or BV$(\Omega)$ if $p=1$.  As in the beginning of the introduction, the space BV$(\Omega)$ refers to the space of functions of bounded variation, and it is no coincidence that the two papers \cite{bourgain2001another, rudin1992nonlinear} utilize this energy space.    Indeed, Gilboa and Osher \cite{gilboa2008nonlocal} in 2008 independently introduced an energy similar to \eqref{BBM_functional}, terming it a non-local total variation, while the connection of the two and introduction of the parameter $q$ is due to Leoni and Spector \cite{LeoniSpector}.  In particular, they show in \cite{LeoniSpectorCorrection} that for $p=1$ and $q=2$ the functionals \eqref{BBM_functional}  $\Gamma$-converge to a constant times the total variation. This result extends previous work by Ponce \cite{ponce2004new} in the case $q=1$ (see also the work of Aubert and Kornprobst \cite{aubertnonlocal} for an application of these results to image processing).

Gilboa and Osher \cite{gilboa2008nonlocal} in fact introduced two forms of non-local total variations, and for our purposes here it will be useful to consider the second.   This alternative involves introducing a non-local gradient operator, defined by
\begin{equation}\label{MS_functional}
G_{n}u(x):=N\int_{\Omega}\frac{u(x)-u(y)}{|x-y|}\frac{x-y}{|x-y|}\rho_{n}(x-y)dy,\quad x\in \Omega,
\end{equation}
for $u \in C^1_c(\Omega)$. Then, one defines the non-local total variation as the $L^1$ norm of \eqref{MS_functional}. The localization analysis of the non-local gradient \eqref{MS_functional} has been performed by Mengesha and Spector in \cite{Spector}, where a more general (and technical) distributional definition is utilized.  Their first observation is that the definition of the non-local gradient via the Lebesgue integral \eqref{MS_functional} extends to spaces of weakly differentiable functions. In this regime they discuss the localization of \eqref{MS_functional}. They prove that the non-local gradient converges to its local
 analogue $D u$ in a topology that corresponds to the regularity of the underlying function $u$.  As a result, they obtain yet another characterization of the spaces $W^{1,p}(\Omega)$ and $\mathrm{BV}(\Omega)$.  Of notable interest for image processing purposes is their result on the $\Gamma$-convergence of the corresponding non-local total variation energies defined via non-local gradients of the form \eqref{MS_functional} to the local total variation.

\subsection*{Background on higher-order regularization} As the total variation and many more contributions in the image processing community have proven, a non-smooth regularization procedure indeed results in a nonlinear smoothing of the image, smoothing more in homogeneous areas of the image domain and preserving characteristic structures such as edges. In particular, the total variation regularizer is tuned towards the preservation of edges and performs very well if the reconstructed image is piecewise constant. The drawback of such a regularization procedure becomes apparent as soon as images or signals (in 1D) are considered which do not only consist of flat regions and jumps, but also possess slanted regions, i.e., piecewise affine parts. The artifact introduced by total variation regularization in this case is called staircasing \cite{ring2000structural}.

In \cite{ChambolleLions} Chambolle and Lions propose a higher-order method by means of an infimal convolution of two convex regularizers. Here, a noisy image is decomposed into three parts $g=u_1+u_2+n$ by solving
\begin{equation}\label{infconv}
\min_{(u_1,u_2)}  \frac{1}{2} \int_{\Omega}(u_1+u_2-g)^{2}dx +\alpha\, \mathrm{TV}(u_{1}) + \beta \,\mathrm{TV}^{2}(u_{2}),
\end{equation}
where $\mathrm{TV}^{2}(u_{2}):=D^2 u_2$ is the distributional Hessian of $u_2$. Then, $u_1$ and $u_{2}$ are the piecewise constant and  piecewise affine parts of $g$ respectively and $n$ the noise (or texture). For recent modifications of this approach in the discrete setting see also \cite{SS08,setzer2011infimal}. Other approaches to combine first and second regularization originate, for instance, from Chan, Marquina, and Mulet \cite{chan2001high} who consider total variation minimisation together with weighted versions of the Laplacian, the Euler-elastica functional \cite{masnou1998level,chan2002euler} which combines total variation regularization with curvature penalisation, and many more \cite{LT06,lairidge,mineJMIV,ipol.2013.40}. Recently Bredies et al. have proposed another interesting higher-order total variation model called Total Generalized Variation (TGV) \cite{TGV}.
The TGV regularizer of order $k$ is of the form
\begin{equation}\label{tgv}
\mathrm{TGV}_{\boldsymbol{\alpha}}^{k}(u)=\sup\left\{\int_\Omega u\mathrm{div}^k\xi\,dx: \xi\in C_c^k(\Omega, \mathrm{Sym}^k(\mathbb{R}^N)),\;
\|\mathrm{div}^l \xi\|_\infty\leq\alpha_l,\; l=0,\ldots,k-1\right\},
\end{equation}
where $\mathrm{Sym}^k(\mathbb{R}^N)$ denotes the space of symmetric tensors of order $k$ with arguments in $\mathbb{R}^N$, and $\boldsymbol{\alpha}=\{\alpha_l\}_{l=0}^{k-1}$ are fixed positive parameters. Its formulation for the solution of general inverse problems was given in \cite{BredValk}.

The idea of pure bounded Hessian regularization is considered by Lysaker et al. \cite{LLT03},  Scherzer et al. \cite{Sch98a,HS06}, Lefkimmiatis et al. \cite{lefkimmiatis2010hessian} and Bergounioux and Piffet \cite{Piffet}. In these works the considered model has the general form
$$
\min_u \frac{1}{2}\int_\Omega (u-g)^2\;dx+\alpha |D^2 u|(\Omega). 
$$
In  \cite{CEP07}, Chan et al. use the squared $L^{2}$ norm of the Laplacian as a regularizer also in combination with the $H^{-1}$ norm in the data fitting term.
Further, in \cite{PS08} minimizers of functionals which are regularized by the total variation of the $(l-1)$-th derivative, i.e.,
$$
|D\nabla^{l-1}u|(\Omega),
$$
are studied. Properties of such regularizers in terms of diffusion filters are further studied in \cite{DiWeBu09}. Therein, the authors consider the Euler-Lagrange equations corresponding to minimizers of functionals of the general type
$$
\mathcal{J}(u) = \int_\Omega (u-g)^2  dx + \alpha\int_\Omega f\left(\sum_{|k|=p} |D^k u|^2\right) dx,
$$
for different non-quadratic functions $f$. There are also works on higher-order PDE methods for image regularization, e.g. \cite{chan2001nontexture,LLT03,BG04,BEG,TVH1_1}. 

Confirmed by all of these works on higher-order total variation regularization the introduction of higher-order derivatives can have a positive effect on artifacts like the staircasing effect inherent to total variation regularization. In \cite{dal2009higher}, for instance, the authors analyse the positive effect of the Chan, Marquina, and Mulet model on the staircasing artifact.
 

\subsection*{Higher-order non-local regularization} From this perspective, the development of a non-local higher-order method promises to be interesting. One approach to a non-local higher-order method would be to approximate the gradients using non-local gradients as developed in \cite{Spector}.  For example, one can define a form of non-local Hessian as
\begin{equation}\label{MS_functional_vector}
G_{n}(\nabla u)(x)=N\int_{\Omega}\frac{\nabla u (x)-\nabla u(y)}{|x-y|}\otimes \frac{x-y}{|x-y|}\rho_{n}(x-y)dy,
\end{equation}
and obtain some straightforward characterization of $W^{2,p}(\mathbb{R}^{N})$ and $\mathrm{BV}^{2}(\mathbb{R}^{N})$, the latter being  the space of $W^{1,1}(\mathbb{R}^{N})$ functions whose second order distributional derivative is a bounded Radon measure. Here $\otimes$ denotes the standard tensor multiplication of vectors.  However, we find it advantageous to introduce and utilize a non-local Hessian that is derivative free, thus only requiring a function $u$ to belong to an $L^{p}$ space (instead of $W^{1,p}$).  In analogy with the theory of distributional derivatives, we first need a notion of non-local Hessian which can be applied to smooth functions.  This definition, via a Lebesgue integral, is as follows.\\

\newtheorem{definition1}{Definition}[section]
\begin{definition1}\label{NlocHessian_smooth_def}
Suppose $u \in C^2_c(\mathbb{R}^N)$.  Then we define the non-local Hessian as the Lebesgue integral
\begin{equation}\label{OUR_functional}
H_{n}u(x):=\frac{N(N+2)}{2}\int_{\mathbb{R}^{N}}\frac{u(x+h)-2u(x)+u(x-h)}{|h|^{2}}\frac{\left( h\otimes h-\frac{|h|^{2}}{N+2}I_{N}\right )}{|h|^{2}}\rho_{n}(h)dh,\quad x\in \mathbb{R}^{N},
\end{equation}
where here $I_{N}$ is the $N\times N$ identity matrix and $\rho_{n}$ is a sequence satisfying \eqref{rho_property_1}--\eqref{rho_property_3}. 
\end{definition1}

We can then define the distributional non-local Hessian via its action on test functions. \\
 
 \newtheorem{definition2}[definition1]{Definition}
\begin{definition2}\label{NlocHessian_distr_def}
Suppose $u \in L^p(\mathbb{R}^N)$.  Then we define the distributional non-local Hessian componentwise as
\begin{align}
<\mathfrak{H}^{ij}_n u,\varphi>:= \int_{\mathbb{R}^N} u H^{ij}_n \varphi\;dx, \label{distributionalNLH}
\end{align}
for $\varphi \in C^\infty_c(\mathbb{R}^N)$.
\end{definition2}

A natural question is then whether these two notions agree.  The following theorem shows that this is the case, provided the Lebesgue integral exists.\\

 \newtheorem{theorem3}[definition1]{Theorem}
\begin{theorem3}[Non-local Integration by Parts]\label{NIByparts}
Suppose that $u \in L^{p}(\mathbb{R}^N)$ for some $1\leq p <+\infty$ and  $\frac{|u(x+h)-2u(x)+u(x-h)|}{|h|^2}\rho_n(h)\in L^{1}(\mathbb{R}^N \times \mathbb{R}^N)$.
Then $\mathfrak{H}_nu \in L^{1}(\mathbb{R}^N,\mathbb{R}^{N \times N}) $ and for any $\varphi \in C_{c}^{2}(\mathbb{R}^N)$ and $i,j=1\ldots N$,
 \begin{align}
<\mathfrak{H}^{ij}_n u,\varphi> = \int_{\mathbb{R}^N}H_n^{ij}u(x) \varphi(x) \;dx \label{intbyparts},
\end{align}
and therefore $H_nu=\mathfrak{H}_nu$ almost everywhere as functions.
\end{theorem3}

We will see in Section \ref{analysis}, Lemmas \ref{W2p_Hessian_label} and \ref{BV2_Hessian_label} that the Lebesgue integral even makes sense for $u \in W^{2,p}(\mathbb{R}^N)$ or BV$^2(\mathbb{R}^N)$, and therefore the distributional definition $\mathfrak{H}_nu$ coincides with the Lebesgue integral for these functions.  

In the following, we undertake the rigorous analysis of the non-local Hessian, proving localization results in various topologies and characterizations of higher-order spaces of weakly differentiable functions.  Our first result is the following theorem concerning the localization in the smooth case.\\

 \newtheorem{theorem4}[definition1]{Theorem}
\begin{theorem4}\label{localization_smooth_label}
Suppose that $u\in C_{c}^{2}(\mathbb{R}^{N})$. Then for any $1\leq p \leq+\infty$, 
\[H_{n}u\to \nabla^{2}u,\quad  \text{in }L^{p}(\mathbb{R}^{N},\mathbb{R}^{N\times N})\quad \text{as }n\to\infty.\]
\end{theorem4}

When less smoothness is assumed on $u$, we have analogous convergence theorems, where the topology of convergence depends on the smoothness of $u$.  When $u \in W^{2,p}(\mathbb{R}^N)$, we have the following.\\

 \newtheorem{theorem5}[definition1]{Theorem}
\begin{theorem5}\label{loc_W2p_label}
Let $1\le p< \infty$. Then for every $u\in W^{2,p}(\mathbb{R}^{N})$ we have that
\[H_{n}u\to \nabla^{2}u\quad \text{in }L^{p}(\mathbb{R}^{N},\mathbb{R}^{N\times N})\quad \text{as }n\to\infty.\]
\end{theorem5}
In the setting of BV$^2(\mathbb{R}^N)$ (see Section \ref{preliminaries} for a definition), we have the following theorem on the localization of the non-local Hessian.\\

 \newtheorem{theorem6}[definition1]{Theorem}
\begin{theorem6}\label{loc_BV2_label}
Let $u\in \mathrm{BV}^{2}(\mathbb{R}^{N})$. Then
\[\mu_{n}\to D^{2}u,\quad \text{weakly}^{\ast},\]
i.e., for every $\phi\in C_{0}(\mathbb{R}^{N},\mathbb{R}^{N\times N})$
\begin{equation}\label{loc_weakstar}
\lim_{n\to\infty}\int_{\mathbb{R}^{N}}H_{n}u(x)\cdot \phi(x)dx=\int_{\mathbb{R}^{N}}\phi(x)\cdot \,d D^{2}u. 
\end{equation}
\end{theorem6}

We have seen that the non-local Hessian is well-defined as a Lebesgue integral and localizes for spaces of weakly differentiable functions.  In fact, it is sufficient to assume that $u \in L^p(\mathbb{R}^N)$ is a function such that the distributions $\mathfrak{H}_nu$ are in $L^p(\mathbb{R}^N,\mathbb{R}^{N \times N})$ with a uniform bound of their $L^{p}$ norms, in order to deduce that $u \in W^{2,p}(\mathbb{R}^N)$ if $1<p<+\infty$ or $u \in \mathrm{BV}^2(\mathbb{R}^N)$ if $p=1$.  Precisely, we have the following theorems characterizing the second order Sobolev and BV spaces.\\ 

 \newtheorem{theorem7}[definition1]{Theorem}
\begin{theorem7}\label{char}
Let $u\in L^{p}(\mathbb{R}^{N})$ for some $1<p<\infty$. Then
\begin{equation}\label{W2p_char}
u\in W^{2,p}(\mathbb{R}^{N}) \quad \iff \quad \liminf_{n\to\infty} \int_{\mathbb{R}^{N}}|\mathfrak{H}_{n}u(x)|^{p}dx<\infty.
\end{equation}
Let now $u\in L^{1}(\mathbb{R}^{N})$. Then
\begin{equation}\label{BV2_char}
u\in \mathrm{BV}^{2}(\mathbb{R}^{N}) \quad \iff \quad \liminf_{n\to\infty} \int_{\mathbb{R}^{N}}|\mathfrak{H}_{n}u(x)|dx<\infty.
\end{equation}
\end{theorem7}


While it is a natural extension of the first-order model and very well suited for a rigorous analysis, the definition of the non-local Hessian in \eqref{OUR_functional} has some drawbacks in practice. In order to utilize the full potential of non-local models in image processing practice, it is natural to allow arbitrary, fixed (non-localizing) weights $\sigma_x(h)$. The difficulty of accommodating this in our current model comes from the fact that \eqref{OUR_functional} has a built-in symmetry that associates triples of points rather than pairs which makes dealing with bounded domains more difficult and reduces the freedom in choosing the weights. Therefore we also propose an alternative formulation of the non-local Hessian that is closely related to \eqref{OUR_functional} but allows more freedom and is more amenable to numerical implementation. It is based on the observation that ultimately all finite-differences discretizations of the second-order derivatives stem from the requirement that they should be compatible with the Taylor expansion for smooth functions. More precisely, we can directly define the \emph{implicit} non-local gradient $G_u(x) \in \R^N$ and Hessian $H_u(x) \in \tmop{Sym}(\RNN)$ that \emph{best explain $u$ around $x$ in terms of a quadratic model}:
\begin{equation}\label{eq:implicitnl-intro}
(G'_u(x),H'_u(x)) =\hspace{-0.1cm}\underset{G_u\in\RN,H_u \in \tmop{Sym}(\RNN)}{\operatorname{argmin}}
 \hspace{0.1cm}\half  \int_{\Om - \{x\}} \left( u ( x + z ) -u ( x ) -G_{u} ^{\top}  z - \frac{1}{2}  z^{\top} H_{u} z
  \right)^{2} \sigma_{x} (z)  d z,
\end{equation}
where $\Omega \setminus \{x\}=\{y-x:\; y\in \Omega\}$.
While more complicated at first, this definition is equivalent to the explicit model \eqref{OUR_functional} under a certain choice of the weights, see Theorem \ref{thm:implexplcircle}.  As the objectives of the minimization problems \eqref{eq:implicitnl-intro} are quadratic, their solutions can be characterized by linear optimality conditions. Thus functionals based on the implicit non-local derivatives can be easily included in usual convex solvers by adding these conditions.
Moreover, the weights $\sigma_x(z)$ between any pair of points $x$ and $y=x+z$ can be chosen arbitrarily, without any restrictions on symmetry. These advantages come at the cost of being considerably harder to treat analytically, and we leave a rigorous definition and analysis of \eqref{eq:implicitnl-intro} for future work.

Of practical relevance is the question of how to choose the weights $\sigma_x$ for a particular purpose. We show an example where the weights are used to construct a regularizer that both favours piecewise affine functions but also allows for jumps in the data. It is motivated by the recent discussion of ``Amoeba'' filters in \cite{Lerallut2007,Welk2011,Welk2012} which combine standard filters such as median filters with non-parametric structuring elements that are based on the data -- that is in long thin objects they would extend along the structure and thus prevent smoothing perpendicular to the structure. In Amoeba filtering, the shape of a structuring element at a point is defined as a unit circle with respect to the geodesic distance on a manifold defined by the image itself. Motivated by their results we propose to extend their method to variational image processing by using the geodesic distance to set the weights $\sigma_x$. This allows to get a very close approximation to true piecewise affine regularization, in many cases improving on the results obtained using TGV.

Recently there has been another approach at defining a higher-order extension of non-local regularization for TGV-based energies \cite{Ranftl2014}. The authors start with the cascading formulation of TGV,
\[
\mathrm{TGV}^2(u) = \inf_{w:\Om\to\RN} \alpha_1 \int_\Omega |D u - w| + \alpha_0 \int_\Omega |D w |,
\]
which reduces the higher-order differential operators that appear in the definition of TGV to an infimal convolution of two terms involving only first-order derivatives. These can then be replaced by classical first-order non-local derivatives, and one obtains an energy of the form
\[
\inf_{w:\Om\to\RN} \int_\Om \int_\Om \alpha_1(x,y) |u(x)-u(y)- w(x)^\top(x-y)| dy dx + \sum_{i=1}^2 \int_\Om \int_\Om \alpha_0(x,y)|w^i(x)-w^i(y)| d y d x.
\]
Comparing this to \eqref{eq:implicitnl-intro}, we see that the left integral is similar to our implicit model for first-order derivatives with $w$ taking the role of $G_u$, but with the square in \eqref{eq:implicitnl-intro} replaced by a $L^1$-type norm and added to the energy as a penalty term rather than strictly enforcing optimality.

Finally, let us mention an important localization result from the perspective of variational image processing, the following theorem asserting the $\Gamma$-convergence \cite{dalmasogamma, Braidesgamma} of the non-local Hessian energies to the energy of the Hessian.\\
  
 \newtheorem{theorem8}[definition1]{Theorem}
\begin{theorem8}\label{thm:explicitgammaconv}
Assume that $\rho_{n}$ satisfy \eqref{rho_property_1}, \eqref{rho_property_2} and \eqref{rho_property_3}.  Then
\[
\Gamma^-_{L^1(\mathbb{R}^N)} \mhyphen \lim_{n \to \infty} \int_{\mathbb{R}^N} |\mathfrak{H}_n u|\;dx = |D^{2}u|(\mathbb{R}^N),
\]
where the $\Gamma \mhyphen$limit is taken with respect to the strong convergence $u_n \to u$ in $L^1(\mathbb{R}^N)$.
\end{theorem8}
The relevance of this theorem in the context of variational problems comes from the fact that $\Gamma$-convergence of the \emph{objective functions} of a sequence of minimization problems implies convergence of the minimizers in a suitable topology \cite[Ch.~1.5]{Braidesgamma}. Therefore Theorem~\ref{thm:explicitgammaconv} guarantees that under a suitable choice of weights, the solutions of a class of \emph{non-local} Hessian-based problems converges to the solution of the \emph{local} Hessian-regularized problem, and thus our notion of ``non-local Hessian'' is justified.

\subsection*{Organization of the paper}
The paper is organized as follows: 
In Section \ref{preliminaries} we recall some preliminary notions and we fix our notation.
Section \ref{analysis} deals with the analysis of the non-local Hessian functional \eqref{OUR_functional}. After a justification of the  introduction of its distributional form, we proceed in Section \ref{localization_smooth} with the localization of \eqref{OUR_functional} to the classical Hessian for smooth functions $u$. The localization of \eqref{OUR_functional} to its classical analogue for $W^{2,p}(\mathbb{R}^{N})$ and BV$^{2}(\mathbb{R}^{N})$ functions is shown in Sections \ref{localisation_W2p} and \ref{localization_BV2} respectively. In Section \ref{characterizations} we provide the non-local characterizations of the spaces $W^{2,p}(\mathbb{R}^{N})$ and BV$^{2}(\mathbb{R}^{N})$ in the spirit of \cite{bourgain2001another}. The $\Gamma$-convergence result, Theorem \ref{thm:explicitgammaconv}, is proved in Section \ref{section:gamma}.

The introduction of the implicit formulation of non-local Hessian \eqref{eq:implicitnl-intro},  along with its connection to the explicit one, is presented in Section \ref{sec:implicitdef}. In Section \ref{sec:jumppreservingweights} we describe how we choose the weights $\sigma_{x}$ in \eqref{eq:implicitnl-intro} in order to achieve jump preservation in the restored images. Finally, in Section \ref{sec:experiments} we present our numerical results, comparing our method with TGV.

\section{Preliminaries and Notation}\label{preliminaries}
For the reader's convenience we recall here some important notions that we are going to use in the following sections and we also fix some notation.

As far as our notation is concerned, whenever a function space has two arguments, the first always denotes the  function's domain while the second denotes its range. Whenever the range is omitted, it is assumed that the functions are real valued.

We use $dx,dy,dz$  for various integrations with respect Lebesgue measure on $\mathbb{R}^N$, while in Section \ref{analysis} we will ave occasion to use the more succinct notation $\mathcal{L}^{N^2}$ to denote integration with respect to the Lebesgue measure in the product space $\mathbb{R}^N \times \mathbb{R}^N$.

 The reader should not confuse the different forms of the letter ``H''. We denote by $\mathcal{H}$ the one-dimensional Hausdorff measure ($\mathcal{H}^{N}$ for the N-dimensional), while $H$ denotes the non-local Hessian when this is a function. As we have already seen, $\mathfrak{H}$ denotes the distributional form of the non-local Hessian.
 
  It is also very convenient to introduce the following notation:
\[\mathpzc{d}^{2}u(x,y):=u(y)-2u(x)+u(x+(x-y)),\]
which can be interpreted a discrete second order differential operator in $x$ at the direction $x-y$. 

We denote by $|\cdot|$ the Euclidean norm (vectors) and Frobenius norm (matrices). 

As usual, we denote by BV$(\Omega)$ the space of functions of bounded variation defined on an open $\Omega\subseteq \mathbb{R}^{N}$. This space consists of all real valued functions $u\in L^{1}(\Omega)$ whose distributional derivative $Du$ can be represented by a finite Radon measure. The total variation TV$(u)$ of a function $u\in\mathrm{BV}(\Omega)$, is defined to be the total variation of the measure $Du$, i.e., TV$(u):=|Du|(\Omega)$. The definition is similar for vector valued functions.
We refer the reader to \cite{AmbrosioBV} for a full account of the theory of BV functions.

We denote by BV$^{2}(\Omega)$ the space of functions of bounded Hessian. These are all the functions that belong to the Sobolev space $W^{1,1}(\Omega)$ such that $\nabla u$ is an $\mathbb{R}^{N}$-valued BV function, i.e., $\nabla u\in \mathrm{BV}(\Omega,\mathbb{R}^{N})$ and we set $D^{2}u:=D(\nabla u)$. We refer the reader to \cite{demengelBH, Piffet, mineJMIV} for more information about this space. Let us however state a theorem that will be useful for our purposes.  It is the analogue result to the strict approximation by smooth functions for the classical $\mathrm{BV}$ case, see \cite{AmbrosioBV}.\\

 \newtheorem{theorem9}[definition1]{Theorem}
\begin{theorem9}[$\mathrm{BV}^{2}$ strict approximation by smooth functions, \cite{demengelBH}]\label{strictBV2_label}
Let $\Omega\subseteq \mathbb{R}^{N}$ open and $u\in\mathrm{BV}^{2}(\Omega)$. Then there exists a sequence $(u_{n})_{n\in\mathbb{N}}\in W^{2,1}(\Omega)\cap C^{\infty}(\Omega)$ that converges to $u$ strictly in $\mathrm{BV}^{2}(\Omega)$, that is to say
\[u_{n}\to u,\quad\text{in }L^{1}(\Omega)\quad \text{and}\quad |D^{2}u_{n}|(\Omega)\to |D^{2}u|(\Omega),\quad \text{as }n\to\infty.\]
\end{theorem9}

We recall also the two basic notions of convergence regarding finite Radon measures. We note that $\mathcal{M}(\Omega,\mathbb{R}^{\ell})$ denotes the space of $\mathbb{R}^{\ell}$-valued finite Radon measures in $\Omega$.
If $(\mu_{n})_{n\in\mathbb{N}}$ and $\mu$ are real valued finite Radon measures defined on an open $\Omega\subseteq \mathbb{R}^{N}$ we say that the sequence $\mu_{n}$ converges weakly$^{\ast}$ to $\mu$ if for all $\phi\in C_{0}(\Omega)$ we have $\int_{\Omega}\phi\,d\mu_{n}\to \int_{\Omega}\phi\,d\mu$ as $n$ goes to infinity. Here $\phi \in C_{0}(\Omega)$ means that $\phi$ is continuous on $\Omega$ and that for every $\epsilon>0$ there exists a compact set $K \subset \Omega$ such that $\sup_{x \in \Omega \setminus K} |\phi(x)| \leq \epsilon$. Note that from the Riesz representation theorem the dual space $(C_{0}(\Omega,\mathbb{R}^{\ell}), \|\cdot\|_{\infty})^{\ast}$ can be identified with $\mathcal{M}(\Omega,\mathbb{R}^{\ell})$.
 We say that the convergence is strict if in addition to that we also have that $|\mu_{n}|(\Omega)\to |\mu|(\Omega)$, i.e., the total variations of $\mu_{n}$ converge to the total variation of $\mu$. The definition is similar for vector and matrix valued measures with all the operations regarded component-wise.

We now remind the reader about some basic facts concerning $\Gamma$-convergence. Let $(X,d)$ be a metric space and $F, F_{n}:X\to\mathbb{R}\cup\{+\infty\}$ for all $n\in \mathbb{N}$. We say that the sequence of functionals $F_{n}$ $\Gamma$-converges to $F$ at $x\in X$ in the topology of $X$ and we write $\Gamma_{X}^{-}\mhyphen \lim_{n\to\infty} F_{n}(x)=F(x)$ if the following two conditions hold\\
\begin{enumerate}
\item For every sequence $(x_{n})_{n\in\mathbb{N}}$ converging to $x$ in $(X,d)$ we have
\[F(x)\le \liminf_{n\to\infty}F_{n}(x_{n}).\]
\item There exists a sequence $(x_{n})_{n\in\mathbb{N}}$ converging to $x$ in $(X,d)$ such that
\[F(x)\ge \limsup_{n\to\infty}F_{n}(x_{n}).\]
\end{enumerate}
It can be proved that $\Gamma_{X}^{-} \mhyphen\lim_{n\to\infty} F_{n}(x)=F(x)$ if the $\Gamma$-lower and $\Gamma$-upper limits of $F_{n}$ at $x$, denoted by $\Gamma_{X}^{-}\mhyphen\liminf_{n\to\infty} F_{n}(x)$ and $\Gamma_{X}^{-}\mhyphen\limsup_{n\to\infty} F_{n}(x)$ respectively, are equal to $F(x)$ where
\begin{align*}
\Gamma_{X}^{-}\mhyphen\liminf_{n\to\infty} F_{n}(x)& = \min \{\liminf_{n\to\infty} F_{n}(x_{n}):\; x_{n}\to x \text{ in } (X,d)\},\\
\Gamma_{X}^{-}\mhyphen\limsup_{n\to\infty} F_{n}(x)&=\min \{\limsup_{n\to\infty} F_{n}(x_{n}):\; x_{n}\to x \text{ in } (X,d)\}.
\end{align*}
Finally, if $F:X\to\mathbb{R}\cup\{+\infty\}$ we denote by $sc_{X}^{-} F$ the lower semicontinuous envelope of  $F$, i.e., the greatest lower semicontinuous function majorized by $F$.
We refer the reader to \cite{dalmasogamma, Braidesgamma} for further details regarding $\Gamma$-convergence and lower semicontinuous envelopes.

\section{Analysis of the non-local Hessian}\label{analysis} 
The precise form we have chosen for the non-local Hessian can be derived from the model case of non-local gradients - the fractional gradient - which has been developed in \cite{SpectorShieh}. Here we prove several results analogous to the first order case, as in \cite{Spector}, for the generalizations involving generic radial weights that satisfy \eqref{rho_property_1}--\eqref{rho_property_3}.  Of primary importance is to first establish that the distributional non-local Hessian defined by equation \eqref{distributionalNLH} is, in fact, a distribution.
 Here we observe that if $u \in L^1(\mathbb{R}^N)$, then
\begin{align*}
|<\mathfrak{H}_n u,\varphi>| \leq C \|u\|_{L^1(\mathbb{R}^N)} \|\nabla^2 \varphi\|_{L^\infty(\mathbb{R}^N)},
\end{align*}
so that $\mathfrak{H}_n u$ is a distribution. Also observe that if $u\in L^p(\mathbb{R}^N)$ for some $1<p<\infty$, then from the estimate \eqref{estimateW2p} below together with the fact that $\varphi$ is of compact support we have
\[|<\mathfrak{H}_n u,\varphi>| \leq C \|u\|_{L^{p}(\mathbb{R}^{N})} \|\nabla^2 \varphi\|_{L^q(\mathbb{R}^N,\mathbb{R}^{N\times N})}\le C \|u\|_{L^p(\mathbb{R}^N)}\|\nabla^2 \varphi\|_{L^\infty(\mathbb{R}^N)},\]
where $1/p+1/q=1$ and thus $\mathfrak{H}_n u$ is indeed again a distribution.
 One observes that the definition is in analogy to the theory of Sobolev spaces, where weak derivatives are defined in terms of the integration by parts formula. Because the Hessian is composed of two derivatives, we observe that there is no change in sign in the definition, preserving some symmetry that will be useful for us in the sequel.  
 
 The second important item to address is the agreement of the distributional non-local Hessian with the non-local Hessian.  The necessary and sufficient condition is the existence of the latter, which is the assertion of Theorem \ref{NIByparts}.  We now substantiate this assertion.\\

{\em Proof of Theorem \ref{NIByparts}.}
Let $1\leq p <+\infty$ and suppose $u \in L^p(\mathbb{R}^N)$ and that $H_n u$ is well-defined as a Lebesgue integral.  Let $\varphi \in C^2_c(\mathbb{R}^N)$, and fix $i,j \in \{1, \ldots, N\}$.  Then it is a consequence of Fubini's theorem and Lebesgue's dominated convergence theorem that
\begin{align*}
&\int_{\mathbb{R}^{N}}H^{ij}_{n} u(x)\varphi(x)dx\\
&=\frac{N(N+2)}{2}\lim_{\epsilon\to 0}\int_{\mathbb{R}^{N}}\int_{\mathbb{R}^{N}\setminus B(x,\epsilon)}\frac{\mathpzc{d}^{2}u(x,y)}{|x-y|^{2}}\frac{\left ((x_i-y_i)\otimes (x_j-y_j)-\frac{|x-y|^{2}}{N+2}I_{N}\right)}{|x-y|^{2}}\rho(x-y) \varphi(x)dydx\\
&=\frac{N(N+2)}{2}\lim_{\epsilon\to 0}\int_{d_{\epsilon}^{N}}\frac{\mathpzc{d}^{2}u(x,y)}{|x-y|^{2}}\frac{\left ((x_i-y_i)\otimes (x_j-y_j)-\frac{|x-y|^{2}}{N+2}I_{N}\right)}{|x-y|^{2}}\rho(x-y) \varphi(x)d(\mathcal{L}^{N})^{2}(x,y),
\end{align*}
where $d_{\epsilon}^{N}:=\mathbb{R}^{N}\times \mathbb{R}^{N}\setminus \{|x-y|<\epsilon\}$.
Similarly we have
\begin{align*}
&\int_{\mathbb{R}^{N}}u(x) H^{ij}_n \varphi (x)dx\\
&=\frac{N(N+2)}{2}\lim_{\epsilon\to 0}\int_{\mathbb{R}^{N}}\int_{\mathbb{R}^{N}\setminus B(x,\epsilon)}\hspace{-0.3cm}u(x)\frac{\mathpzc{d}^{2}\varphi(x,y)}{|x-y|^{2}}\frac{(x_{i}-y_{i})(x_{j}-y_{j})-\delta_{ij}\frac{|x-y|^{2}}{N+2}}{|x-y|^{2}}\rho(x-y)dydx\\
&=\frac{N(N+2)}{2}\lim_{\epsilon\to 0}\int_{d_{\epsilon}^{N}}u(x)\frac{\mathpzc{d}^{2}\varphi(x,y)}{|x-y|^{2}}\frac{(x_{i}-y_{i})(x_{j}-y_{j})-\delta_{ij}\frac{|x-y|^{2}}{N+2}}{|x-y|^{2}}\rho(x-y)d(\mathcal{L}^{N})^{2}(x,y),
\end{align*}
where, for notational convenience, we used the standard convention
\[
\delta_{ij}=
\begin{cases}
1& \text{if}\quad i=j,\\
0& \text{if}\quad i\ne j.
\end{cases}
\]
Thus, it suffices to show that for every $i,j$ and $\epsilon>0$ we have
\begin{eqnarray}\label{suffices_parts2}
&&\int_{d_{\epsilon}^{N}}\frac{\mathpzc{d}^{2}u(x,y)}{|x-y|^{2}}\frac{(x_{i}-y_{i}) (x_{j}-y_{j})-\delta_{ij}\frac{|x-y|^{2}}{N+2}}{|x-y|^{2}}\rho(x-y)\varphi(x)d(\mathcal{L}^{N})^{2}(x,y)=\\
&&\int_{d_{\epsilon}^{N}}u(x)\frac{\mathpzc{d}^{2}\varphi(x,y)}{|x-y|^{2}}\frac{(x_{i}-y_{i})(x_{j}-y_{j})-\delta_{ij}\frac{|x-y|^{2}}{N+2}}{|x-y|^{2}}\rho(x-y)d(\mathcal{L}^{N})^{2}(x,y).\nonumber
\end{eqnarray}
In order to show \eqref{suffices_parts2}, it suffices to prove 
\begin{eqnarray}\label{suffices_part2_1}
&&\int_{d_{\epsilon}^{N}}\frac{u(y)\varphi(x)}{|x-y|^{2}}\frac{(x_{i}-y_{i})(x_{j}-y_{j})-\delta_{ij}\frac{|x-y|^{2}}{N+2}}{|x-y|^{2}}\rho(x-y)d(\mathcal{L}^{N})^{2}(x,y)=\\
&&\int_{d_{\epsilon}^{N}}\frac{u(x)\varphi(y)}{|x-y|^{2}}\frac{(x_{i}-y_{i})(x_{j}-y_{j})-\delta_{ij}\frac{|x-y|^{2}}{N+2}}{|x-y|^{2}}\rho(x-y)d(\mathcal{L}^{N})^{2}(x,y), \nonumber
\end{eqnarray}
and 
\begin{eqnarray}\label{suffices_part2_2}
&&\int_{d_{\epsilon}^{N}}\frac{u(x+(x-y))\varphi(x)}{|x-y|^{2}}\frac{(x_{i}-y_{i})(x_{j}-y_{j})-\delta_{ij}\frac{|x-y|^{2}}{N+2}}{|x-y|^{2}}\rho(x-y)d(\mathcal{L}^{N})^{2}(x,y)=\\
&&\int_{d_{\epsilon}^{N}}\frac{u(x)\varphi(x+(x-y))}{|x-y|^{2}}\frac{(x_{i}-y_{i})(x_{j}-y_{j})-\delta_{ij}\frac{|x-y|^{2}}{N+2}}{|x-y|^{2}}\rho(x-y)d(\mathcal{L}^{N})^{2}(x,y). \nonumber
\end{eqnarray}
Equation \eqref{suffices_part2_1} can be easily showed by alternating $x$ and $y$ and using the symmetry of the domain. Finally equation \eqref{suffices_part2_2} can be proved by employing the substitution $u=2x-y$, $v=3x-2y$, noting that $x-y=v-u$ and that the determinant of the Jacobian of this substitution is $-1$. 
\endproof

 Having established that the notion of distributional non-local Hessian and non-local Hessian agree whenever the latter exists, it is a natural question to ask when this is the case.  It is a simple calculation to verify that the Lebesgue integral \eqref{OUR_functional} exists whenever $u \in C^2_c(\mathbb{R}^N)$.  The following several theorems, whose proofs we also defer until after the analysis in the smooth case, show that this is also the case for functions in the spaces $W^{2,p}(\mathbb{R}^N)$ and BV$^2(\mathbb{R}^N)$.\\

 \newtheorem{lemma10}[definition1]{Lemma}
\begin{lemma10}\label{W2p_Hessian_label}
Suppose that $u\in W^{2,p}(\mathbb{R}^{N})$, where $1\le p<\infty$. Then $H_nu$ is well-defined as a Lebesgue integral, $H_{n}u\in L^{p}(\mathbb{R}^{N},\mathbb{R}^{N\times N})$, and 
\begin{equation}\label{estimateW2p}
\int_{\mathbb{R}^{N}}|H_{n}u(x)|^{p}dx\le M \|\nabla^{2}u\|_{L^{p}(\mathbb{R}^{N})}^{p},
\end{equation}
where the constant $M$ depends only on $N$ and $p$.
\end{lemma10}

 \newtheorem{lemma11}[definition1]{Lemma}
\begin{lemma11}\label{BV2_Hessian_label}
Suppose that $u\in \mathrm{BV}^{2}(\mathbb{R}^{N})$. Then $H_{n}u\in L^{1}(\mathbb{R}^{N},\mathbb{R}^{N\times N})$
with
\begin{equation}\label{estimateBV2}
\int_{\mathbb{R}^{N}}|H_{n}u(x)|dx\le M|D^{2}u|(\mathbb{R}^{N}),
\end{equation}
where the constant $M$ depends only on $N$.
\end{lemma11}

We finally point out that all our initial analysis is done in $\mathbb{R}^{N}$. Defining \eqref{OUR_functional} in a general $\Omega\subseteq \mathbb{R}^{N}$ is problematic as the term $x+(x-y)$ does not necessarily belong to $\Omega$ whenever $x,y\in \Omega$. 
\subsection{Localization -- Smooth case} \label{localization_smooth}

We are now ready to prove the localization of $H_{n}u$ to $\nabla^{2}u$ for smooth functions.\\

{\em Proof of Theorem \ref{localization_smooth_label}.}

\textbf{Case} $1\leq p <+\infty$\\
Let us assume that we have shown the case $p=+\infty$.  Then we must show that the uniform convergence $H_n v \to \nabla^2v$ for $v \in C^2_c(\mathbb{R}^N)$, implies convergence in $L^p(\mathbb{R}^N,\mathbb{R}^{N\times N})$ for any $1\leq p<+\infty$.  We claim that this will follow from the following uniform estimate on the \emph{tails} of the nonlocal Hessian.  Suppose $supp\;v \subset B(0,R)$, where $supp\;v $ denotes the support of $v$.  Then for any $1\leq p<+\infty$ and $\epsilon>0$ there exists a $L = L(\epsilon,p)>1$ such that
\begin{align}
\sup_n \int_{B(0,LR)^c} |H_nv(x)|^p\;dx \leq \epsilon. \label{estimatetails}
\end{align}
If this were the case, we would estimate the $L^p$-convergence as follows
\begin{align*}
\int_{\mathbb{R}^N} |H_nv(x) - \nabla^2v(x)|^p\;dx = \int_{B(0,LR)} |H_nv(x) - \nabla^2v(x)|^p\;dx +  \int_{B(0,LR)^c} |H_nv(x) |^p\;dx,
\end{align*}
from which equation \eqref{estimatetails} implies
\begin{align*}
\limsup_{n\to \infty}\int_{\mathbb{R}^N} |H_nv(x) - \nabla^2v(x)|^p\;dx \leq \limsup_{n \to \infty} \int_{B(0,LR)} |H_nv(x) - \nabla^2v(x)|^p\;dx +   \epsilon.
\end{align*}
The conclusion then follows, since the first term vanishes from the uniform convergence assumed, after which $\epsilon>0$ is arbitrary.  We will therefore show the estimate \eqref{estimatetails}.   We have, by Jensen's inequality with respect to the measure $\rho_n$, which has $\int_{\mathbb{R}^N} \rho_n(x)dx = 1$, that
\begin{align*}
\int_{B(0,LR)^c} |H_nv(x)|^p\;dx &\leq \int_{B(0,LR)^c}\int_{\mathbb{R}^{N}}\frac{|v(y)-2v(x)+v(x+x-y)|^{p}}{|x-y|^{2p}}\rho_{n}(x-y)dydx \\
&= \int_{B(0,LR)^c}\int_{y \in B(0,R)}\frac{|v(y)|^{p}}{|x-y|^{2p}}\rho_{n}(x-y)dydx \\
&\;\;+ \int_{B(0,LR)^c}\int_{x+x-y \in B(0,R)}\frac{|v(x+x-y)|^{p}}{|x-y|^{2p}}\rho_{n}(x-y)dydx.
\end{align*}
Letting $z=x+x-y$ (which means that $x-y=z-x$), we obtain
\begin{align*}
\int_{B(0,LR)^c}\int_{x+x-y \in B(0,R)}\frac{|v(x+x-y)|^{p}}{|x-y|^{2p}}\rho_{n}(x-y)dydx =  \int_{B(0,LR)^c} \int_{z \in B(0,R)}\frac{|v(z)|^{p}}{|z-x|^{2p}}\rho_{n}(z-x)dzdx,
\end{align*}
 and therefore by symmetry of $\rho_n$ we have
 \begin{align*}
\int_{B(0,LR)^c} |H_nv(x)|^p\;dx &\leq 2\int_{B(0,LR)^c}\int_{y \in B(0,R)}\frac{|v(y)|^{p}}{|x-y|^{2p}}\rho_{n}(x-y)dydx \\
&\leq \frac{2}{(L-1)^{2p}} \int_{B(0,LR)^c} \int_{y \in B(0,R)}|v(y)|^{p}\rho_{n}(x-y)dydx \\
&\leq \frac{2}{(L-1)^{2p}} \|\rho_n\|_{L^1(\mathbb{R}^N)} \|v\|_{L^p(\mathbb{R}^N)}^p.
\end{align*}
Again using $\int_{\mathbb{R}^{N}} \rho_n(x)dx=1$, the claim, and therefore the case $1\leq p<+\infty$, then follows by choosing $L$ sufficiently large.

\textbf{Case} $p=+\infty$\\
It therefore remains to show that the convergence in $L^\infty(\mathbb{R}^N,\mathbb{R}^{N\times N})$ is true.  Precisely, we will show that
\[\left |H_{n}u - \nabla^{2}u \right | \to 0,\quad \text{uniformly},\]
for which it suffices to prove the convergence component wise, i.e., $\left |\left (H_{n}u - \nabla^{2}u \right )_{(i_{0},j_{0})}\right |\to 0$ by considering two cases $i_{0}\neq j_{0}$ and $i_{0}=j_{0}$.

\textbf{Subcase} $\mathbf{i_{0}\ne j_{0}}:$\\
Let us first observe that Proposition \ref{prop:cij} in the Appendix and the assumption that $\int_{\mathbb{R}^{N}} \rho_n(x)dx=1$ for all $n \in \mathbb{N}$ can be used to deduce that
\begin{align}
\int_{\mathbb{R}^N} \frac{z_{i_0}^{2}z_{j_{0}}^{2}}{|z|^{4}}\rho_{n}(z)dz &= \int_{0}^\infty \rho_n(t)t^{N-1}\;dt \int_{\Ss^{N-1}} \nu_{i_{0}}^{2} \nu_{j_{0}}^{2}  d
    \mathcal{H}^{N-1} ( x )  \label{oneoverthree}\\
    &=  \frac{1 }{N ( N+2 )} \cdot \left\{ \begin{array}{ll}
      1, & i_{0} \neq j_{0} ,\\
      3, & i_{0} =j_{0} .
    \end{array} \right.\nonumber
  \end{align}
Therefore, in the case that $i_0 \neq j_0$, we have
 \begin{align*}
\bigg\vert \big(H_{n}u& - \nabla^{2}u \big)_{(i_{0},j_{0})}\bigg\vert \\
&=\frac{N(N+2)}{2}\left |\int_{\mathbb{R}^{N}} \frac{\mathpzc{d}^{2}u(x,y)}{|x-y|^{2}}\frac{(x_{i_{0}}-y_{i_{0}})(x_{j_{0}}-y_{j_{0}})}{|x-y|^{2}}\rho_{n}(x-y) dy\right.-\left.2\frac{\partial u}{\partial x_{i_{0}}\partial x_{j_{0}}}(x) \int_{\mathbb{R}^{N}} \frac{z_{i_{0}}^{2}z_{j_{0}}^{2}}{|z|^{4}}\rho_{n}(z)dz\right |.
\end{align*}
However, again utilizing the radial symmetry of $\rho_n$, we have that the following integrals vanish:
\begin{align*}
\mathbf{i_{0}=j_{0}}&: \int_{\mathbb{R}^N}\frac{z_{i}z_{j_{0}}^{3}}{|z|^{4}}\rho_{n}(z)dz=0, \quad \text{for }i\ne j_{0},\\
\mathbf{i_{0}=j_{0}}&: \int_{\mathbb{R}^N}\frac{z_{i} z_{j}z_{j_{0}}^{2}}{|z|^{4}}\rho_{n}(z)dz=0, \quad \text{for }i\ne j_{0}, j \ne j_{0}, i\ne j,\\
\mathbf{i_{0}\ne j_{0}}&: \int_{\mathbb{R}^N}\frac{z_{i} z_{j}z_{i_{0}}z_{j_{0}}}{|z|^{4}}\rho_{n}(z)dz=0,\quad \text{for }i\ne i_{0}, j\ne j_{0},
\end{align*}
which implies that
\begin{align*}
 \sum_{i,j=1}^N\frac{\partial u}{\partial x_{i}\partial x_{j}}(x) \int_{\mathbb{R}^{N}} \frac{z_iz_{i_{0}}z_jz_{j_{0}}}{|z|^{4}}\rho_{n}(z)dz = 2\frac{\partial u}{\partial x_{i_{0}}\partial x_{j_{0}}}(x) \int_{\mathbb{R}^{N}} \frac{z_{i_{0}}^{2}z_{j_{0}}^{2}}{|z|^{4}}\rho_{n}(z)dz.
\end{align*}
Thus, introducing these factors of zero, and writing in a more compact way, we have that
 \begin{align*}
\bigg\vert \big(H_{n}u& - \nabla^{2}u \big)_{(i_{0},j_{0})}\bigg\vert \\
&=\frac{N(N+2)}{2} \left |\int_{\mathbb{R}^{N}}\frac{\mathpzc{d}^{2}u(x,y)-(x-y)^{T}\nabla^{2}u(x)(x-y)}{|x-y|^{2}}\frac{(x_{i_0}-y_{i_0})(x_{j_0}-y_{j_0})}{|x-y|^{2}}\rho_{n}(x-y)dy\right |.
\end{align*}
We want to show that the right hand side tends to zero as $n \to \infty$, and therefore we define now the following quantity for $\delta>0$:    
\begin{equation}\label{Q_delta}
Q_{\delta}u(x)=\left |\int_{B(x,\delta)}\frac{\mathpzc{d}^{2}u(x,y)-(x-y)^{T}\nabla^{2}u(x)(x-y)}{|x-y|^{2}}\frac{(x_{i_0}-y_{i_0})(x_{j_0}-y_{j_0})}{|x-y|^{2}}\rho_{n}(x-y)dy \right |.
\end{equation}
We then claim that we can make $Q_{\delta}u(x)$ as small as we want, independently of $x$ and $n$, by choosing sufficiently small $\delta>0$.  If this is the case, then the case $i_0 \neq j_0$ would be completed, since we would then have that
 \begin{eqnarray*}
\left |\left (H_{n}u -\nabla^{2}u \right )_{(i_{0},j_{0})}\right |&\leq& Q_\delta(x)+\frac{N(N+2)}{2}\int_{|z|\ge \delta} \frac{ |u(x+z)-2u(x)+u(x-z) |}{|z|^{2}}\rho_{n}(z) dz\\
&& +\frac{N(N+2)}{2}\left |\nabla^2u(x) \right |\int_{|z|\ge \delta} \rho_{n}(z)dz \\
&\le&\frac{N(N+2)}{2}\epsilon+\frac{N(N+2)}{2}\left(\frac{4\|u\|_{\infty}}{\delta^{2}}+\|\nabla^{2} u\|_{\infty}\right)\int_{|z|\ge \delta}\rho_{n}(z)dz\\
&<&N(N+2)\epsilon,
\end{eqnarray*}
 for $n$ large enough, and the result follows from sending $\epsilon \to 0$.

We therefore proceed to make estimates for \eqref{Q_delta}.  Since we have assumed $u \in C^2_c(\mathbb{R}^N)$, we have that given $\epsilon>0$, there is a $\delta>0$ such that for every $i,j=1,\ldots,N$ we have
\[\left|\frac{\partial u}{\partial x_{i}\partial {x_{j}}}(x)-\frac{\partial u}{\partial x_{i}\partial {x_{j}}}(y)\right |<\epsilon,\quad \text{whenever }|x-y|<\delta.\]
Using \eqref{mean_value_duxy} we can estimate
 \begin{align}
 Q_{\delta}u(x)&=\left |\int_{B(x,\delta)}\frac{\mathpzc{d}^{2}u(x,y)-(x-y)^{T}\nabla^{2}u(x)(x-y)}{|x-y|^{2}}\frac{(x_{i_0}-y_{i_0})(x_{j_0}-y_{j_0})}{|x-y|^{2}}\rho_{n}(x-y)dy \right | \label{uniform_cts}\\
                                &=\left |\int_{B(x,\delta)}\frac{(x-y)^{T}\left ( \int_{0}^{1}\int_{0}^{1}\nabla^{2}u(x+(s+t-1)(y-x))-\nabla^{2}u(x)dsdt \right )(x-y)}{|x-y|^{2}}\right.\\
                                &\quad\left.\times\frac{((x_{i_0}-y_{i_0})(x_{j_0}-y_{j_0})}
                              {|x-y|^{2}}\rho_{n}(x-y)dy \right |\nonumber\\
                                &\le N \int_{B(x,\delta)}\frac{|x-y|\epsilon|x-y|}{|x-y|^{2}}\frac{|x_{i_0}-y_{i_0}||x_{j_0}-y_{j_0}|}{|x-y|^{2}}\rho_{n}(x-y)dy\nonumber\\
                                &\le\epsilon N. \nonumber
\end{align}
Here, we have used the mean value theorem for scalar and vector valued functions to write
\begin{equation}\label{mean_value_duxy}
\mathpzc{d}^{2}u(x,y)=(x-y)^{T} \left (\int_{0}^{1}\int_{0}^{1}\nabla^{2}u(x+(t+s-1)(y-x))dsdt \right ) (x-y),
\end{equation}
and the fact that $\int_{\mathbb{R}^{N}} \rho_n(x)dx=1$ for all $n \in \mathbb{N}$.  This completes the proof in the case $i_0 \neq j_0$.

 \textbf{Subcase} $\mathbf{i_{0}= j_{0}}:$
 
 Let us record several observations before we proceed with this case. In fact, the same argument shows that for a single $i \in \{1,\ldots,N\}$
 \begin{align}
I^n_{i}(x) := \left |\int_{\mathbb{R}^{N}}\frac{\mathpzc{d}^{2}u(x,y)-(x-y)^{T}\nabla^{2}u(x)(x-y)}{|x-y|^{2}}\frac{(x_{i}-y_{i})^2}{|x-y|^{2}}\rho_{n}(x-y)dy\right | \to 0 \label{convergestozero}
 \end{align}
 uniformly in $x$ as $n \to \infty$, and therefore by summing in $i$ we deduce that
 \begin{align}
  \left |\int_{\mathbb{R}^{N}}\frac{\mathpzc{d}^{2}u(x,y)-(x-y)^{T}\nabla^{2}u(x)(x-y)}{|x-y|^{2}}\rho_{n}(x-y)dy\right | \to 0.  \label{laplacetozero}
 \end{align}
 Moreover, we observe that the same formula from Proposition \ref{prop:cij} and cancellation of odd powers implies
 \begin{align*}
\int_{\mathbb{R}^{N}} \frac{(x-y)^T\nabla^2u(x)(x-y) (x_{i_0}-y_{i_0})^2}{|x-y|^{4}}\rho_{n}(x-y)dy
&= \sum_{j=1}^N\frac{\partial^2 u}{\partial x^2_{j}}(x) \int_{\mathbb{R}^{N}} \frac{z^2_jz^2_{i_{0}}}{|z|^{4}}\rho_{n}(z)dz \\
&= \frac{1}{N(N+2)} \Delta u + \frac{2}{3}\frac{\partial^2 u}{\partial x^2_{i_{0}}}(x) \int_{\mathbb{R}^{N}} \frac{z_{i_{0}}^{4}}{|z|^{4}}\rho_{n}(z)dz \\
&= \frac{2}{N(N+2)} \left( \frac{1}{2}\Delta u +\frac{\partial^2 u}{\partial x^2_{i_{0}}}(x) \right),
\end{align*}
while we also have that
\begin{align*}
\int_{\mathbb{R}^{N}}\frac{(x-y)^{T}\nabla^{2}u(x)(x-y)}{|x-y|^{2}}\rho_{n}(x-y)dy &= \sum_{j=1}^N \frac{\partial^2 u}{\partial x^2_{j}}(x)  \int_{\mathbb{R}^{N}} \frac{z^2_j}{|z|^{2}}\rho_{n}(z)dz \\
&= \frac{1}{N} \Delta u(x).
\end{align*}
Thus, we can estimate
\begin{align*}
\bigg\vert \big(H_{n}u - \nabla^{2}u \big)_{(i_{0},i_{0})}\bigg\vert &\leq I_{i_0}^n(x) + \left |\frac{N}{2} \int_{\mathbb{R}^{N}}\frac{\mathpzc{d}^{2}u(x,y)}{|x-y|^2}\rho_{n}(x-y)dy -\int_{\mathbb{R}^N}\frac{\Delta u(x)}{2} \rho_{n}(x-y)dy\right| \\
&=  I_{i_0}^n(x) + \left |\frac{N}{2} \int_{\mathbb{R}^{N}}\frac{\mathpzc{d}^{2}u(x,y)-(x-y)^T\nabla^2u(x)(x-y)}{|x-y|^2}\rho_{n}(x-y)dy \right|,
\end{align*}
and the proof is completed by invoking the convergences established in equations \eqref{convergestozero} and \eqref{laplacetozero}.

\endproof

\subsection{Localization --  $W^{2,p}(\mathbb{R}^{N})$ case}\label{localisation_W2p}

The objective of this section is to show that if $u\in W^{2,p}(\mathbb{R}^{N})$, $1\le p<\infty$ then the non-local Hessian $H_{n}u$ converges to $\nabla^{2}u$ in $L^{p}$.  The first step is show that indeed in that case $H_{n}u$ is indeed an $L^{p}$ function. This follows from Lemma \ref{W2p_Hessian_label} which we prove next.

{\em Proof of Lemma \ref{W2p_Hessian_label}.}
Let us begin by making estimates for a function $v\in C^{\infty}(\mathbb{R}^{N})\cap W^{2,p}(\mathbb{R}^{N})$. From the definition of the non-local Hessian, and utilizing Jensen's inequality, equation \eqref{mean_value_duxy} as well as Fubini's theorem we have the following successive estimates (the constant is always denoted with $M(N,p)$):
\begin{align}
&\int_{\mathbb{R}^{N}}|H_{n}v(x)|_{\mathbb{R}^{N\times N}}^{p}dx\\
&= \left(\frac{N(N+2)}{2} \right)^{p}\int_{\mathbb{R}^{N}}\left |\int_{\mathbb{R}^{N}}\frac{\mathpzc{d}^{2}v(x,y)}{|x-y|^{2}}\frac{\left ((x-y)\otimes (x-y)-\frac{|x-y|^{2}}{N+2}I_{N}\right)}{|x-y|^{2}}\rho_{n}(x-y)dy \right |^{p}dx \nonumber\\
&\le M(N,p) \int_{\mathbb{R}^{N}}\left (\int_{\mathbb{R}^{N}}\frac{|\mathpzc{d}^{2}v(x,y)|}{|x-y|^{2}}\rho_{n}(x-y)dy \right )^{p}dx \nonumber
\end{align}
\begin{align}
&\le M(N,p)\int_{\mathbb{R}^{N}}\left ( \int_{\mathbb{R}^{N}}\frac{|\mathpzc{d}^{2}v(x,y)|^{p}}{|x-y|^{2p}}\rho_{n}(x-y)dy\right ) \left ( \underbrace{\int_{\mathbb{R}^{N}}\rho(x-y)dy}_{= 1}\right )^{p/p'}dx\label{key}\\
&\le M(N,p)\int_{\mathbb{R}^{N}}\left ( \int_{\mathbb{R}^{N}}\frac{|\int_{0}^{1}\nabla v (x+t(y-x)) -\nabla v(x+(t-1)(y-x))dt|^{p}|x-y|^{p}}{|x-y|^{2p}}\rho_{n}(x-y)dy\right )dx \nonumber\\
&\le M(N,p)\int_{\mathbb{R}^{N}}\left ( \int_{\mathbb{R}^{N}}\frac{\int_{0}^{1}|\nabla v (x+t(y-x)) -\nabla v(x+(t-1)(y-x))|^{p}dt}{|x-y|^{p}}\rho_{n}(x-y)dy\right )dx \nonumber\\
&\le M(N,p)\int_{\mathbb{R}^{N}}\left ( \int_{\mathbb{R}^{N}}\left (\int_{0}^{1}\int_{0}^{1}\left |\nabla^{2}v(x+(t+s-1)(y-x))\right |^{p} dsdt\right )\rho_{n}(x-y)dy\right )dx \nonumber\\
&= M(N,p)\int_{0}^{1}\int_{0}^{1}\int_{\mathbb{R}^{N}}\left ( \int_{\mathbb{R}^{N}}\left (\left |\nabla^{2}v(x+(t+s-1)(y-x))\right |^{p} \right )\rho_{n}(x-y)dy\right )dx dsdt \nonumber\\
&=M(N,p)\int_{0}^{1}\int_{0}^{1}\int_{\mathbb{R}^{N}}\left ( \int_{\mathbb{R}^{N}}\left (\left |\nabla^{2}v(x+(t+s-1)\xi)\right |^{p} \right )\rho_{n}(\xi)d\xi\right )dx dsdt, \nonumber\\
&= M(N,p)\int_{0}^{1}\int_{0}^{1}\int_{\mathbb{R}^{N}}\left ( \rho_{n}(\xi)\int_{\mathbb{R}^{N}}\left (\left |\nabla^{2}v(x+(t+s-1)\xi)\right |^{p} \right )dx\right )d\xi dsdt, \nonumber\\
&= M(N,p)\int_{0}^{1}\int_{0}^{1}\int_{\mathbb{R}^{N}} \rho_{n}(\xi)\|\nabla^{2}v\|_{L^{p}(\mathbb{R}^{N})}^{p} d\xi dsdt \nonumber\\
&=M(N,p)\|\nabla^{2}v\|_{L^{p}(\mathbb{R}^{N})}^{p} \nonumber.
\end{align}
Consider now a sequence $(v_{k})_{k\in \mathbb{N}}$ in $C^{\infty}(\mathbb{R}^{N})\cap W^{2,p}(\mathbb{R}^{N})$ approximating $u$ in $W^{2,p}(\mathbb{R}^{N})$. 
We already have from above that
\begin{equation}\label{keyforsmooth}
\int_{\mathbb{R}^{N}}\left ( \int_{\mathbb{R}^{N}}\frac{|\mathpzc{d}^{2}v_{k}(x,y)|^{p}}{|x-y|^{2p}}\rho_{n}(x-y)dy\right )dx\le M \|\nabla^{2}v_{k}\|_{L^{p}(\mathbb{R}^{N})},\quad \forall k\in \mathbb{N},
\end{equation}
or
\begin{equation}\label{keyforsmooth2}
\int_{\mathbb{R}^{N}}\left ( \int_{\mathbb{R}^{N}}\frac{|v_{k}(x+h)-2v_{k}(x)+v_{k}(x-h)|^{p}}{|h|^{2p}}\rho_{n}(h)dh\right )dx\le M \|\nabla^{2}v_{k}\|_{L^{p}(\mathbb{R}^{N})},\quad \forall k\in \mathbb{N}.
\end{equation}
Since $v_{k}$ converges to $u$ in $L^{p}(\mathbb{R}^{N})$ we have  that there exists a subsequence $v_{k_{\ell}}$ converging to $u$ almost everywhere.
From an application of Fatou's lemma we get that for every $x\in \mathbb{R}^{N}$
\[\underbrace{\int_{\mathbb{R}^{N}}\frac{|u(x+h)-2u(x)+u(x-h)|^{p}}{|h|^{2p}}\rho_{n}(h)dh}_{F(x)}\le M \liminf_{\ell\to\infty}\underbrace{\int_{\mathbb{R}^{N}}\frac{|v_{k_{\ell}}(x+h)-2v_{k_{\ell}}(x)+v_{k_{\ell}}(x-h)|^{p}}{|h|^{2p}}\rho_{n}(h)dh}_{F_{k_{\ell}}(x)}.\]
Applying one more time Fatou's Lemma to $F_{k_{\ell}}$, $F$ we get that
\begin{align*}
&\int_{\mathbb{R}^{N}}\int_{\mathbb{R}^{N}}\frac{|u(x+h)-2u(x)+u(x-h)|^{p}}{|h|^{2p}}\rho_{n}(h)dhdx\\
&\le M \liminf_{\ell\to\infty}\int_{\mathbb{R}^{N}}\int_{\mathbb{R}^{N}}\frac{|v_{k_{\ell}}(x+h)-2v_{k_{\ell}}(x)+v_{k_{\ell}}(x-h)|^{p}}{|h|^{2p}}\rho_{n}(h)dhdx\\
&\le M  \liminf_{\ell\to\infty}\|\nabla^{2}v_{k_{\ell}}\|_{L^{p}(\mathbb{R}^{N})}\\
&= M\|\nabla^{2}u\|_{L^{p}(\mathbb{R}^{N})}.
\end{align*}
This argument, along with Jensen's inequality, allows us to conclude that the conditions of Theorem \ref{NIByparts} are satisfied, in particular that $H_n u$ is well-defined as a Lebesgue integral, so that the estimate \eqref{key} holds for $W^{2,p}$ functions as well, thus completing the proof.
\endproof

We now have the necessary tools to prove the localization for $W^{2,p}$ functions.\\

{\em Proof of Theorem \ref{loc_W2p_label}.}
The result holds for functions $v\in C_{c}^{2}(\mathbb{R}^{N})$ since from Theorem \ref{localization_smooth_label} we have that $H_{n}v\to \nabla^{2}v$ in $L^p(\mathbb{R}^N,\mathbb{R}^{N\times N})$.  We use now the fact that that $C_{c}^{\infty}(\mathbb{R}^{N})$ and hence $C_{c}^{2}(\mathbb{R}^{N})$ is dense in $W^{2,p}(\mathbb{R}^{N})$, see for example \cite{brezis1983analyse}. Let $\epsilon>0$, then from density we have that there exists  a function $v\in C_{c}^{2}(\mathbb{R}^{N})$ such that
\[\|\nabla^{2}u-\nabla^{2}v\|_{L^{p}(\mathbb{R}^{N})}\le \epsilon.\]
Thus using also Lemma \ref{W2p_Hessian_label} we have
\begin{eqnarray*}
\|H_{n}u-\nabla^{2}u\|_{L^{p}(\mathbb{R}^{N})}&\le& \|H_{n}u-H_{n}v\|_{L^{p}(\mathbb{R}^{N})}+\|H_{n}v-\nabla^{2} v\|_{L^{p}(\mathbb{R}^{N})}+\|\nabla^{2}v-\nabla^{2}u\|_{L^{p}(\mathbb{R}^{N})}\\
&\le& C\epsilon +\|H_{n}v-\nabla^{2} v\|_{L^{p}(\mathbb{R}^{N})}+\epsilon.
\end{eqnarray*}
 Taking limits as $n\to \infty$ we get\[\limsup_{n\to\infty}\|H_{n}u-\nabla^{2}u\|_{L^{p}(\mathbb{R}^{N})} \le (C+1)\epsilon,\]
and thus we conclude that
\[\lim_{n\to\infty}\|H_{n}u-\nabla^{2}u\|_{L^{p}(\mathbb{R}^{N})}=0.\]
\endproof

\subsection{Localization -- $\mathrm{BV}^{2}(\mathbb{R}^{N})$ case}\label{localization_BV2}
Analogously with the first order case in \cite{Spector}, we can define a second order non-local divergence that corresponds to $H_{n}$, and we can also derive a second order non-local integration by parts formula which is an essential tool for the proofs of this section. The second order non-local divergence is defined for a function $\phi=(\phi_{ij})_{i,j=1}^{N}$ as 
\begin{equation}\label{D2}
\mathcal{D}_{n}^{2}\phi(x)=\frac{N(N+2)}{2}\int_{\mathbb{R}^{N}}\frac{\phi(y)-2\phi(x)+\phi(x+(x-y))}{|x-y|^{2}}\cdot \frac{\left ((x-y)\otimes (x-y)-\frac{|x-y|^{2}}{N+2}I_{N}\right)}{|x-y|^{2}}\rho_{n}(x-y)dy.
\end{equation}
where $A\cdot B=\sum_{i,j=1}^{N}A_{ij}B_{ij}$ for two $N\times N$ matrices $A$ and $B$. Notice that \eqref{D2} is well-defined for $\phi\in C_{c}^{2}(\mathbb{R}^{N},\mathbb{R}^{N\times N})$.\\

 \newtheorem{theorem12}[definition1]{Theorem}
\begin{theorem12}[Second order non-local integration by parts formula]\label{byparts2_label}
Suppose that $u\in L^{1}(\mathbb{R}^{N})$,\\ $\tfrac{|\mathpzc{d}^{2}u(x,y)|}{|x-y|^{2}}\rho_{n}(x-y)\in L^{1}(\mathbb{R}^{N}\times \mathbb{R}^{N})$ and let $\phi\in C_{c}^{2}(\mathbb{R}^{N},\mathbb{R}^{N\times N})$. Then
\begin{equation}\label{byparts2formula}
\int_{\mathbb{R}^{N}}H_{n}u(x)\cdot \phi(x)dx=\int_{\mathbb{R}^{N}}u(x)\mathcal{D}_{n}^{2} \phi(x)dx.
\end{equation} 
\end{theorem12}
In fact, this theorem can be deduced as a consequence of Theorem \ref{NIByparts} through a component by component application and collection of terms.
The following lemma shows the convergence of the second order non-local divergence to the continuous analogue $\mathrm{div}^{2}\phi$, where $\phi\in C_{c}^{2}(\mathbb{R}^{N},\mathbb{R}^{N\times N})$ and
\[\mathrm{div}^{2}\phi:=\sum_{i,j=1}^{N}\frac{\partial \phi_{ij}}{\partial x_{i}\partial x_{j}}.\]

 \newtheorem{lemma13}[definition1]{Lemma}
\begin{lemma13}\label{D2_uniform_label}
Let $\phi\in C_{c}^{2}(\mathbb{R}^{N},\mathbb{R}^{N\times N})$. Then for every $1\le p\le \infty$ we have
\begin{equation}\label{Dtodiv}
\lim_{n\to\infty} \|\mathcal{D}_{n}^{2}\phi - \mathrm{div}^{2}\phi\|_{L^{p}(\mathbb{R}^{N})}=0.
\end{equation}
\end{lemma13}
\begin{proof}
The proof follows immediately from Theorem \ref{localization_smooth_label} and \ref{loc_W2p_label}.
\end{proof}

Next we prove Lemma \ref{BV2_Hessian_label} that shows that the non-local Hessian \eqref{OUR_functional} is well defined for $u\in\mathrm{BV}^{2}(\mathbb{R}^{N})$. It is the analogue of Lemma \ref{W2p_Hessian_label} for functions in $\mathrm{BV}^{2}(\mathbb{R}^{N})$ this time.\\

{\em Proof of Lemma \ref{BV2_Hessian_label}.}
Let $(u_{k})_{k\in\mathbb{N}}$ be a sequence of functions in $C^{\infty}(\mathbb{R}^{N})$ that converges strictly in $\mathrm{BV}^{2}(\mathbb{R}^{N})$.
By the same calculations as in the proof of Lemma \ref{W2p_Hessian_label} we have for every $k\in \mathbb{N}$
\[\int_{\mathbb{R}^{N}}|H_{n}u_{k}(x)|dx\le M(N,1)\|\nabla^{2}u_{k}\|_{L^{1}(\mathbb{R}^{N})}.\]
Using Fatou's Lemma in a similar way as in Lemma \ref{W2p_Hessian_label}, we can establish that $H_nu$ is well-defined as a Lebesgue integral, along with the estimate
\begin{eqnarray*}
\int_{\mathbb{R}^{N}}|H_{n}u(x)|dx&\le& M(N,1)\liminf_{k\to\infty}|D^{2}u_{k}|(\mathbb{R}^{N})\\                                                                      &=&     M(N,1)|D^{2}u|(\mathbb{R}^{N}).
\end{eqnarray*}
However, employing the strict convergence of $D^2u_k$ to $D^2u$, the result is demonstrated.
\endproof

We can now proceed in proving the localization result for $\mathrm{BV}^{2}$ functions. We define $\mu_{n}$ to be the $\mathbb{R}^{N\times N}$-valued finite Radon measures $\mu_{n}:=H_{n}u\mathcal{L}^{N}$.\\

{\em Proof of Theorem \ref{loc_BV2_label}.}
Since $C_{c}^{\infty}(\mathbb{R}^{N},\mathbb{R}^{N\times N})$ is dense in $C_{0}(\mathbb{R}^{N},\mathbb{R}^{N\times N})$ it suffices to prove \eqref{loc_weakstar} for every $\psi\in C_{c}^{\infty}(\mathbb{R}^{N},\mathbb{R}^{N\times N})$. Indeed, suppose we have done this, let $\epsilon>0$ and let  $\phi \in C_{0}(\mathbb{R}^{N},\mathbb{R}^{N\times N})$
and $\psi\in C_{c}^{\infty}(\mathbb{R}^{N},\mathbb{R}^{N\times N})$ such that $\|\phi-\psi\|_{\infty}<\epsilon$. Then, using also the estimate \eqref{estimateBV2}, we have
\begin{align*}
\left | \int_{\mathbb{R}^{N}}H_{n}u(x)\cdot\phi(x)dx-\int_{\mathbb{R}^{N}}\phi(x) \,d D^{2}u\right |&\le \left |\int_{\mathbb{R}^{N}}H_{n}u(x)\cdot(\phi(x)-\psi(x))dx \right |\\
&\;\;\;\;+\left |\int_{\mathbb{R}^{N}}H_{n}u(x)\cdot\psi(x)dx-\int_{\mathbb{R}^{N}}\psi(x) \,d D^{2}u \right |\\
&\;\;\;\;+\left | \int_{\mathbb{R}^{N}}(\phi(x)-\psi(x)) \,d D^{2}u\right |\\
&\le \epsilon \int_{\mathbb{R^{N}}}|H_{n}u(x)|dx+\left |\int_{\mathbb{R}^{N}}H_{n}(x)\cdot\psi(x)dx-\int_{\mathbb{R}^{N}}\psi(x) \,d D^{2}u \right |\\
&\;\;\;\;+\epsilon |D^{2}u|(\mathbb{R}^{N})\\
&\le M\epsilon |D^{2}u|(\mathbb{R}^{N}) +\left |\int_{\mathbb{R}^{N}}H_{n}(x)\cdot\psi(x)dx-\int_{\mathbb{R}^{N}}\psi(x) \,d D^{2}u \right |\\
&\;\;\;\;+\epsilon |D^{2}u|(\mathbb{R}^{N}).
\end{align*}
 Taking the limit $n\to\infty$ from both sides of the above inequality we get
\[\limsup_{n\to\infty}\left | \int_{\mathbb{R}^{N}}H_{n}u(x)\cdot\phi(x)dx-\int_{\mathbb{R}^{N}}\phi(x) \,d D^{2}u\right |\le \tilde{M}\epsilon\]
and since $\epsilon$ is arbitrary we have \eqref{loc_weakstar}. We thus proceed to prove \eqref{loc_weakstar} for $C_{c}^{\infty}$ functions. From the estimate \eqref{estimateBV2} we have that $(|\mu_{n}|)_{n\in\mathbb{N}}$ is bounded, thus there exists a subsequence $(\mu_{n_{k}})_{k\in\mathbb{N}}$ and a $\mathbb{R}^{N\times N}$-valued Radon measure $\mu$ such that $\mu_{n_{k}}$ converges to $\mu$ weakly$^{\ast}$. This means that for every $\psi\in C_{c}^{\infty}(\mathbb{R}^{N},\mathbb{R}^{N\times N})$ we have
\[\lim_{k\to\infty}\int_{\mathbb{R}^{N}}H_{n_{k}}u(x)\cdot \psi(x)dx=\int_{\mathbb{R}^{N}}\psi(x)\cdot d\mu.\]
On the other hand from the integration by parts formula \eqref{byparts2formula} and Lemma \ref{D2_uniform_label} we get
\begin{eqnarray*}
\lim_{k\to\infty}\int_{\mathbb{R}^{N}}H_{n_{k}}u(x)\cdot \psi(x)dx&=&\lim_{k\to\infty} \int_{\mathbb{R}^{N}} u(x)\mathcal{D}_{n_{k}}^{2}\psi(x)dx\\
&=&\int_{\mathbb{R}^{N}}u(x)\mathrm{div}^{2}\psi(x)dx\\
&=&\int_{\mathbb{R}^{N}}\psi(x)\cdot dD^{2}u.
\end{eqnarray*}
This means that $\mu=D^{2}u$. Observe now that since we actually deduce that  every subsequence of $(\mu_{n})_{n\in\mathbb{N}}$ has a further subsequence that converges to $D^{2}u$ weakly$^{\ast}$, then that the initial sequence $(\mu_{n})_{n\in \mathbb{N}}$ converges to $D^{2}u$ weakly star.
\endproof

Let us note here that in the case $N=1$, we can also prove strict convergence of the measures $\mu_{n}$ to $D^{2}u$, that is, in addition to \eqref{loc_weakstar} we also have
\[|\mu_{n}|(\mathbb{R})\to |D^{2}u|(\mathbb{R}).\]
 \newtheorem{theorem14}[definition1]{Theorem}
\begin{theorem14}
$\,$Let $N=1$. Then the sequence $(\mu_{n})_{n\in\mathbb{N}}$ converges to $D^{2}u$ strictly as measures, i.e.,
\begin{equation}\label{1Dwstar}
\mu_{n}\to D^{2}u,\quad \text{weakly}^{\ast},
\end{equation}
and
\begin{equation}\label{1Dstrict}
|\mu_{n}|(\mathbb{R})\to |D^{2}u|(\mathbb{R}).
\end{equation}
\end{theorem14}
\begin{proof}
The weak${^\ast}$ convergence was proven in Theorem \ref{loc_BV2_label}. Since in the space of finite Radon measures, the total variation norm is lower semicontinuous with respect to the weak$^{\ast}$ convergence, we also have
\begin{equation}\label{strict1D_liminf}
 |D^{2}u|(\mathbb{R})\le \liminf_{n\to\infty} |\mu_{n}|(\mathbb{R}).
\end{equation}
Thus it suffices to show that
\begin{equation}\label{strict1D_limsup}
\limsup_{n\to\infty} |\mu_{n}|(\mathbb{R})\le  |D^{2}u|(\mathbb{R}).
\end{equation}
Note that in dimension one the non-local Hessian formula is
\begin{equation}\label{NL_hessian_1d}
H_{n}u(x)=\int_{\mathbb{R}}\frac{u(y)-2u(x)+u(x+(x-y))}{|x-y|^{2}}\rho_{n}(x-y)dy.
\end{equation}
Following the proof of Lemma \ref{W2p_Hessian_label}, we can easily verify that for $v\in C^{\infty}(\mathbb{R})\cap \mathrm{BV}^{2}(\mathbb{R})$ we have
\[\int_{\mathbb{R}}|H_{n}v(x)|dx\le \|\nabla^{2}v\|_{L^{1}(\mathbb{R})},\]
i.e., the constant $M$ that appears in the estimate \eqref{estimateW2p} is equal to 1. Using Fatou's Lemma and the $\mathrm{BV}^{2}$ strict approximation of $u$ by smooth functions we get that
\[|\mu_{n}|(\mathbb{R})=\int_{\mathbb{R}}|H_{n}u(x)|dx\le|D^{2}u|(\mathbb{R}),\]
from where \eqref{strict1D_limsup} straightforwardly follows.
\end{proof}

\subsection{Characterization of higher-order Sobolev and BV spaces} \label{characterizations}
Characterization of Sobolev and $\mathrm{BV}$ spaces in terms of non-local, derivative-free energies has been done so far only in the first order case, see \cite{bourgain2001another, ponce2004new, mengesha2012nonlocal, Spector}. Here we characterize the spaces $W^{2,p}(\mathbb{R}^{N})$ and $\mathrm{BV}^{2}(\mathbb{R}^{N})$ using our definition of non-local Hessian.\\

{\em Proof of Theorem \ref{char}.}
Firstly, we prove \eqref{W2p_char}. Suppose that $u\in W^{2,p}(\mathbb{R}^{N})$. Then, Lemma \ref{W2p_Hessian_label} 
gives 
\[\liminf_{n\to\infty}\int_{\mathbb{R}^{N}}|H_{n}u(x)|^{p}dx\le M\|\nabla^{2}u\|_{L^{p}(\mathbb{R}^{N})}^{p}<\infty.\]
Suppose now conversely that
\[\liminf_{n\to\infty} \int_{\mathbb{R}^{N}}|\mathfrak{H}_{n}u(x)|^{p}dx<\infty.\]
This means that up to a subsequence, $\mathfrak{H}_{n}u$ is bounded in $L^{p}(\mathbb{R}^{N},\mathbb{R}^{N\times N})$, thus there exists a subsequence $(\mathfrak{H}_{n_{k}}u)_{k\in\mathbb{N}}$ and $v\in L^{p}(\mathbb{R}^{N},\mathbb{R}^{N\times N})$ such that $\mathfrak{H}_{n_{k}}u\rightharpoonup v$ weakly in $L^{p}(\mathbb{R}^{N},\mathbb{R}^{N\times N})$. Thus, using the definition of $L^{p}$ weak convergence together with the definition of $\mathfrak{H}_nu$, we have for every $\psi\in C_{c}^{\infty}(\mathbb{R}^{N})$,
\begin{align*}
\int_{\mathbb{R}^{N}}v^{ij}(x)\cdot \psi(x)dx&=\lim_{k\to\infty}\int_{\mathbb{R}^{N}}\mathfrak{H}^{ij}_{n_{k}}u(x)\cdot \psi(x)\\
&= \lim_{k\to\infty}\int_{\mathbb{R}^{N}}u(x)H^{ij}_n\psi(x)\\
&= \int_{\mathbb{R}^{N}}u(x)\frac{\partial^2 \psi(x)}{\partial x_i\partial x_j}dx,
\end{align*}
something that shows that $v=\nabla^{2}u$ is the second order weak derivative of $u$. Now since $u \in L^p(\mathbb{R}^N)$ and the second order distributional derivative is a function, mollification of $u$ and the Gagliardo-Nirenberg inequality (see \cite[p. 128, Equation 2.5]{Nirenberg})
\begin{equation}
\|\nabla u\|_{L^{p}(\mathbb{R}^{N},\mathbb{R}^{N})}\le C \|\nabla^{2}u\|_{L^{p}(\mathbb{R}^{N},\mathbb{R}^{N\times N})}^{\frac{1}{2}}\|u\|_{L^{p}(\mathbb{R}^{N})}^{\frac{1}{2}}, \label{GN}
\end{equation}
implies that the first distributional derivative and it belongs to $L^{p}(\mathbb{R}^N,\mathbb{R}^N)$ and thus $u\in W^{2,p}(\mathbb{R}^{N})$. 

We now proceed in proving \eqref{BV2_char}. Again supposing that $u\in \mathrm{BV}^{2}(\mathbb{R}^{N})$ we have that Lemma \ref{BV2_Hessian_label} gives us
\[\liminf_{n\to\infty}\int_{\mathbb{R}^{N}}|H_{n}u(x)|dx\le C|D^{2}u|(\mathbb{R}^{N}).\]
Suppose now that 
\[\liminf_{n\to\infty} \int_{\mathbb{R}^{N}}|\mathfrak{H}_{n}u(x)|dx<\infty.\]
Considering again the measures $\mu_{n}=\mathfrak{H}_{n}u\mathcal{L}^{N}$ we have that that there exists a subsequence $(\mu_{n_{k}})_{k\in\mathbb{N}}$ and a finite Radon measure $\mu$ such that $\mu_{n_{k}}\overset{*}{\rightharpoonup}\mu$ weakly$^{\ast}$. Then for every $\psi \in C_{c}^{\infty}(\mathbb{R}^{N})$ we have, similarly as before,
\begin{eqnarray*}
\int_{\mathbb{R}^{N}}\psi d\mu^{ij}&=&\lim_{k\to\infty}\int_{\mathbb{R}^{N}}\mathfrak{H}^{ij}_{n_{k}}u(x)\cdot \psi(x)dx\\
&=&\lim_{k\to\infty}\int_{\mathbb{R}^{N}}u(x) H^{ij}_n\psi(x)dx\\
&=&\int_{\mathbb{R}^{N}}u(x)\frac{\partial^2 \psi(x)}{\partial x_i \partial x_j} dx.
\end{eqnarray*}
something that shows that $\mu=D^{2}u$. Again, by first mollifying and then passing the limit, the inequality \eqref{GN} implies that $D u \in \mathcal{M}(\mathbb{R}^{N},\mathbb{R}^{N})$.  However, $D u \in \mathcal{M}(\mathbb{R}^{N},\mathbb{R}^{N})$ and $D^2 u \in \mathcal{M}(\mathbb{R}^{N},\mathbb{R}^{N\times N})$ implies that $Du$ is an $L^1(\mathbb{R}^N, \mathbb{R}^N)$ function (see \cite[Exercise 3.2]{AmbrosioBV}), and we therefore conclude that $u\in\mathrm{BV}^{2}(\mathbb{R}^{N})$. \endproof

\subsection{Gamma convergence} \label{section:gamma}

For notational convenience we define the functional
\begin{equation}
F_n(u):= \int_{\mathbb{R}^N} |\mathfrak{H}_n u|\;dx.
\end{equation}

{\em Proof of Theorem \ref{thm:explicitgammaconv}.}
The computation of the Gamma limit consists of two inequalities.  For the lower bound, we must show that
\begin{equation*}
|D^{2}u|(\mathbb{R}^{N})\leq \liminf_{n\to \infty} F_n(u_n)
\end{equation*}
for every sequence $u_n \to u$ in $L^1(\mathbb{R}^N)$.  Without loss of generality we may assume that
\begin{equation*}
C:=\liminf_{n\to \infty} F_n(u_n) < +\infty.
\end{equation*}
which implies that 
\begin{equation*}
\sup _\varphi \liminf_{n \to \infty} \left|\int_{\mathbb{R}^N} \mathfrak{H}^{ij}_n u_n \varphi \;dx \right| \leq C.
\end{equation*}
where the supremum is taken over $\varphi \in C^\infty_c(\mathbb{R}^N)$ such that $\|\varphi\|_{L^\infty(\mathbb{R}^N)} \leq1$.  Now, the definition of the distributional non-local Hessian and the convergence $u_n \to u$ in $L^1(\mathbb{R}^N)$ imply that
\begin{align*}
\lim_{n \to \infty} \int_{\mathbb{R}^N} \mathfrak{H}^{ij}_n u_n  \varphi \;dx &=  \lim_{n \to \infty} \int_{\mathbb{R}^N} u_n  H^{ij}_n \varphi \;dx  \\
&= \int_{\mathbb{R}^N} u \frac{ \partial^2 \varphi}{\partial x_i \partial x_j} \;dx.
\end{align*}
We thus conclude that
\begin{equation*}
\sup _\varphi \liminf_{n \to \infty} \left|\int_{\mathbb{R}^N} u \frac{ \partial^2 \varphi}{\partial x_i \partial x_j}\;dx \right| \leq C,
\end{equation*}
which, arguing as in the previous section, says that $u \in \mathrm{BV}^2(\mathbb{R}^N)$, in particular that $D^2 u \in \mathcal{M}(\mathbb{R}^{N},\mathbb{R}^{N\times N})$ and
\begin{equation*}
|D^{2}u|(\mathbb{R}^{N}) \leq \Gamma^-_{L^1(\mathbb{R}^N)} \mhyphen \liminf_{n\to\infty} F_{n}(u)
\end{equation*}
for every $u \in L^1(\mathbb{R}^N)$.

For the upper bound we observe that if $u \in C^2_c(\mathbb{R}^N)$, we have by the uniform convergence of Theorem \ref{localization_smooth_label} and the fact that $u$ is sufficiently smooth with compact support that
\begin{equation*}
 \lim_{n \to \infty} F_n(u)= |D^{2}u|(\mathbb{R}^{N}).
\end{equation*}
Then choosing $u_n=u$ we conclude
\begin{align*}
\Gamma^-_{L^1(\mathbb{R}^N)} \mhyphen \limsup_{n\to\infty} F_{n}(u) &\leq  \lim_{n \to \infty} F_n(u)\\
&= |D^{2}u|(\mathbb{R}^{N}).
\end{align*}
Now, taking the lower semicontinuous envelope with respect to $L^1(\mathbb{R}^N)$ strong convergence, using $\Gamma^-_{L^1(\mathbb{R}^N)}\mhyphen \limsup$ is lower semicontinuous on $L^1(\mathbb{R}^N)$ (c.f. \cite[Proposition 6.8]{dalmasogamma}),  we deduce that
\begin{align*}
\Gamma^-_{L^1(\Omega)} \mhyphen \limsup_{n\to\infty} F_{n}(u) &\leq sc^-_{L^1(\Omega)}|D^{2}u|(\mathbb{R}^{N}). \\
&= |D^{2}u|(\mathbb{R}^{N})
\end{align*}
for all $u \in L^1(\mathbb{R}^N)$.  \endproof

\section{Extensions and Applications}
\subsection{An asymmetric extension}\label{sec:implicitdef}
In the previous sections we have shown that our non-local definition of $H_n$ as in \eqref{OUR_functional} localizes to the classical distributional Hessian for a specific choice of the weights $\rho_n$, and thus can be rightfully called a non-local Hessian. In numerical applications, however, the strength of such non-local models lies in the fact that the weights can be chosen to have non-local interactions and model specific patterns in the data. A classic example is the non-local total variation \cite{gilboa2008non-local}:
\beq\label{eq:oshernltv}
J_{NL-TV}(u) = \int_{\Om} \int_{\Om} |u(x) - u(y)| \sqrt{w(x,y)} d y d x.
\eeq
A possible choice is to choose $w(x,y)$ large if the \emph{patches} (neighbourhoods) around $x$ and $y$ are similar with respect to a \emph{patch distance} $d_a$, such as a weighted $\ell^2$ norm, and small if they are not. In \cite{gilboa2008non-local} this is achieved by setting 
$w(x,y) = 1$ if the neighbourhood around $y$ is one of the $K \in \N$ closest to the neighbourhood around $x$ in a search window, and $w(x,y) = 0$ otherwise. In effect, if the image contains a repeating pattern with a defect that is small enough to not throw off the patch distances too much, it will be repaired as long as most similar patterns do not show the defect.

Computing suitable weights is much less obvious in the case of $H$. We can formally extend \eqref{OUR_functional} using arbitrary pairwise weights $\rho:\Rn\times\Rn \to\R$,
\beq\label{eq:hrho}
H_\rho u(x)=C\int_{\mathbb{R}^{N}}\frac{u(x+z)-2u(x)+u(x-z)}{|z|^{2}}\frac{\left( z z^{\top}-\frac{|z|^{2}}{N+2}I_{N}\right )}{|z|^{2}}\rho_x(z) dz,
\eeq
and use it to create non-local generalisations of functionals such as TV$^2$, for example to minimise the non-local $L^2-\TV^2$ model
\beq
f(u) \assign \int_{\RN} | u - g |^2 d x +\alpha\int_{\RN} | \HH_\rho | d x.
\eeq
However, apart from being formulated on $\RN$ instead of $\Om$, formulation \eqref{eq:hrho} has an important drawback compared to the first-order formulation \eqref{eq:oshernltv}: while the weights are defined between two points $x$ and $y$, the left part of the integrand uses the values of $u$ not only at $x$ and $y$, but also at the ``mirrored'' point $x+(x-y)$. In fact we can replace the weighting function by the symmetrised version $\half\{\rho_{_{x}}(y-x) + \rho_{x}(x-y)\}$, which in effect relates three points instead of two, and limits the choice of possible weighting functions.

In this section we therefore introduce a more versatile extension of \eqref{OUR_functional} that allows for full non-symmetric weights. We start with the realisation that the finite-differences integrand in \eqref{eq:hrho} effectively comes from canceling the first-order differences in the Taylor expansion of $u$ around $x$, which couples the values of $u$ at $x$, $y$, and $x+(x-y)$ into one term. Instead, we can avoid this coupling by directly defining the non-local gradient $G'_u(x)$ and Hessian looking for a non-local gradient $G_u(x) \in \R^N$ and Hessian $H_u(x) \in \tmop{Sym}(\RNN)$ that \emph{best explain $u$ around $x$ in terms of a quadratic model}, i.e., that take the place of the gradient and Hessian in the Taylor expansion:
\begin{equation}\label{eq:implicitnl}
(G'_u(x),H'_u(x)) \assign \underset{G_u\in\RN,H_u \in \tmop{Sym}(\RNN)}{\operatorname{argmin}} \half  \int_{\Om - \{x\}} \left( u ( x + z ) -u ( x ) -G_{u} ^{\top}  z - \frac{1}{2}  z^{\top} H_{u} z
  \right)^{2} \sigma_{x} (z)  d z.
\end{equation}
Here the variable $x+z$ takes the place of $y$ in $\eqref{OUR_functional}$. We denote definition \eqref{eq:implicitnl} the \emph{implicit} non-local Hessian as opposed to the \emph{explicit} formulation \eqref{eq:hrho}.

The advantage is that any terms involving $\sigma_x(z)$ are now only based on the two values of $u$ at $x$ and $y=x+z$, and (in particular bounded) domains other than $\RN$ are naturally dealt with, which is important for a numerical implementation. We also note that this approach allows to incorporate non-local first-order terms as a side-effect, and can be naturally extended to third- and higher-order derivatives, which we leave to further work.

With respect to implementation, the implicit model \eqref{eq:implicitnl} does not add much to the overall difficulty: it is enough to add the non-local gradient and Hessian $G'_u(x) \in \RN$ and $H'_u(x)\in\RNN$ as additional variables to the problem, and couple them to $u$ by adding the optimality conditions for \eqref{eq:implicitnl} to the problem. Since \eqref{eq:implicitnl} is a finite-dimensional quadratic minimization problem, the optimality conditions are linear, i.e., of the form $A_{u,x} H'_u(x) = b_{u,x}$. Such linear constraints can be readily added to most solvers capable of solving the local problem. Alternatively the matrices $A_{u,x}$ can be inverted explicitly in a pre-computation step, however we did not find this to increase overall performance.

While the implicit and explicit model look different at first glance, from the considerations about the Taylor expansion we expect them to be closely related, and in fact this is true as shown in the following series of propositions. The first proposition explicitly computes the constants $C_{i j}$ and is needed for the proofs:\\

\newtheorem{proposition15}[definition1]{Proposition}
\begin{proposition15}
  \label{prop:cij} For $N \geqs 1$ we have the following definition and identity for $C$:
  \begin{eqnarray*}
    C_{i_{0}  j_{0}} \assign \int_{\Ss^{N-1}} \nu_{i_{0}}^{2} \nu_{j_{0}}^{2}  d
    \mathcal{H}^{N-1} ( x ) & = & \frac{|S^{N-1}| }{N ( N+2 )} \cdot \left\{ \begin{array}{ll}
      1, & i_{0} \neq j_{0} ,\\
      3, & i_{0} =j_{0} .
    \end{array} \right.\label{eq:cijexplicit}
  \end{eqnarray*}
  The matrix $C= ( C_{i_{0}  j_{0}} )_{i_{0} ,j_{0} =1, \ldots,N}$ is
  \begin{eqnarray*}
    C & = & C_{1 2}  ( E+2I ) ,
  \end{eqnarray*}
  where $E \in \RNN$ is the all-ones matrix.
\end{proposition15}

\begin{proof}See Appendix.
\end{proof}

The following theorem shows that for radially symmetric $\rho_x$ and $\sigma_x$, if we restrict ourselves to a circle, the \emph{implicit} non-local Hessian is (up to a factor) the \emph{explicit} non-local Hessian as defined earlier.\\

\newtheorem{theorem16}[definition1]{Theorem}
\begin{theorem16}\label{thm:implexplcircle}
Let $\Om \assign \RN$, $h > 0$, $x \in \RN$. Consider the minimization problem
\beqa\label{eq:minprobsphere}
&&\min_{G_{u,h}\in\RN,H_{u,h} \in \tmop{Sym}(\RNN)} \half  \int_{h \Ss^{N-1}} \left( u ( x + z ) -u ( x ) -G_{u,h} ^{\top}  z - \frac{1}{2}  z^{\top} H_{u,h} z
  \right)^{2} \dHHN(z),
\eeqa
and assume that the integral is finite for every $G_{u,h} \in \R^N$ and $H_{u,h} \in \R^{N\times N}$. Then the non-local Hessian part of the minimizing pair $(G'_{u,h},H'_{u,h})$
is given by
\begin{equation}
H'_{u,h} = \frac{1}{2} C_{1 2}^{-1} h^{-(N-1)}  \int_{h \Ss^{N-1}} \left( \frac{u
  ( x+z ) -2u ( x ) +u ( x-z )}{|z|^{2}} \right)  \left( \frac{z z^{\top}}{|z|^2} - \frac{1}{N+2} I \right)  \dHHN(z). \label{eq:phpimplicit}
\end{equation}
\end{theorem16}
\begin{proof}
The optimality conditions for (\ref{eq:minprobsphere}) in terms of $G_{u} = G_{u,h}$ and $H_{u} = H_{u,h}$, leaving the assumption on the symmetry
of $H_{u}$ aside for the moment as we will enforce them explicitly later, are
\begin{eqnarray*}
  0 & = &  \int_{h \Ss^{N-1}} - \left( u ( x+z ) -u ( x ) -G_{u}^{\top} z-
  \frac{1}{2} z^{\top} H_{u} z \right) z \dHHN(z),\\
  0 & = & \int_{h \Ss^{N-1}} - \frac{1}{2}  \left( u ( x+z ) -u ( x )
  -G_{u}^{\top} z- \frac{1}{2} z^{\top} H_{u} z \right) zz\dHHN(z),
\end{eqnarray*}
Thus from the first line we get $N$ equations
\begin{eqnarray*}
  &  & \left\langle \int_{h \Ss^{N-1}} z_{i} z\dHHN(z),G_{u} \right\rangle + \left\langle \int_{h \Ss^{N-1}} z_{i}  \left(
  \frac{1}{2} zz^{\top} \right)   \dHHN(z),H_{u} \right\rangle\\
  & = & \int_{h \Ss^{N-1}} z_{i}  ( u ( x+z ) -u ( x ) )\dHHN(z),  \hspace{1em} i=1, \ldots ,N.
\end{eqnarray*}
Similarly, from the second line we get $N^{2}$ equations (of which some may be
redundant)
\begin{eqnarray*}
  &  & \left\langle \int_{h \Ss^{N-1}} \frac{1}{2} z_{i} z_{j} z^{\top}
  \dHHN(z),G_{u} \right\rangle + \left\langle
  \int_{h \Ss^{N-1}} \frac{1}{2} z_{i} z_{j}  \left( \frac{1}{2} zz^{\top}
  \right)  \dHHN(z),H_{u} \right\rangle\\
  & = & \int_{h \Ss^{N-1}} \frac{1}{2} z_{i} z_{j}  ( u ( x+z ) -u ( x )
  )\dHHN(z), \hspace{1em} i,j=1, \ldots ,N.
\end{eqnarray*}
Note that $G_{u}$ and $H_{u}$ are elements of $\RN$ and
$\RNN$, and both sides of the inner products are
finite-dimensional vectors respective matrices.

If we collect the entries of $G_u$ and $H_u$ in a vector $p = (p_G, p_H)$, these two sets of equations can be rewritten as a linear system with an $m \times m$ block matrix,
\begin{eqnarray*}
  \left(\begin{array}{ll}
    A & V^{\top}\\
    V & B
  \end{array}\right)  p & = & \left(\begin{array}{c}
    a\\
    b
  \end{array}\right) .
\end{eqnarray*}
The entries in the sub-matrices V are all of the form
\begin{eqnarray*}
  \int_{h \Ss^{N-1}} z_{i}  ( z_{j} z_{k} )\dHHN(z), &  &
  \quad i,j,k \in \{ 1, \ldots ,n \} .
\end{eqnarray*}
No matter what the choice of $i,j,k$ is, there is always at least one index
with an odd power, so every single one of these integrals is zero due to
symmetry. This means that the conditions on the gradient and
the Hessian parts of $p$ decouple, i.e., the problem is
\begin{eqnarray}
  Ap_{G} =a, &  & \quad Bp_{H} =b.  \label{eq:bph}
\end{eqnarray}
We can therefore look at the isolated problem of computing the Hessian part $p_H$, or equivalently $H_u$, without interference from the gradient part.
The matrix $B$ is of the form
\begin{eqnarray*}
  \int_{h \Ss^{N-1}} \frac{1}{2} z_{i} z_{j}  \frac{1}{2} z_{i_{0}}
  z_{j_{0}} \dHHN(z). &  & 
\end{eqnarray*}
Again due to symmetry, the only way that this integral is non-zero is if there
are no odd powers, so either all indices are the same or there are two pairs:
\begin{eqnarray*}
  &  & \frac{1}{4}  \int_{h \Ss^{N-1}} z_{k}^{4} \dHHN(z),
  \hspace{1em} \text{for some } k,\\
  \text{or }\quad&  & \frac{1}{4}  \int_{h \Ss^{N-1}} z_{k}^{2}
  z_{l}^{2} \dHHN(z), \hspace{1em} \text{for some } k \neq l.
\end{eqnarray*}
The vector $b$ in (\ref{eq:bph}) is
\begin{eqnarray*}
  &  & \int_{h \Ss^{N-1}} ( u ( x+z ) -u ( x ) )  \frac{1}{2} z_{i} z_{j}
   \dHHN(z) = \frac{1}{2}  \int_{h \Ss^{N-1}} ( u ( x+z ) -2u ( x ) +u ( x-z ) )
  \frac{1}{2} z_{i} z_{j}   \dHHN(z).
\end{eqnarray*}
Now the optimality conditions on the $p_{H} = ( p_{i j} )_{i,j=1, \ldots ,N}$ part of $p$
are of the form (for all $i_{0} ,j_{0}$)
\begin{eqnarray}
  &  &  \sum_{i,j} p_{i j} \int_{h \Ss^{N-1}} \frac{1}{2} z_{i} z_{j} 
  \frac{1}{2} z_{i_{0}} z_{j_{0}}  \dHHN(z) \nonumber\\
  & = & \frac{1}{2}  \int_{h \Ss^{N-1}} ( u ( x+z ) -2u ( x ) +u ( x-z )
  )  \frac{1}{2} z_{i_{0}} z_{j_{0}}   \dHHN(z), \hspace{1em} i_{0}
  ,j_{0} =1, \ldots ,N,\label{eq:phstartingpoint}
\end{eqnarray}
or alternatively
\begin{eqnarray*}
  &  &  \sum_{i,j} p_{i j} \int_{\Ss^{N-1}} \frac{1}{2} h \nu_{i} h \nu_{j} 
  \frac{1}{2} h \nu_{i_{0}} h \nu_{j_{0}}  \dHHN(\nu)\\
  & = & \frac{1}{4}  \int_{\Ss^{N-1}} ( u ( x+h \nu ) -2u ( x ) +u ( x-h \nu )
  ) h \nu_{i_{0}} h \nu_{j_{0}} \dHHN(\nu), \hspace{1em} i_{0} ,j_{0} =1, \ldots
  ,N.  
\end{eqnarray*}
We can divide everything by the constant $\frac{1}{4} h^{4}$ to get
\begin{eqnarray}
  &  &  \sum_{i,j} p_{i j} \int_{\Ss^{N-1}} \nu_{i} \nu_{j} \nu_{i_{0}}
  \nu_{j_{0}} \dHHN(\nu) \nonumber\\
  & = & \int_{\Ss^{N-1}} \left( \frac{u ( x+h \nu ) -2u ( x ) +u ( x-h \nu
  )}{h^{2}} \right) \nu_{i_{0}} \nu_{j_{0}}  \dHHN(\nu) , \hspace{1em} i_{0}
  ,j_{0} =1, \ldots ,N.  \label{eq:nuij}
\end{eqnarray}
We consider the different cases for $i_{0}, j_{0}$ separately.
\begin{itemize}
  \item If $i_{0} \neq j_{0}$, then on the left-hand side all terms vanish
  except for $( i,j ) = ( i_{0} ,j_{0} )$ and $( i,j ) = ( j_{0} ,i_{0} )$.
  Overall we get the condition
  \begin{eqnarray*}
    ( p_{i_{0}  j_{0}} +p_{j_{0}  i_{0}} )  \int_{\Ss^{N-1}} \nu_{i_{0}}^{2}
    \nu_{j_{0}}^{2}  \dHHN(\nu) & = & \int_{\Ss^{N-1}} \left( \frac{u ( x+h \nu )
    -2u ( x ) +u ( x-h \nu )}{h^{2}} \right) \nu_{i_{0}} \nu_{j_{0}}  \dHHN(\nu).
  \end{eqnarray*}
  Since we require $H_u = p_{i j}$ to be symmetric, we get an
  explicit solution for $p_{i_{0} j_{0}}$:
  \begin{eqnarray*}
    p_{i_{0} j_{0}} & = & \frac{1}{2} C_{i_{0}  j_{0}}^{-1}  \int_{\Ss^{N-1}}
    \left( \frac{u ( x+h \nu ) -2u ( x ) +u ( x-h \nu )}{h^{2}} \right)
    \nu_{i_{0}} \nu_{j_{0}}  \dHHN(\nu) ,
  \end{eqnarray*}
  where $C$ is the scalar as in Prop.~\ref{prop:cij}:
  \begin{eqnarray*}
    C_{i_{0}  j_{0}} & \assign & \int_{\Ss^{N-1}} \nu_{i_{0}}^{2}
    \nu_{j_{0}}^{2}  \dHHN(\nu) .
  \end{eqnarray*}
  \item If $i_{0} =j_{0}$, the on the left-hand side all terms vanish except
  for $i=j$ and we get
  \begin{eqnarray*}
    \sum_{i} p_{i i}  \int_{\Ss^{N-1}} \nu_{i}^{2} \nu_{i_{0}}^{2}  \dHHN(\nu) & = &
    \int_{\Ss^{N-1}} \left( \frac{u ( x+h \nu ) -2u ( x ) +u ( x-h \nu )}{h^{2}}
    \right) \nu_{i_{0}} \nu_{j_{0}}  \dHHN(\nu) .
  \end{eqnarray*}
  With the definition for $C$ we can simplify the left-hand side to
  \begin{eqnarray*}
    \sum_{i} p_{i i} C_{i i_{0}} & = & \int_{\Ss^{N-1}} \left( \frac{u ( x+h \nu
    ) -2u ( x ) +u ( x-h \nu )}{h^{2}} \right) \nu_{i_{0}} \nu_{i_{0}}  \dHHN(\nu)
    , \hspace{1em} i_{0} =1, \ldots ,n.
  \end{eqnarray*}
  Consequentially, for the diagonal vector $p' = ( p_{11} , \ldots ,p_{N N} )$
  and $y= ( \tmop{rhs} )_{i_{0}}$ we get
from Proposition~\ref{prop:cij}  \begin{eqnarray*}
    C_{1 2}  ( E+2I ) p' & = & y
  \end{eqnarray*}
  with $E$ the all-ones matrix in $\mathbbm{R}^{N \times N}$. We have $( E+2I
  )^{-1} = \frac{1}{2}  \left( I- \frac{1}{N+2} E \right)$, thus we can
  rewrite
  \begin{eqnarray*}
    C_{1 2} p' & = & \frac{1}{2}  \left( I- \frac{1}{N+2} E \right) y.
  \end{eqnarray*}
  The $E$ term in the right-hand side is (since $| \nu |^{2} =1$ on $\Ss^{N-1}$)
  \begin{eqnarray*}
    &  & - \frac{1}{N+2}  \sum_{i} \int_{\Ss^{N-1}} \left( \frac{u ( x+h \nu )
    -2u ( x ) +u ( x-h \nu )}{h^{2}} \right) \nu_{i}^{2}  \dHHN(\nu)\\
    & = & - \frac{1}{N+2} _{} \int_{\Ss^{N-1}} \left( \frac{u ( x+h \nu ) -2u (
    x ) +u ( x-h \nu )}{h^{2}} \right) 1\dHHN(\nu) .
  \end{eqnarray*}
\end{itemize}
Putting both cases together, we can now solve the least-squares problem for a
fixed radius $h$ in closed form for $p_{H}$ and therefore the first part of the claimed identity:
\begin{equation}
  p_{H} ( h ) =p  =  \frac{1}{2} C_{1 2}^{-1}  \int_{\Ss^{N-1}} \left( \frac{u
  ( x+h \nu ) -2u ( x ) +u ( x-h \nu )}{h^{2}} \right)  \left( \nu \nu^{\top}
  - \frac{1}{N+2} I \right)  \dHHN(\nu)
\end{equation}
The last part follows with the substitution $z = h \nu$.
\end{proof}

The following theorem shows that the explicit form $H_u$ is the same as \emph{separately} solving
the least-squares problem on all spheres for all $h > 0$ and then taking a weighted mean based on $\rho$.\\

\newtheorem{proposition17}[definition1]{Proposition}
\begin{proposition17}\label{prop:implexplconn}
Assume $\rho_x$ is radially symmetric, i.e., $\rho_x(z) = \rho_x(|z|)$. Then the explicit non-local Hessian $H_u$ can be written as
\beqa
H_u & = & \int_{0}^{\infty} H'_{u,h} \int_{hS^{N-1}} \rho_{x} ( z )  d \mathcal{H}^{N-1} ( z )  d h.\label{eq:implexplconn}
\eeqa
\end{proposition17}
\begin{proof}
This is a direct consequence of Theorem.~\ref{thm:implexplcircle}:
\begin{eqnarray*}
  H_u & =&  \frac{N(N+2)}{2}\int_{\RN}\frac{u(x+z)-2u(x)+u(x-z)}{|z|^{2}}\left(\frac{z z^{\top}}{|z|^2} -\frac{1}{N+2}I_{N}\right )\rho_x(z)dz\\
  & =&  \frac{N(N+2)}{2} \int_{0}^{\infty} \rho_x(h) \int_{h \Ss^{N-1}} \frac{u(x+z)-2u(x)+u(x-z)}{|z|^{2}}\left(\frac{z z^{\top}}{|z|^2} -\frac{1}{N+2}I_{N}\right ) dz dh \\
  & \overset{\text{Thm.~\ref{thm:implexplcircle}}}{=} & \frac{N(N+2)}{2} \int_{0}^{\infty} \rho_x(h) 2 C_{1 2} h^{N-1} H_{u,h} d h\\
  & \overset{\text{Prop.~\ref{prop:cij}}}{=}&  N(N+2)\frac{2 \pi^{N/2} }{N ( N+2 ) \Gamma ( N/2 )} \int_{0}^{\infty} \rho_x(h) h^{N-1} H_{u,h} d h\\
  & = & \frac{2 \pi^{N/2} }{\Gamma ( N/2 )} \int_{0}^{\infty} h^{N-1} |h \Ss^{N-1}|^{-1} H_{u,h} \int_{h \Ss^{N-1}} \rho_x(z) d h\\    
  & = & \frac{2 \pi^{N/2} }{\Gamma ( N/2 )} |\Ss^{N-1}|^{-1} \int_{0}^{\infty}  H_{u,h} \int_{h \Ss^{N-1}} \rho_x(z)  d z d h\\
  & = & \int_{0}^{\infty}  H_{u,h} \int_{h \Ss^{N-1}} \rho_x(z)  d z d h.
\end{eqnarray*}
\end{proof}

Note that while Proposition~\ref{prop:implexplconn} requires radial symmetry of $\rho$, this symmetry was generally assumed throughout Section~\ref{analysis}. Therefore Section~\ref{analysis} can also be seen as providing localization results for \emph{implicit} models of the form \eqref{eq:implexplconn}.

\subsection{Choosing the weights for jump preservation}\label{sec:jumppreservingweights}
A characteristic of non-local models is that they are extremely flexible due to the many degrees of freedom in choosing the weights. In this work we will focus on improving on the question of how to reconstruct images that are piecewise quadratic but may have jumps. The issue here is that one wants to keep the Hessian sparse in order to favor piecewise affine functions, but doing it in a straightforward way, such as by adding $|D^2 u|(\Omega)$ as a regularizer, enforces too much first-order regularity \cite{Sch98a,LLT03,LT06,HS06}.

There have been several attempts of overcoming this issue, most notably approaches based on combined first- and higher-order functionals (see \cite{mineJMIV} and the references therein),
infimal convolution \cite{ChambolleLions}, and Total Generalized Variation \cite{TGV,setzer2011infimal}. Here we propose another strategy making use of the non-local formulation \eqref{eq:implicitnl}.

\begin{figure}[tp]
\bc
\includegraphics[height=.25\linewidth]{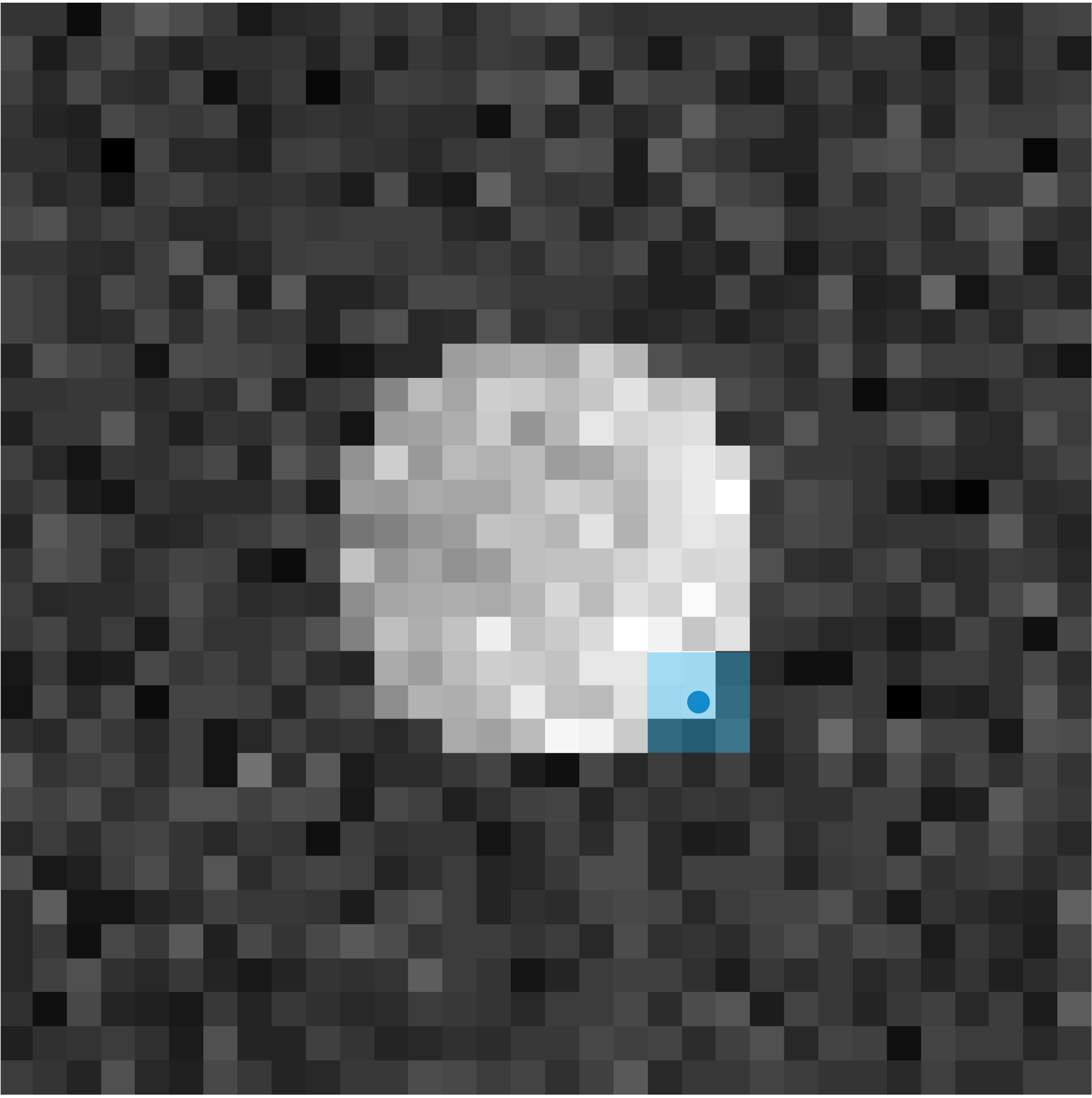}
\includegraphics[height=.25\linewidth]{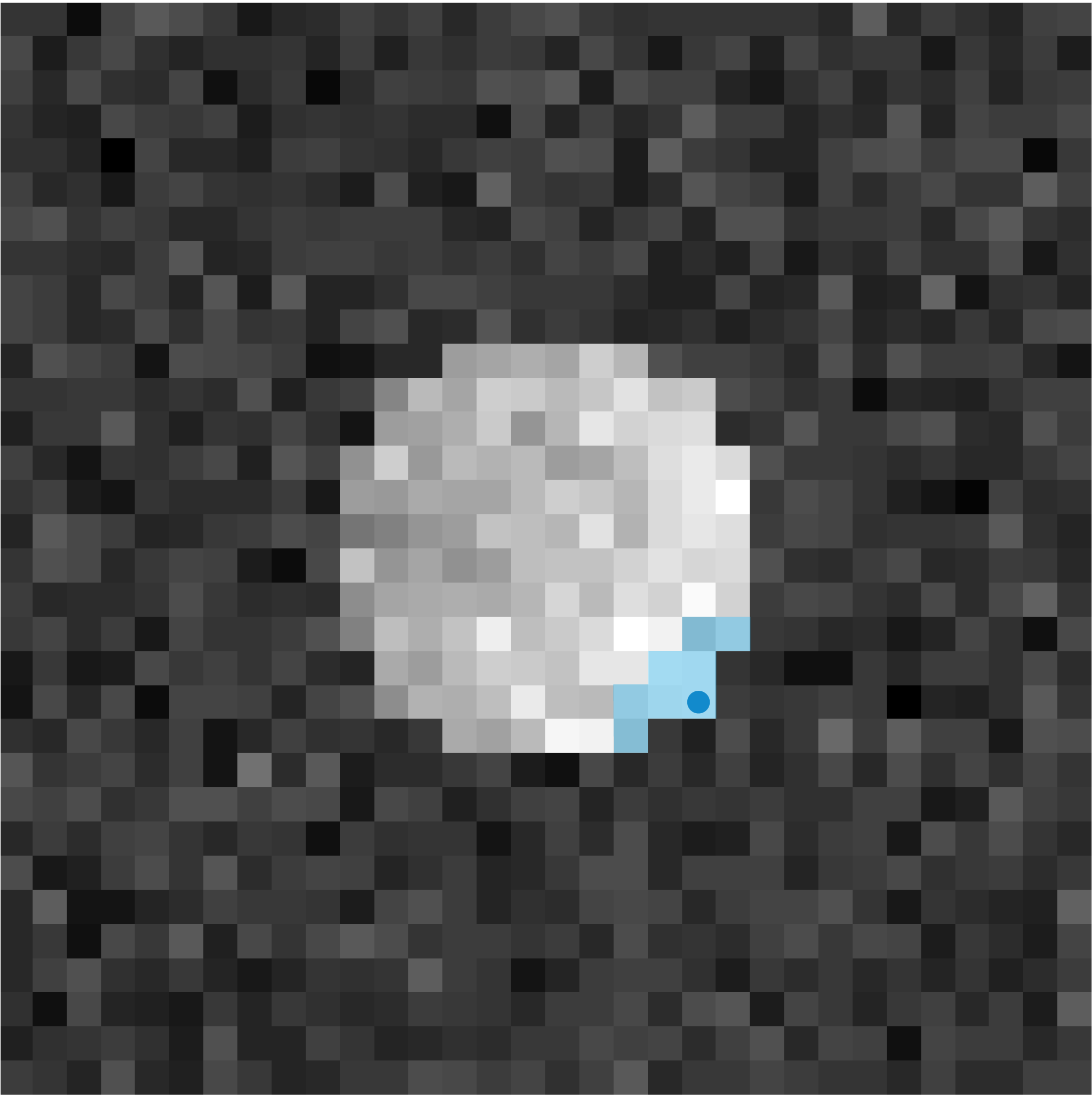}
\includegraphics[height=.25\linewidth]{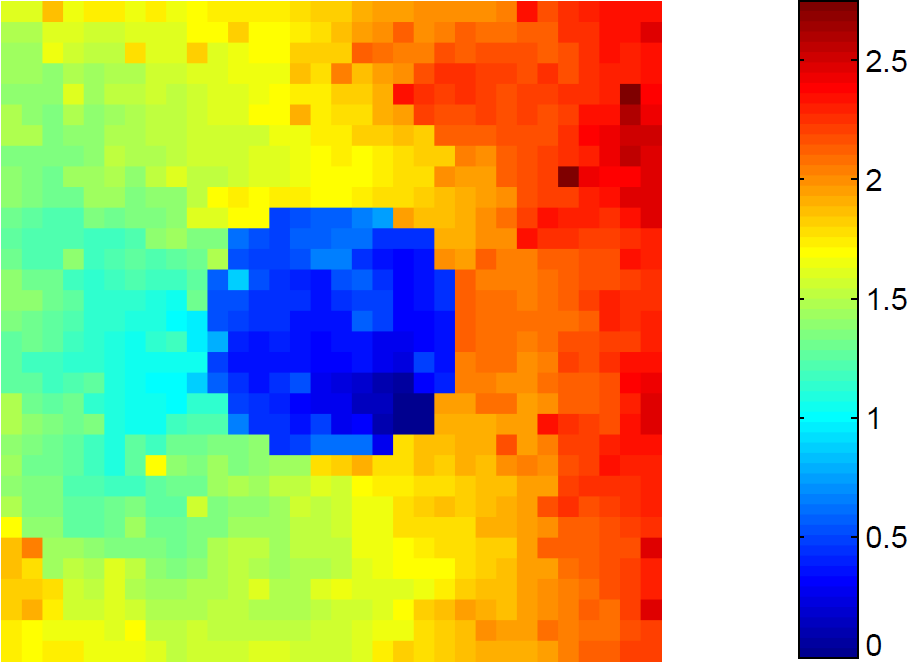}

\ec
\caption{Adaptive choice of the neighbourhood for the image of a disc with constant slope and added Gaussian noise.
\textbf{Left:} A standard discretization of the second-order derivatives at the marked point uses all points in a $3 \times 3$ neighbourhood.
\textbf{Center:} With a suitable choice of the weights $\sigma_x$, the discretization of the (non-local) Hessian at the marked point only (or mostly) involved points that are likely to belong to the same affine region, here the inside of the disc. This allows to use a straight-forward second-order regularizer such as the norm of the Hessian $\|H'_u\|$, while preserving jumps.
\textbf{Right:} Geodesic distance $d_\MM$ from the marked point $x$ in the lower right. Points that are separated from $x$ by a strong edge have a large geodesic distance to $x$, and are therefore not included in the neighbourhood of $x$ that is used to define the non-local Hessian at $x$.}\label{fig:neighbourhood}
\end{figure}

We draw our motivation for choosing the weights partly from a recent discussion of non-local ``Amoeba'' filters \cite{Lerallut2007,Welk2011,Welk2012}. Amoeba filters use classical techniques such as iterated median filtering, but the structuring element is chosen in a highly adaptive local way that can follow the image structures, instead of being restricted to a small parametrized set of shapes. In the following we propose to extend this idea to the higher-order energy minimization framework (Figure~\ref{fig:neighbourhood}).

Given a noisy image $g:\Om \to \R$, we first compute its (equally noisy) gradient image $\nabla g$ (in all of this section we consider only the discretized problem, so we can assume that the gradient exists). We then define the Riemannian manifold $\MM$ on the points of $\Omega$ with the usual Riemannian metric, weighted by $\varphi(\nabla g)$ for some function $\varphi: \RN \to \R_{+}$. In our experiments we used
\beqa
\varphi(g) \assign |g|^2 + \gamma
\eeqa
for small $\gamma > 0$, but other choices are equally possible. The choice of $\gamma$ controls how strongly the edge information is taken into account: for $\gamma = 0$ the weights are fully anisotropic, while for large $\gamma$ the magnitude of the gradient becomes irrelevant, and the method will reduce to an isotropic regularization.

With these definitions, the \emph{geodesic distance} $d_\MM$ between two points $x,y \in \Om$ now has a powerful interpretation: If $d_\MM(x,y)$ is large, this indicates that $x$ and $y$ are separated by a strong edge and should therefore not appear in the same regularization term. On the other hand, small $d_\MM(x)$ indicates that $x$ and $y$ are part of a more homogeneous region, and it is reasonable to assume that they should be part of the same affine part of the reconstructed image. 

For a given point $x \in \Om$, we sort the neighbours $y^1,y^2,\ldots$ in order of ascending distance, i.e., $d_\MM(x,y^1) \leq d_\MM(x,y^2) \leq \ldots$. We choose a neighbourhood size $M \in \N$ and set the weights to
\beqa
\sigma_x(y) \assign \begin{cases}1, & i \leq M, \\
0 & i > M. \\
\end{cases}
\eeqa
In other words, the non-local Hessian at $x$ is computed using its $M$ closest neighbours with respect to the geodesic distance through the gradient image.

The geodesic distances $\sigma_x(y)$ can be efficiently computed using the Fast Marching Method \cite{sethian1999level, fedkiw2003level} by solving the Eikonal equation
\beqa
|\nabla c(y)| & = & \varphi(\nabla g(y)),\\
c(x) & = & 0,
\eeqa
and setting $d_\MM(x,y) = c(x)$. Although it is necessary to process this step for every point $x \in \Om$, it is in practice a relatively cheap operation: the Fast Marching Method visits the neighbours of $x$ in the order of ascending distance $d_\MM$,
which means it can be stopped after $M$ points, with $M$ usually between $5$ and $20$.
If $M$ is chosen too small, one risks that the linear equation system that defines the non-local Hessian in \eqref{eq:implicitnl} becomes underdetermined.
In our experiments we found $M=12$ to be a good compromise, but the choice does not appear to be a very critical one.
In the experiments we used the $L^{1}$-non-local TV$^2$ model
\beqa
&\min_{u : \Om \to \R} &\sum_{x \in \Om} |u(x) - g(x)|^p d x + \alpha \sum_{x \in \Om} \omega (x) |H'_u(x)| d x\\
&\text{s.t.}& A H'_u = B u,
\eeqa
where  $\alpha > 0$ is the regularization strength, and $p \in \{1,2\}$. The linear constraints implement the optimality conditions for $H'_u$.

The local weight $\omega$ is set as $\omega(x) = M/|\{y\in\Om |B_{y,x} \neq 0 \}|$. While the approach does work with uniform weights $\omega = 1$, we found that in some cases it can erroneously leave single outlier points intact. We believe that this is caused by a subtle issue: by construction of $\sigma$, outlier points are usually close neighbour to fewer points. Therefore they appear in fewer of the regularization terms $|H'_u(x)|$, which effectively decreases regularization strength at outliers. The local weight $\omega$ counteracts this imbalance by dividing by the total number of terms that a particular value $u(x)$ appears in.

\subsection{Numerical results}\label{sec:experiments}
All experiments were performed on an Intel Xeon E5-2630 at 2.3 GHz with 64GB of RAM, MATLAB R2014a running on Scientific Linux 6.3, GCC 4.4.6, and Mosek 7. Run times were between several seconds for the geometric examples to several minutes for the full-size images, the majority of which was spent at the solution stage. The computation of the geodesic distances $d_\MM$ using the Fast Marching Method only took a few milliseconds in all cases, total preprocessing time including building the sparse matrix structures $A$ and $B$ took less than $5$ seconds for the full-size images.

The solution of the Eikonal equation and computation of the weights as well as system matrix uses a custom C++ implementation. For solving the minimization problems we used the commercial Mosek solver with the CVX interface. Compared to the recent first-order approaches this allows to compute solutions with very high accuracy and therefore to evaluate the model without the risk of accidentally comparing only approximate solutions.

\begin{figure}[tp]
\bc
\begin{tabular}{cccc}
\includegraphics[width=.2\linewidth]{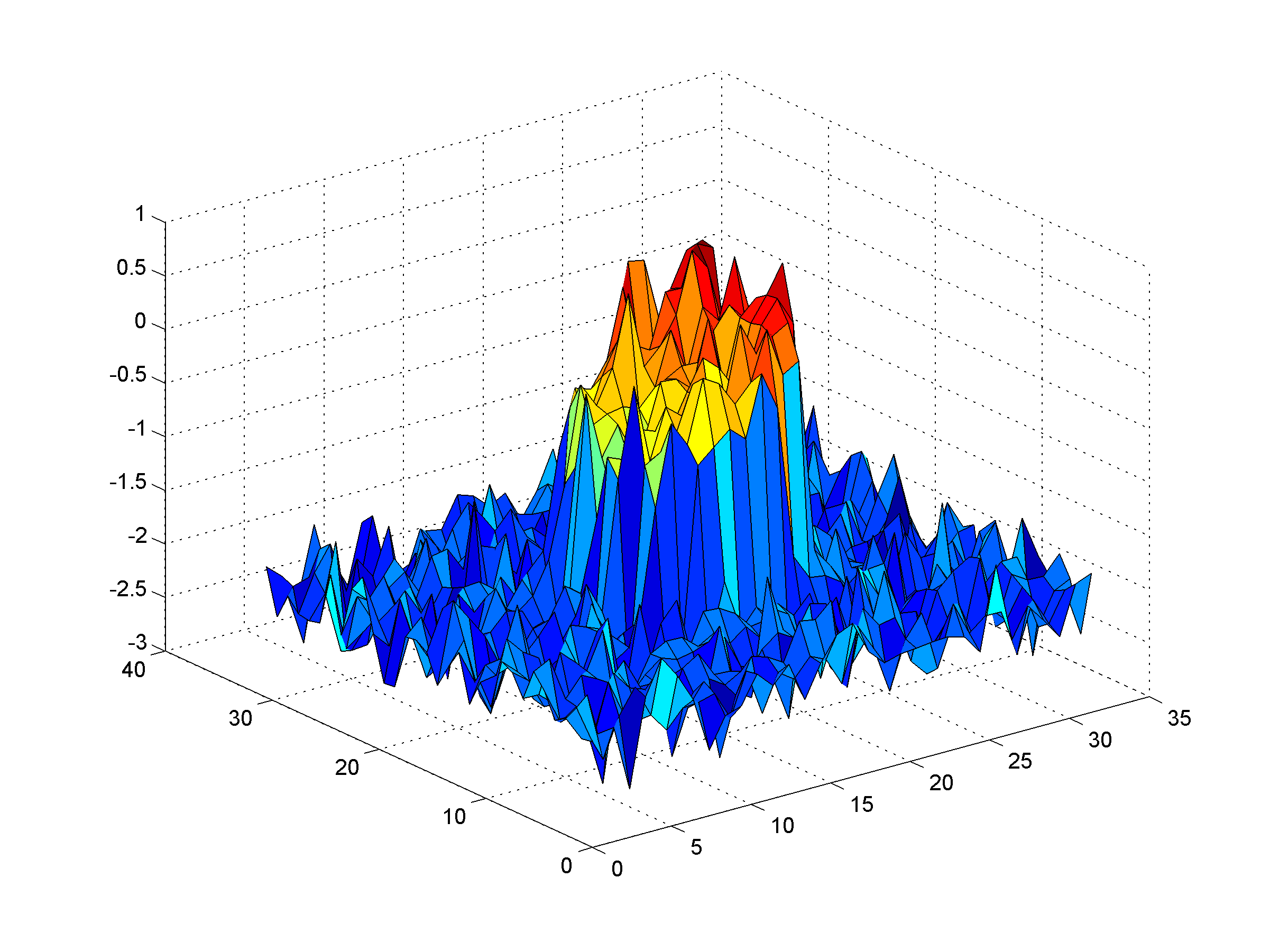} &
\includegraphics[width=.2\linewidth]{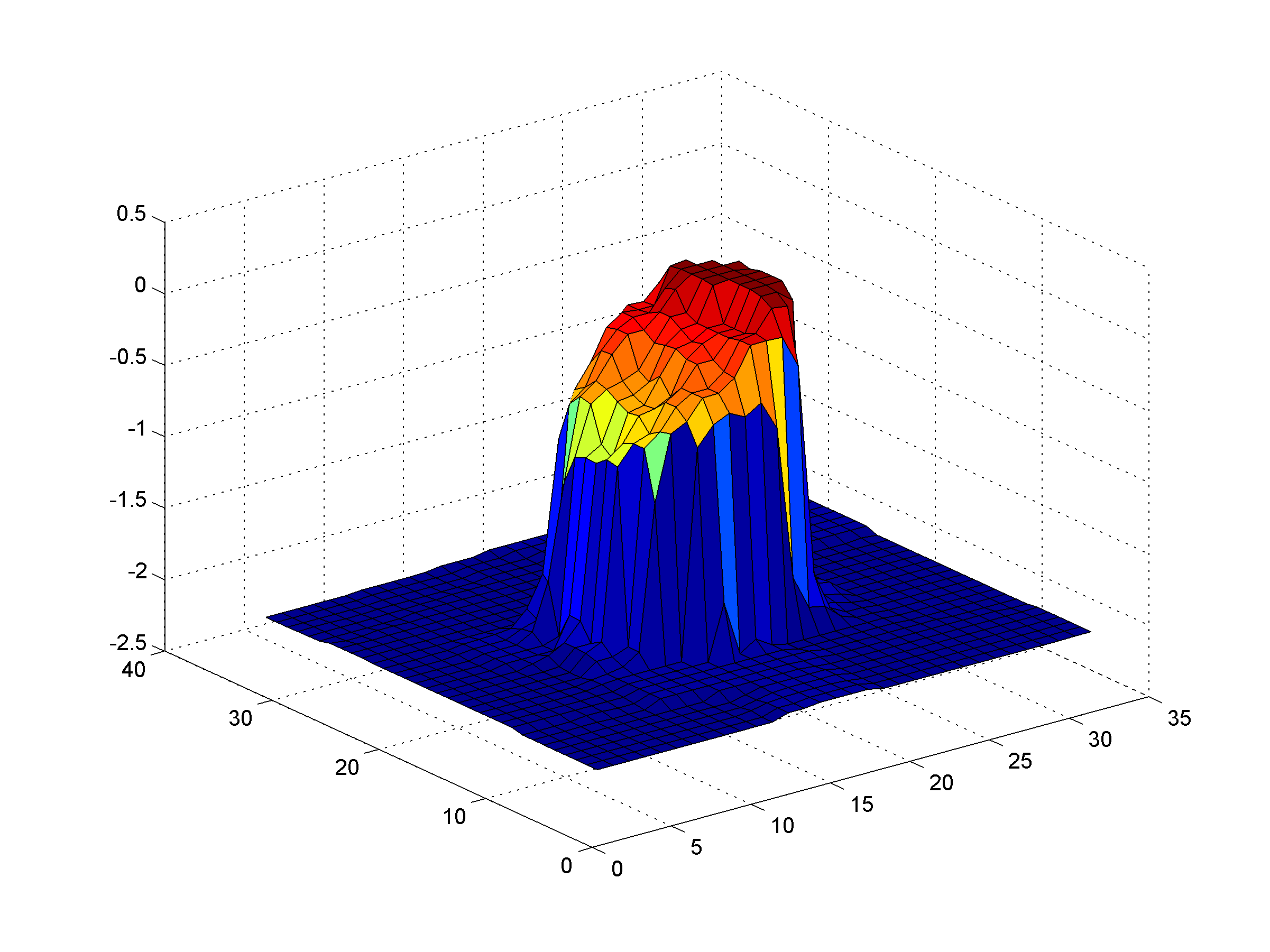} &
\includegraphics[width=.2\linewidth]{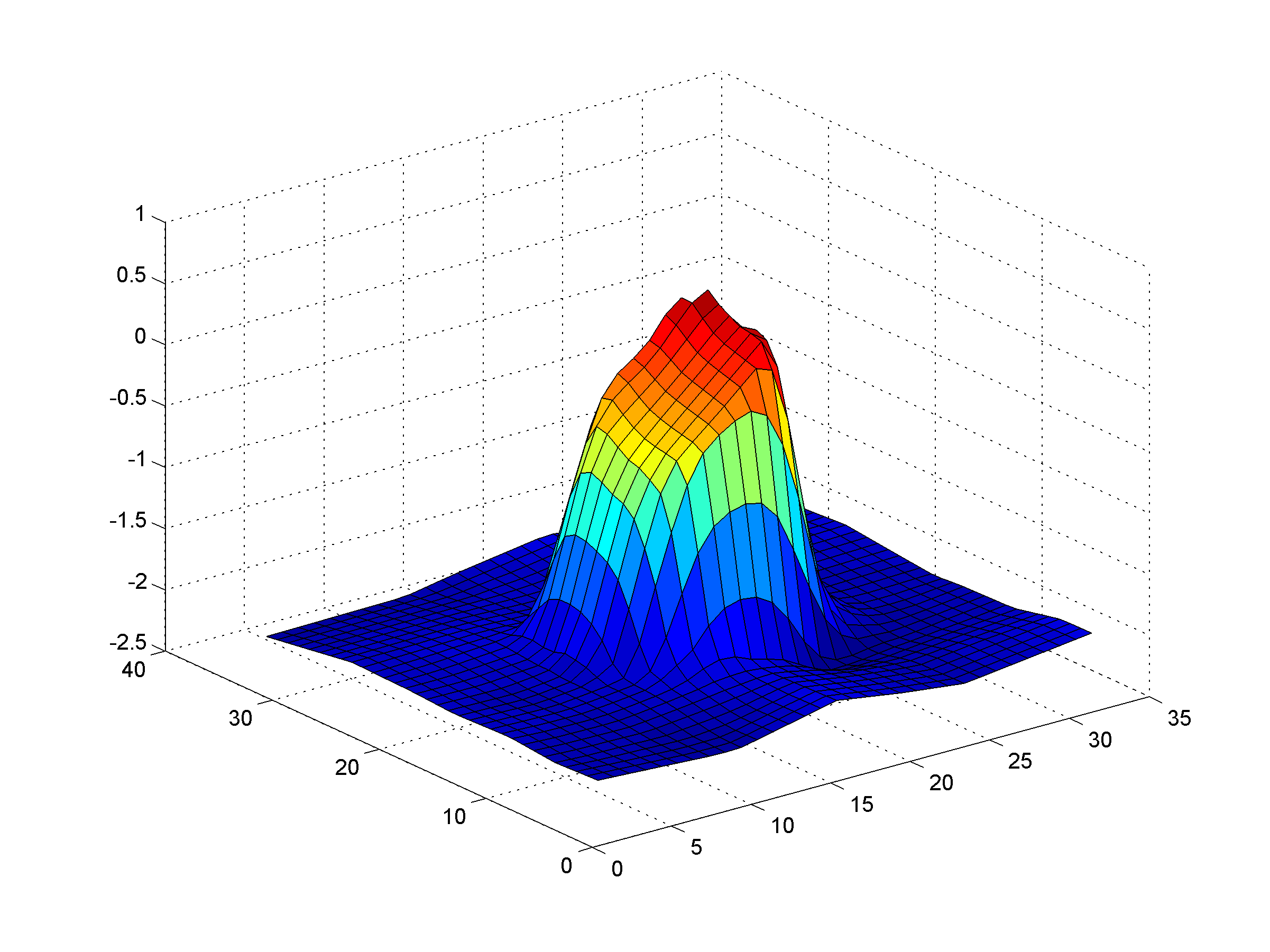} &
\includegraphics[width=.2\linewidth]{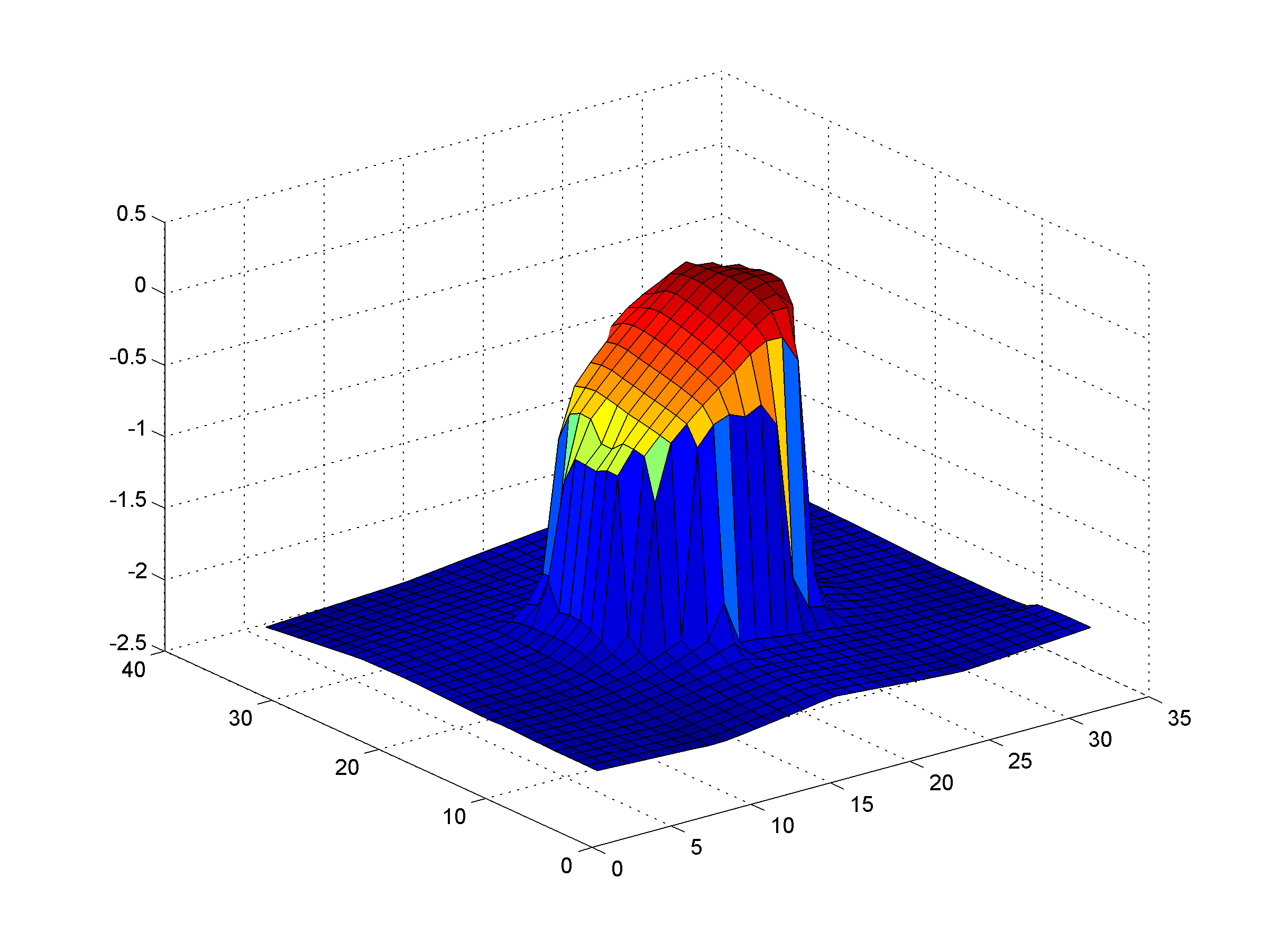}\\
input & $\TV$ & $\TV^2$ & TGV
\end{tabular}
\ec
\caption{Classical local regularization. The input consists of a disc-shaped slope with additive Gaussian noise, $\sigma = 0.25$. Shown is the result of denoising the input with $L^1$-data term. Total variation ($\TV, \alpha = 1.1$) regularization generates the well-known staircasing effect. Regularization of the Hessian ($\TV^2, \alpha=0.8$) avoids this problem, at the cost of over-smoothing jumps. Total Generalized Variation (TGV, $\alpha_0 = 1.5, \alpha_1 = 1.1$) performs best, but still clips the slope at the top.}\label{fig:disc-local}
\end{figure}

Figure~\ref{fig:disc-local} illustrates the effect of several classical local regularizers including $\TV$, $\TV^2$, and TGV. As expected,
$\TV$ generates the well-known staircasing effect, while $\TV^2$ leads to oversmoothing of the jumps. TGV with hand-tuned parameters performs reasonably well, however it exhibits a characteristic pattern of clipping sharp local extrema. This behavior has also been analytically confirmed in \cite{papafitsoros2013study, poschl2013exact} for the one-dimensional case.

\begin{figure}[tp]
\bc
\begin{tabular}{ccc}
\includegraphics[width=.2\linewidth]{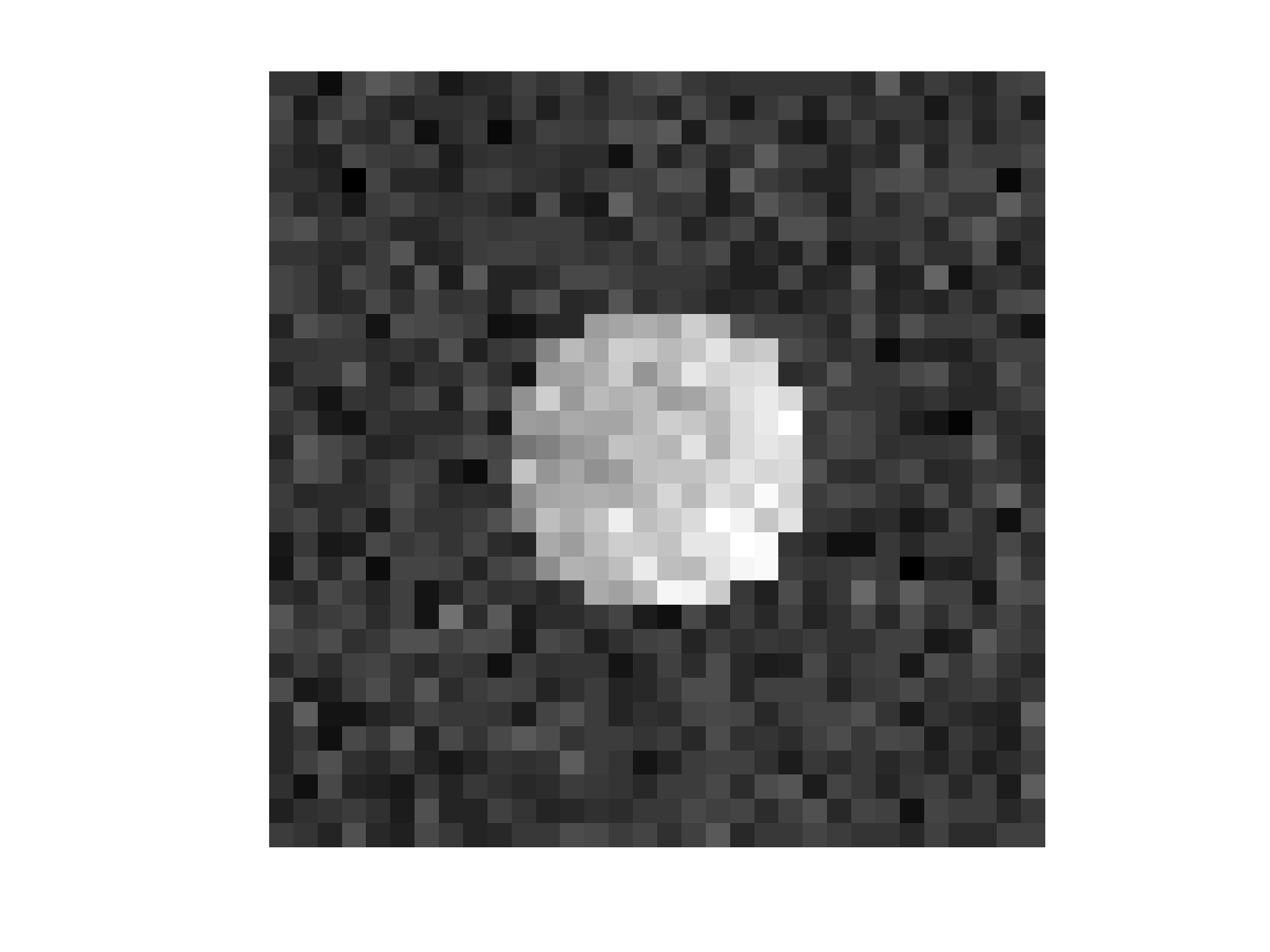} &
\includegraphics[width=.2\linewidth]{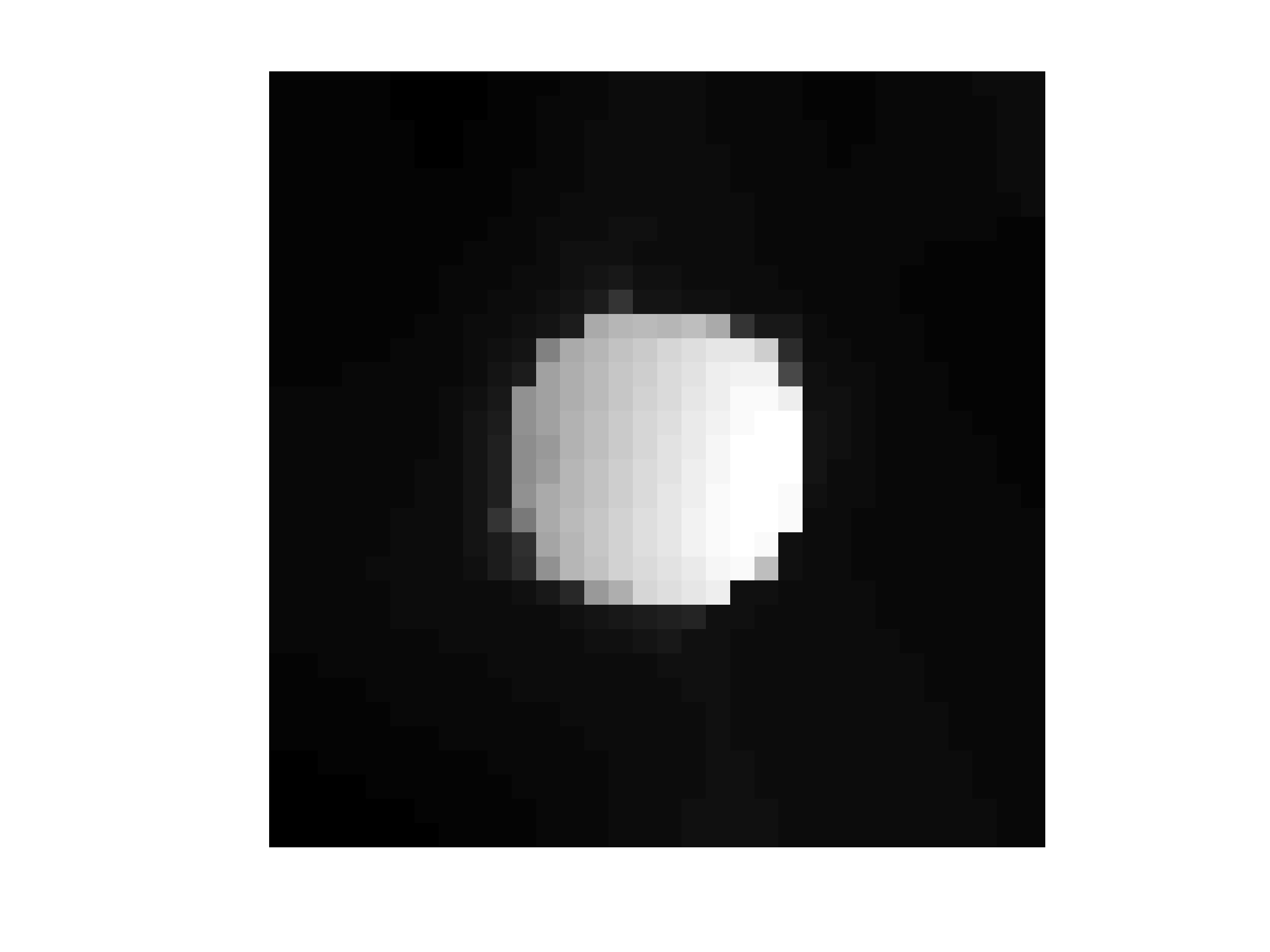} &
\includegraphics[width=.2\linewidth]{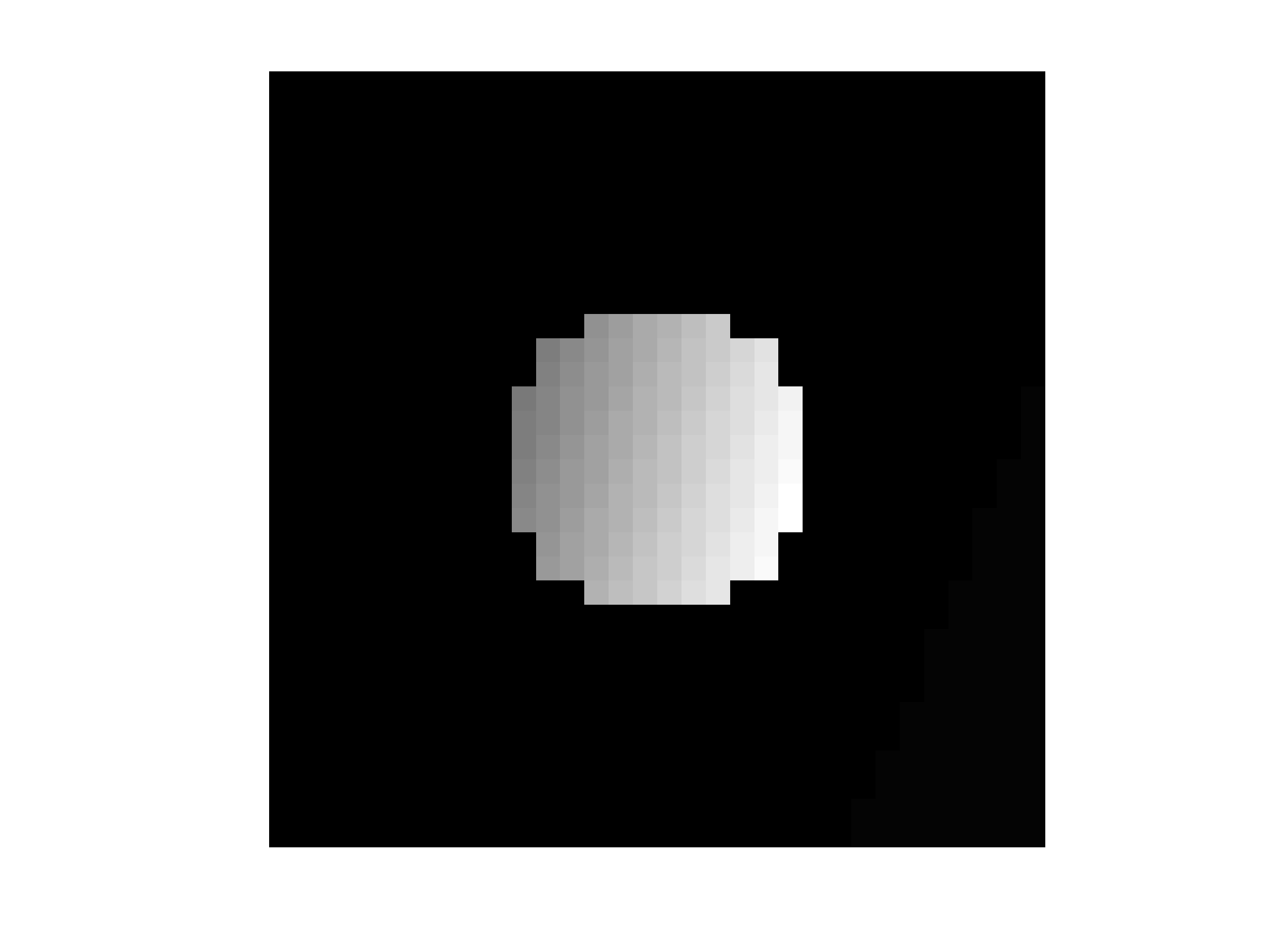}\\
\includegraphics[width=.2\linewidth]{input-surf.png} &
\includegraphics[width=.2\linewidth]{result-surftgv.png} &
\includegraphics[width=.2\linewidth]{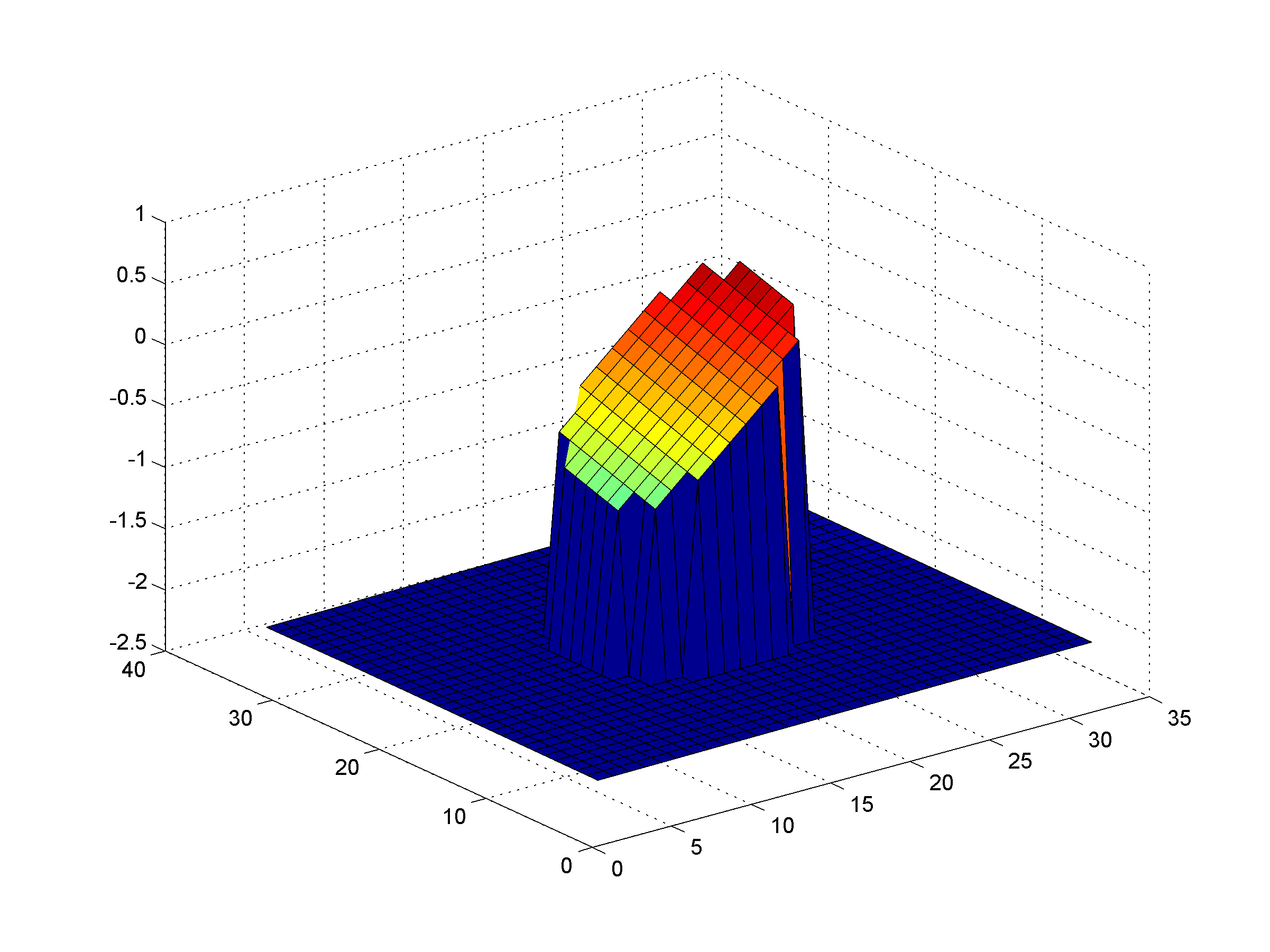}\\
input & TGV & non-local Hessian \\
& & (proposed)
\end{tabular}
\ec
\caption{Non-local regularization of the problem in Figure~\ref{fig:disc-local}. The adaptive choice of the neighbourhood and weights together with the non-local Hessian preserves the jumps, clip the top of the slope, and allows to perfectly reconstruct the piecewise affine signal.}\label{fig:disc-non-local}
\end{figure}

Figure~\ref{fig:disc-non-local} shows the effect of our non-local Hessian-based method for the same problem. The piecewise affine signal is almost perfectly recovered. This is mostly due to the fact that the jumps are relatively large, which means that after computing the neighborhoods the circular region and the background are not coupled in terms of regularization. Therefore the regularization weight $\alpha$ can be chosen very large, which results in virtually affine regions.

\begin{figure}[tp]
\bc
\begin{tabular}{cccc}
\includegraphics[width=.2\linewidth]{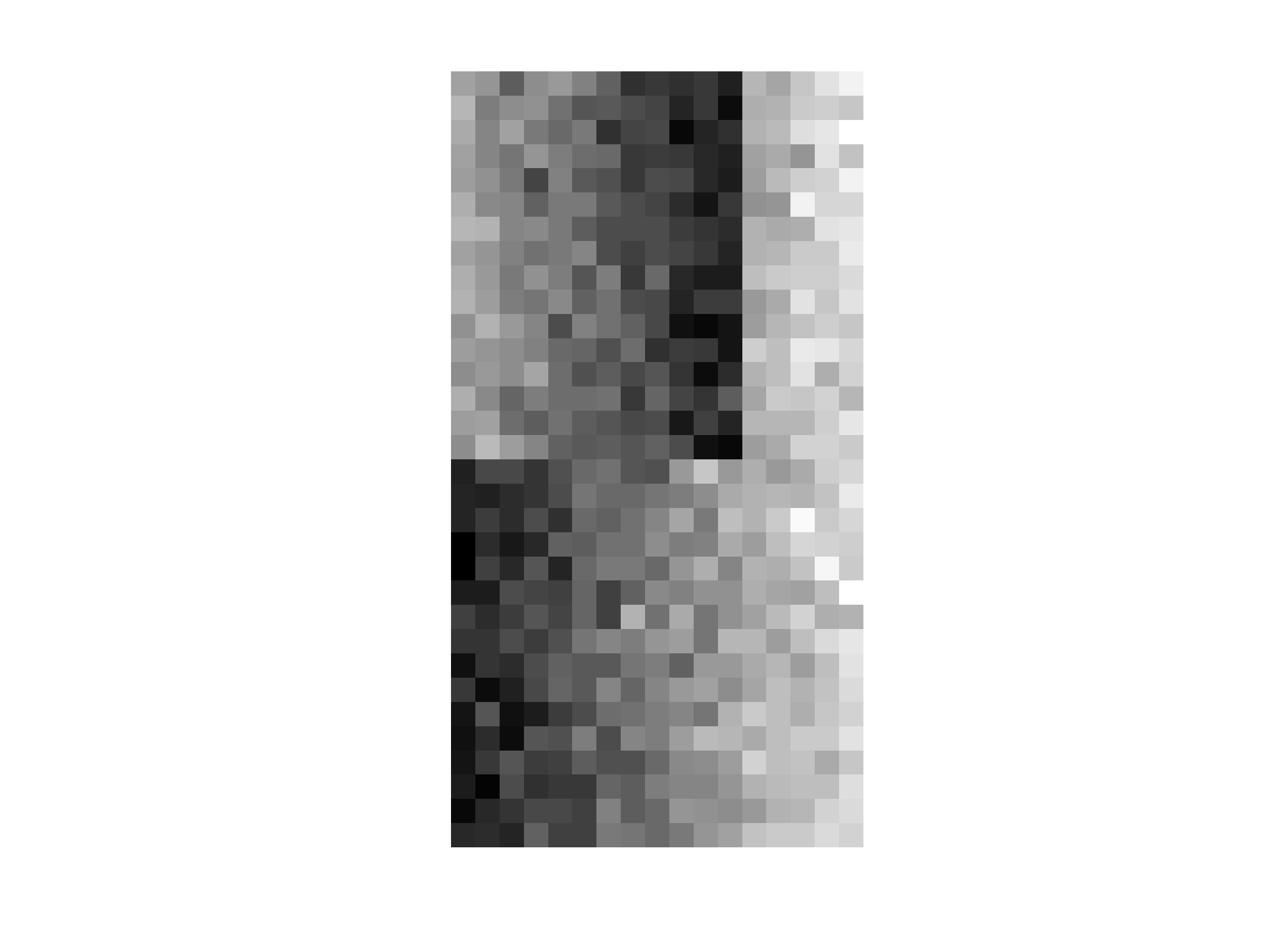} &
\includegraphics[width=.2\linewidth]{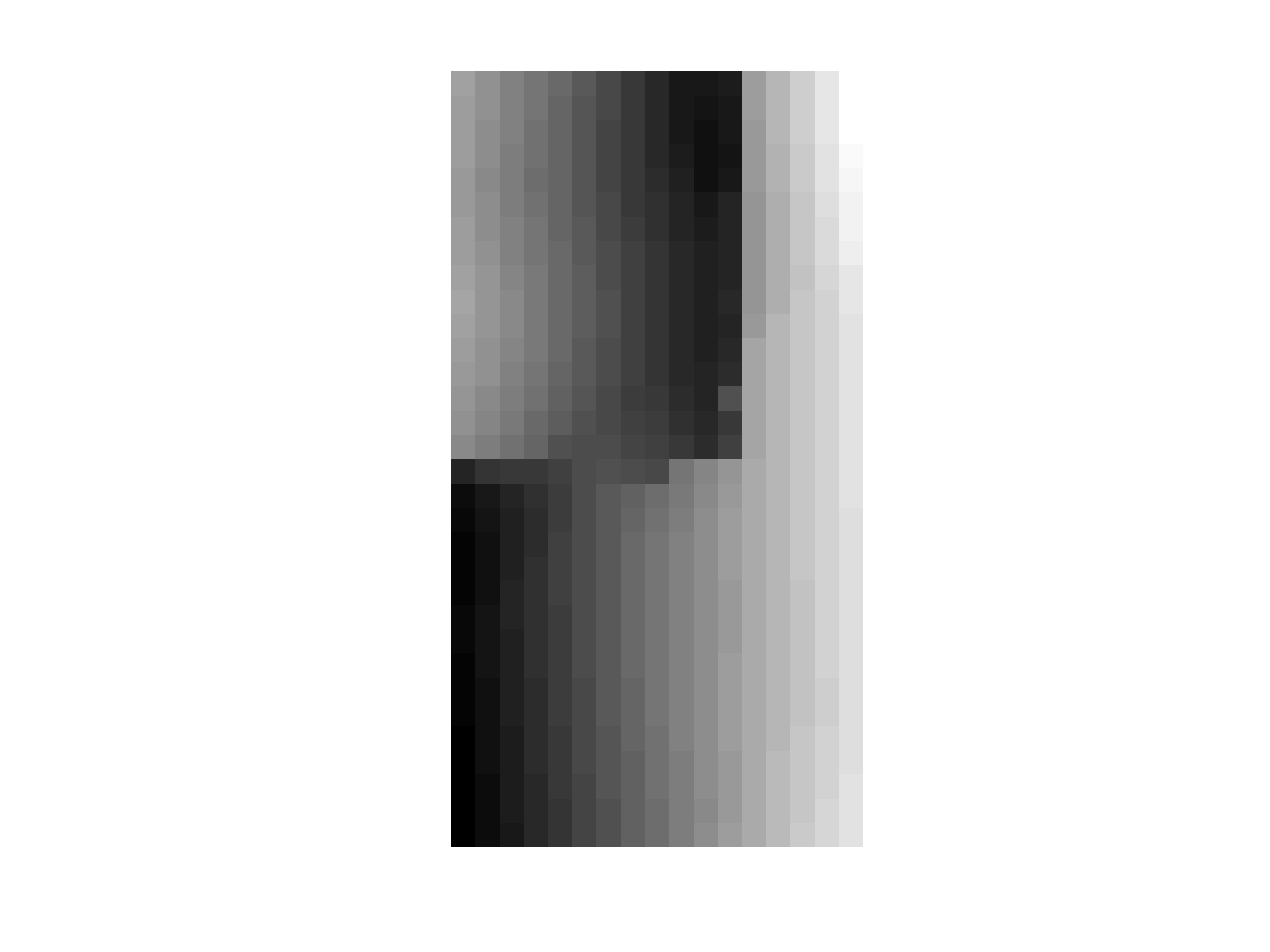} &
\includegraphics[width=.2\linewidth]{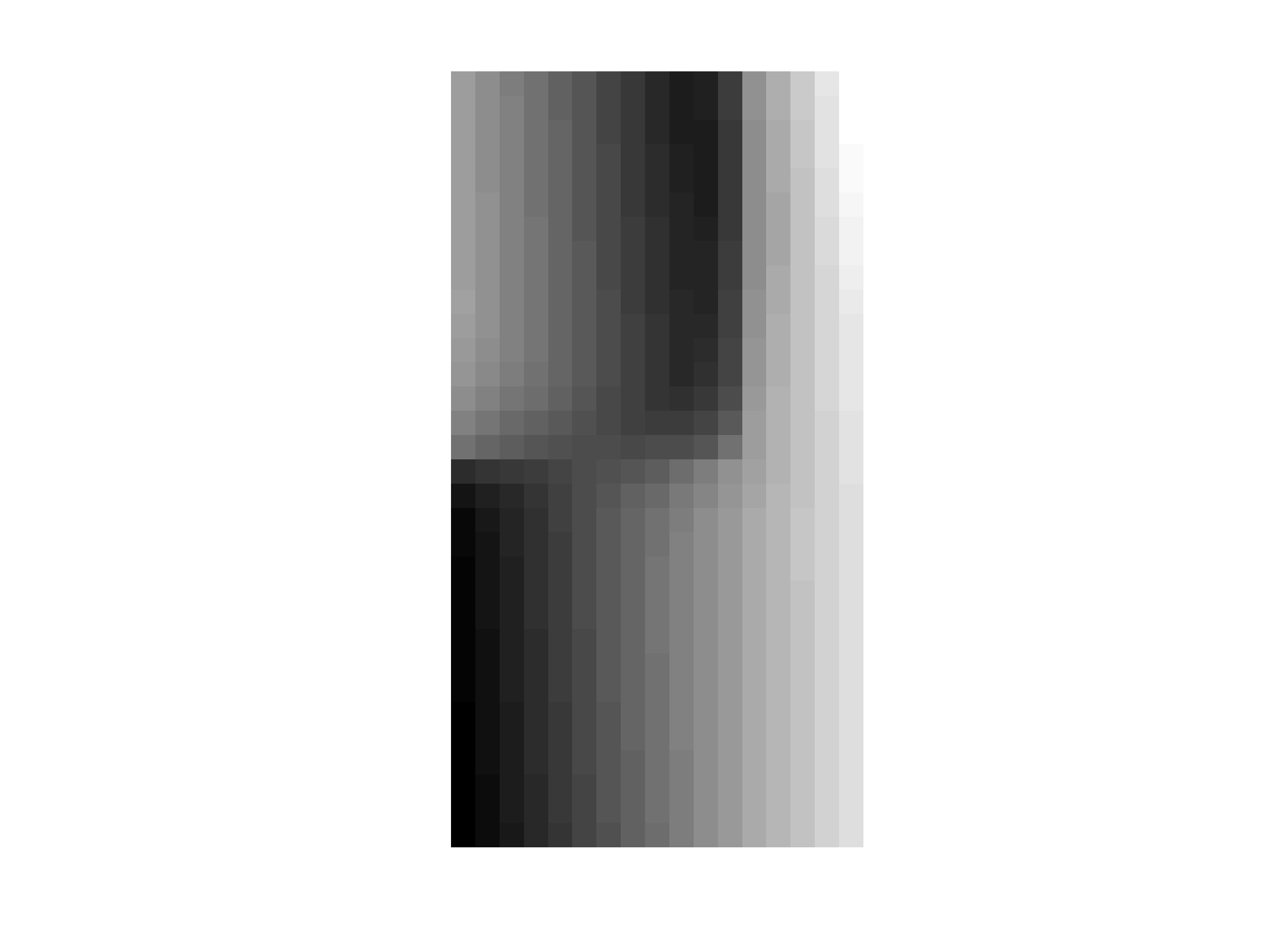} &
\includegraphics[width=.2\linewidth]{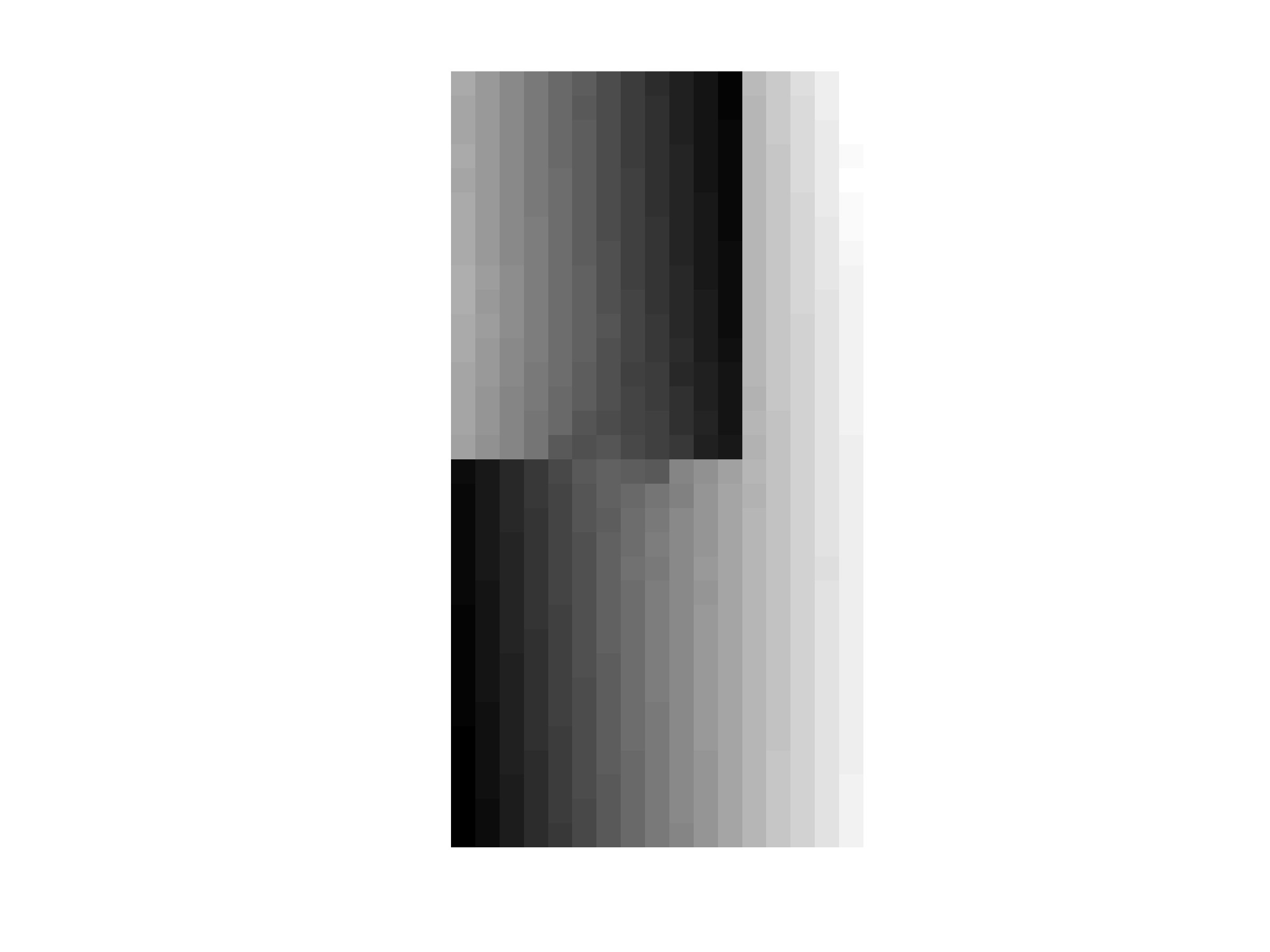}\\
\includegraphics[width=.2\linewidth]{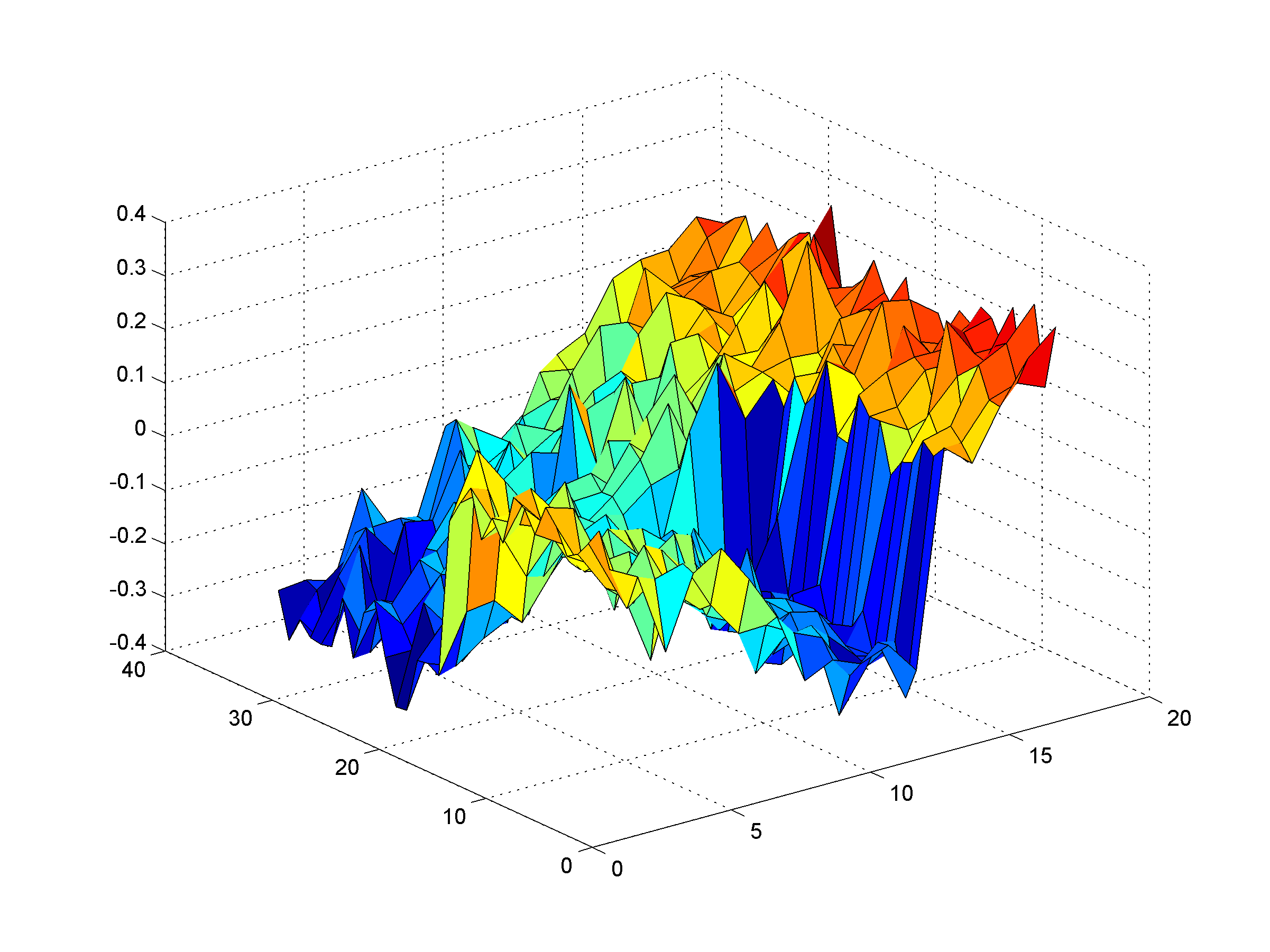} &
\includegraphics[width=.2\linewidth]{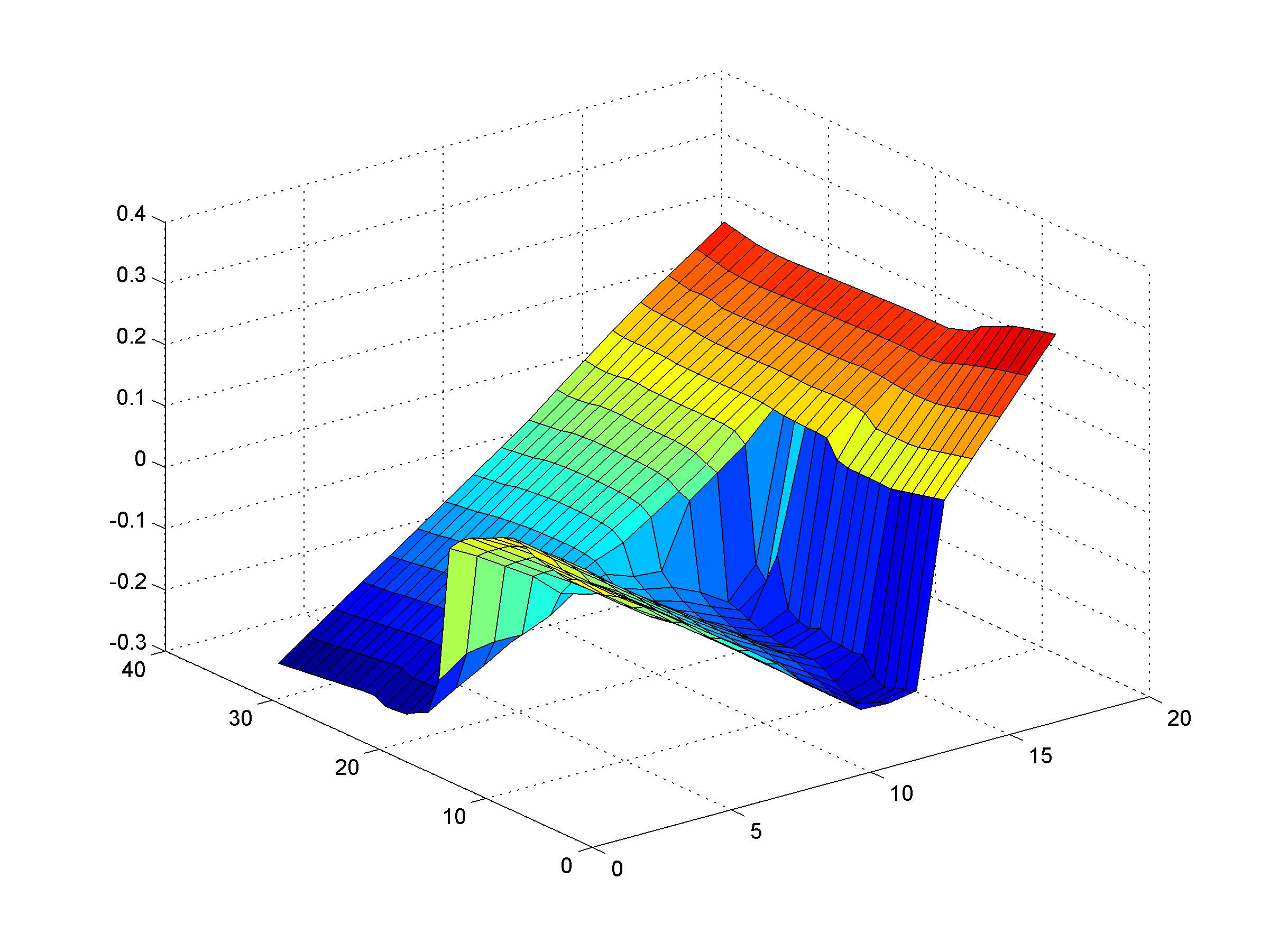} &
\includegraphics[width=.2\linewidth]{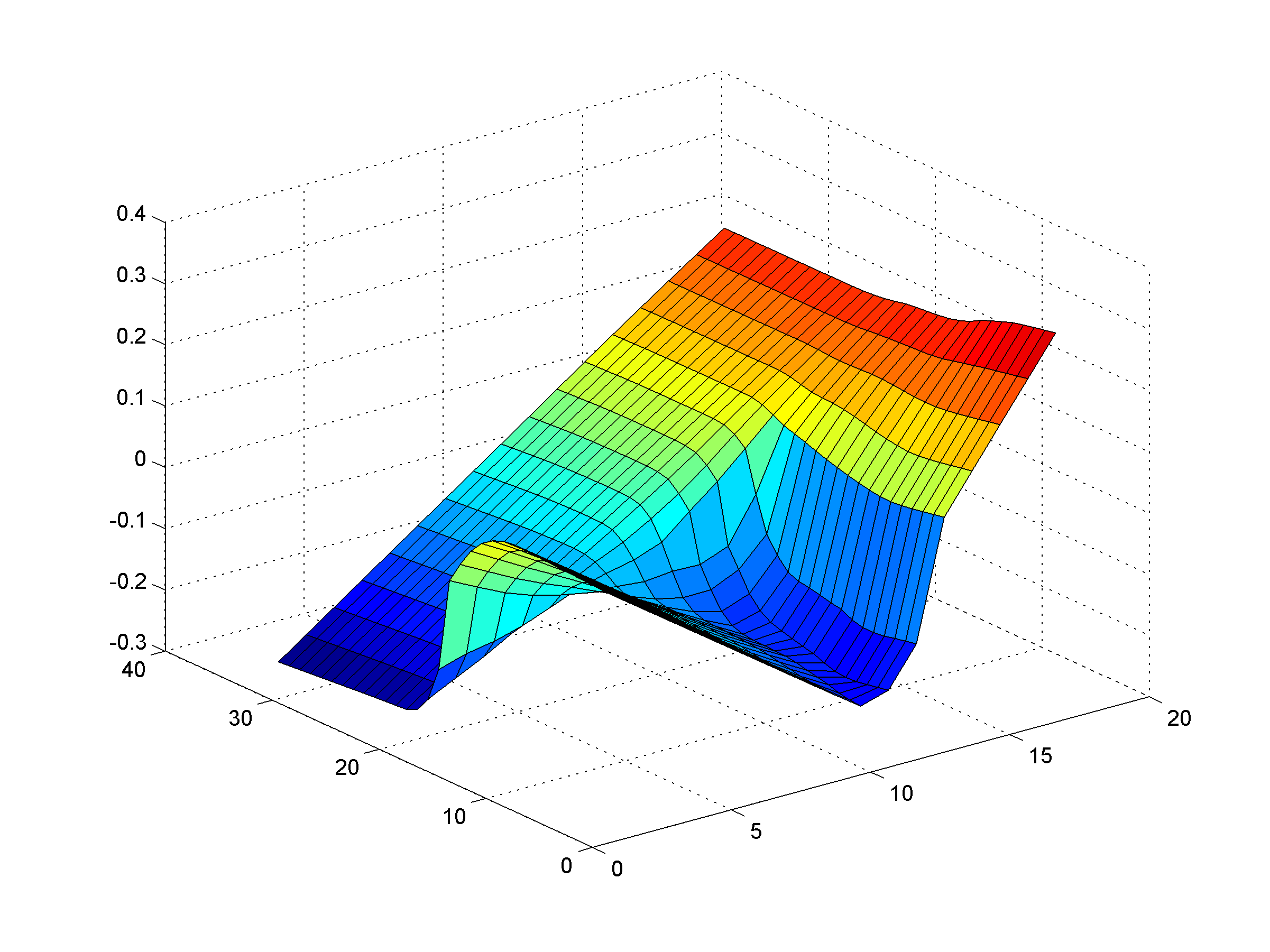} &
\includegraphics[width=.2\linewidth]{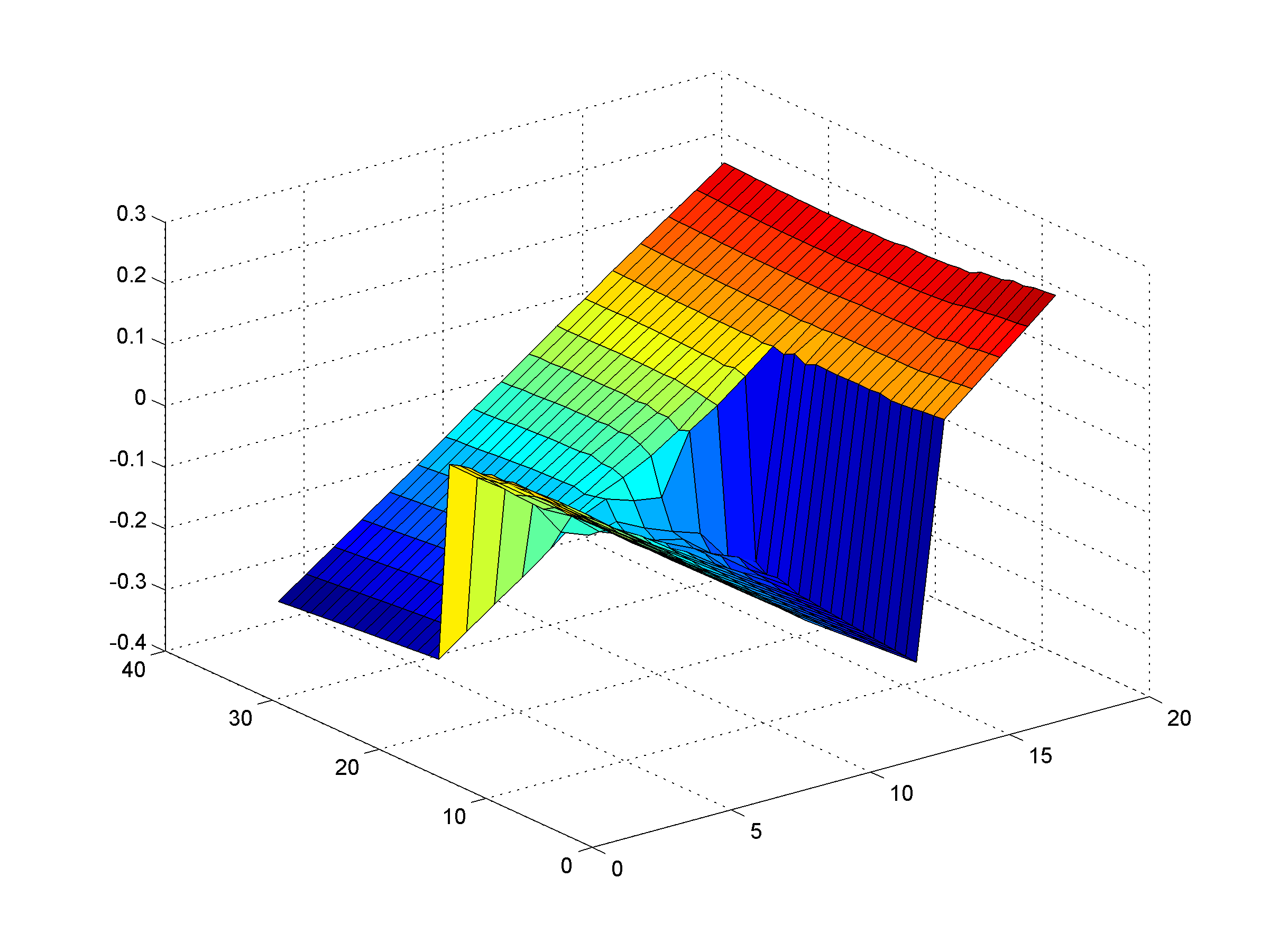}\\
input & TGV low & TGV high & non-local Hessian \\
& & &
\end{tabular}
\ec
\caption{Denoising results for the ``opposing slopes'' image. The small jump at the crossing causes a slight amount of smoothing for both the TGV and non-local approaches. TGV with low regularization strength (TGV low, $\alpha_0 = \alpha_1 = 0.8$) does reasonably well at preserving the strong jump on the right, but does not reconstruct the constant slope well. If one increases the regularization strength (TGV high, $\alpha_0=\alpha_1=1.5$), the jump is smoothed out. The non-local Hessian regularization ($\alpha=10^{-3}$, note the different scales on the axes) fully preserves the large jump by design, and results in an overall cleaner reconstruction. The small amplitude differences at the crossing point cause a slight blur in both approaches.}\label{fig:slopes-non-local}
\end{figure}

To see what happens with smaller jumps, we generated a pattern of opposing slopes (Figure~\ref{fig:slopes-non-local}). As expected both TGV as well as the non-local Hessian approach struggle when the jump is small. This shows the limitations of our approach for choosing the weights -- while it adds some structural information to the regularizer, this information is still restricted to a certain neighbourhood of each point, and does not take into account the full global structure.

\begin{figure}[tp]
\bc
\begin{tabular}{ccc}
\includegraphics[height=.2\linewidth]{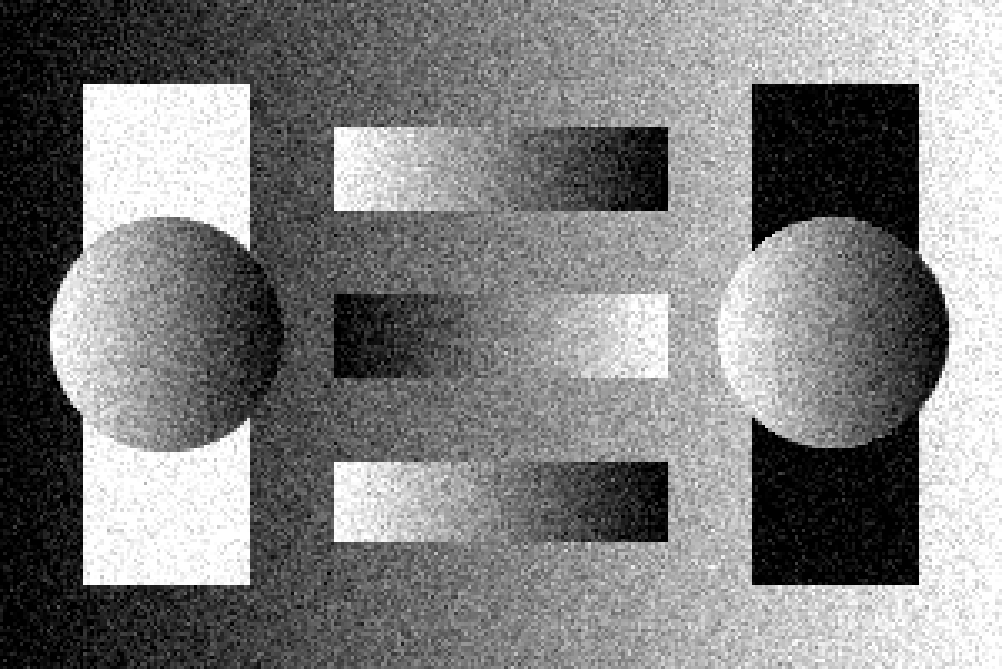} &
\includegraphics[height=.2\linewidth]{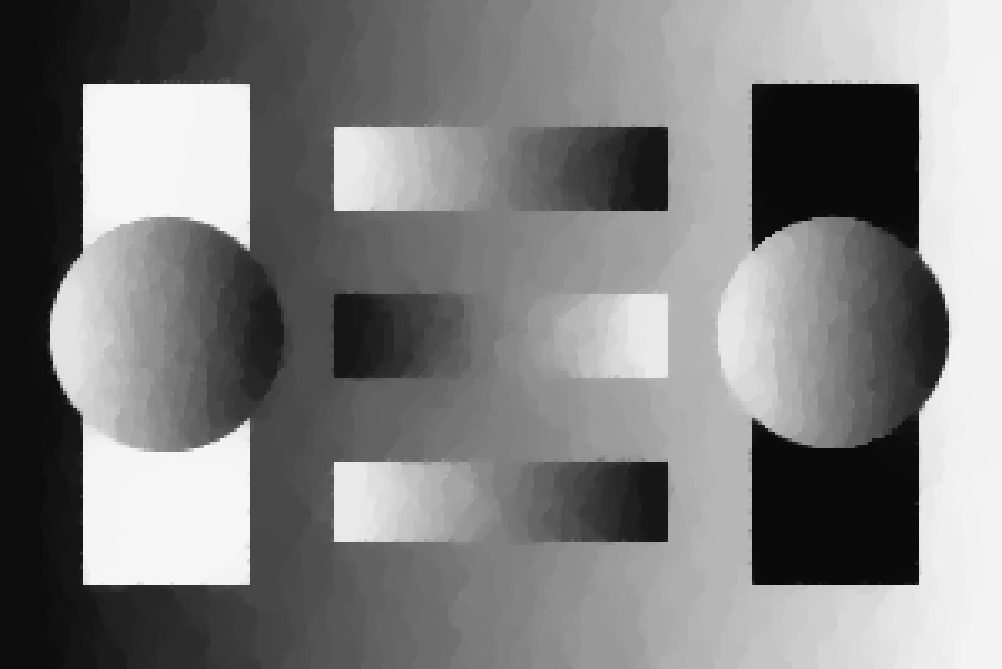} &
\includegraphics[height=.2\linewidth]{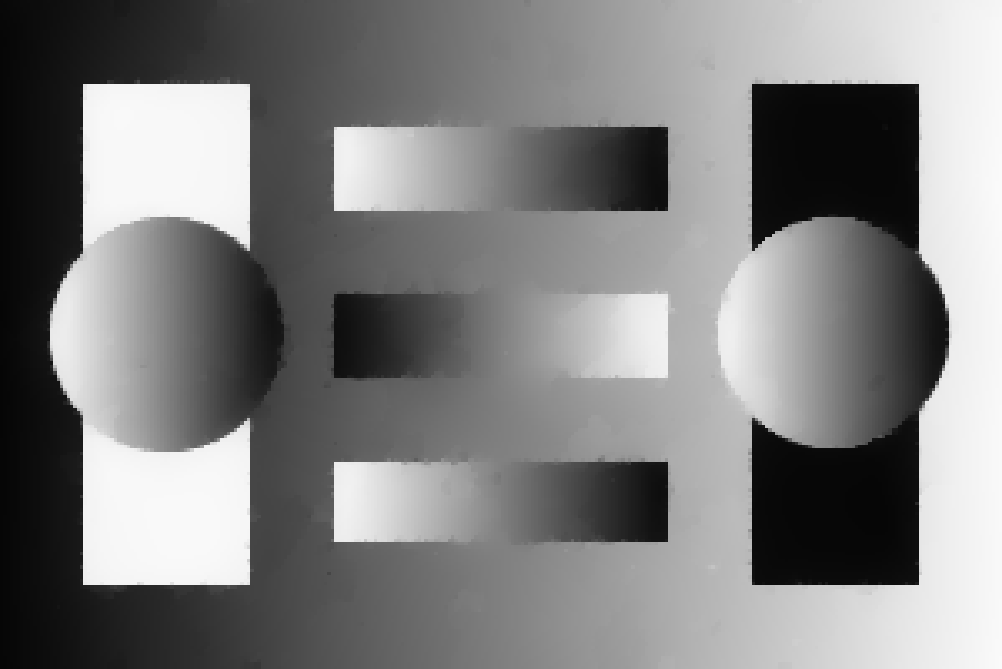}\\
input image &
$L^2-\TV$ (ROF, $\alpha = 0.24$) &
$L^2-\TGV (\alpha_0=0.38, \alpha_1=0.16)$\\
PSNR = $22.82$, SSIM = $0.33$ &
PSNR = $31.83$, SSIM = $0.90$ &
PSNR = $33.28$, SSIM = $0.93$\\
\includegraphics[height=.2\linewidth]{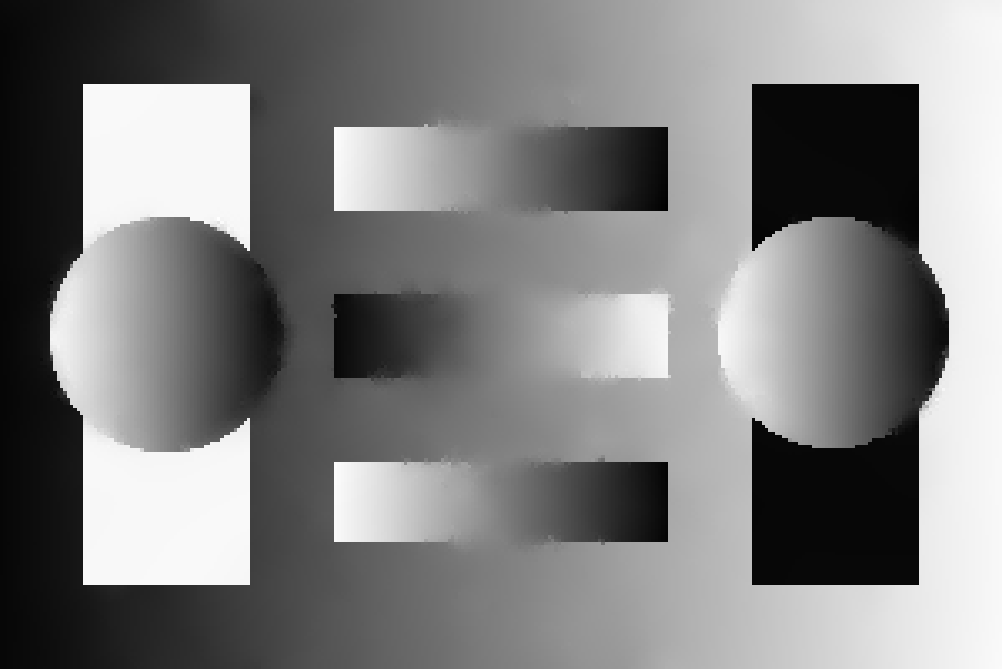} &
\includegraphics[height=.2\linewidth]{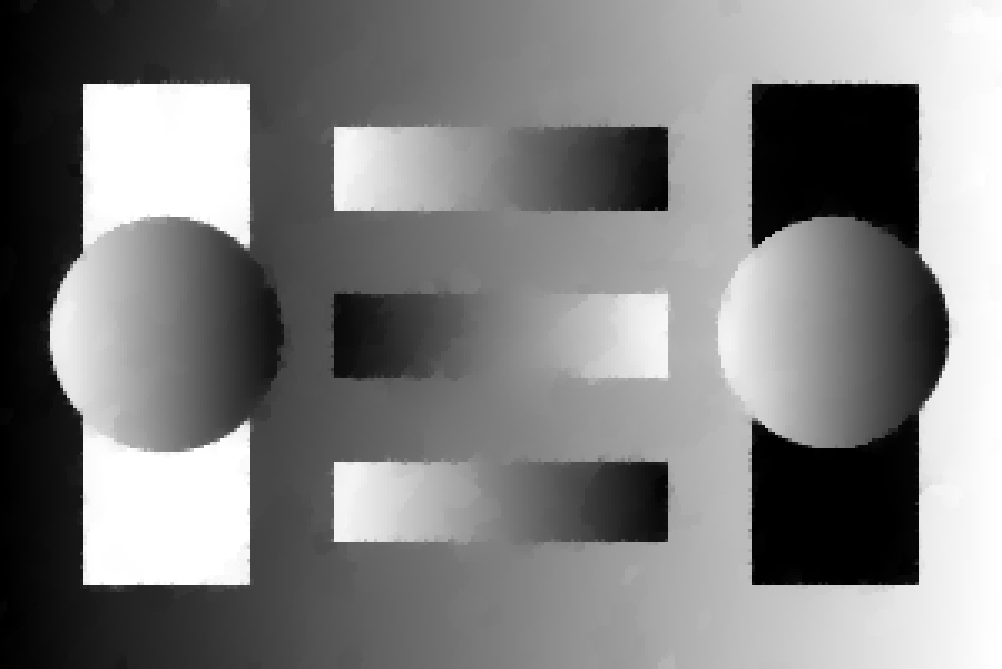} &
\includegraphics[height=.2\linewidth]{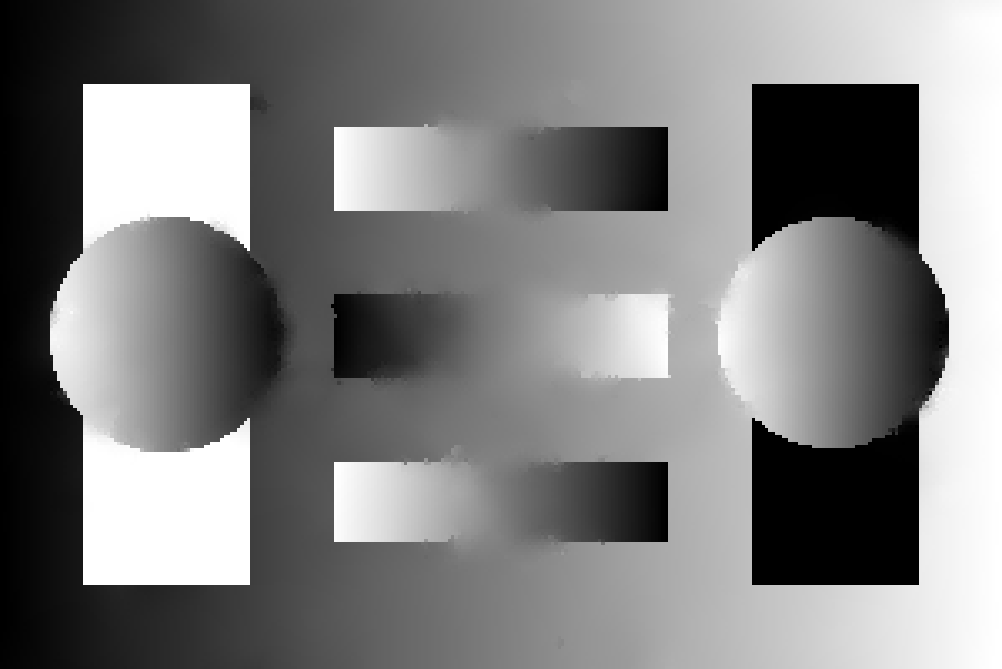}\\
$L^2-$non-local Hessian & $L^1-\TGV (\alpha_0=2.95, \alpha_1=1.0)$ & $L^1-$ NL Hessian ($\gamma = 0.005, \alpha = 1.65$) \\
 PSNR = $34.09$, SSIM = $0.93$ &
 PSNR = $35.39$, SSIM = $0.97$ &
 PSNR = $36.06$, SSIM = $0.98$\\
\end{tabular}
\ec
\caption{Denoising results for a geometric test image (see Figure~\ref{fig:shapes-cut} for the discussion).}\label{fig:shapes-img}
\end{figure}

\begin{figure}[tp]
\bc
\begin{tabular}{ccc}
\includegraphics[height=.2\linewidth]{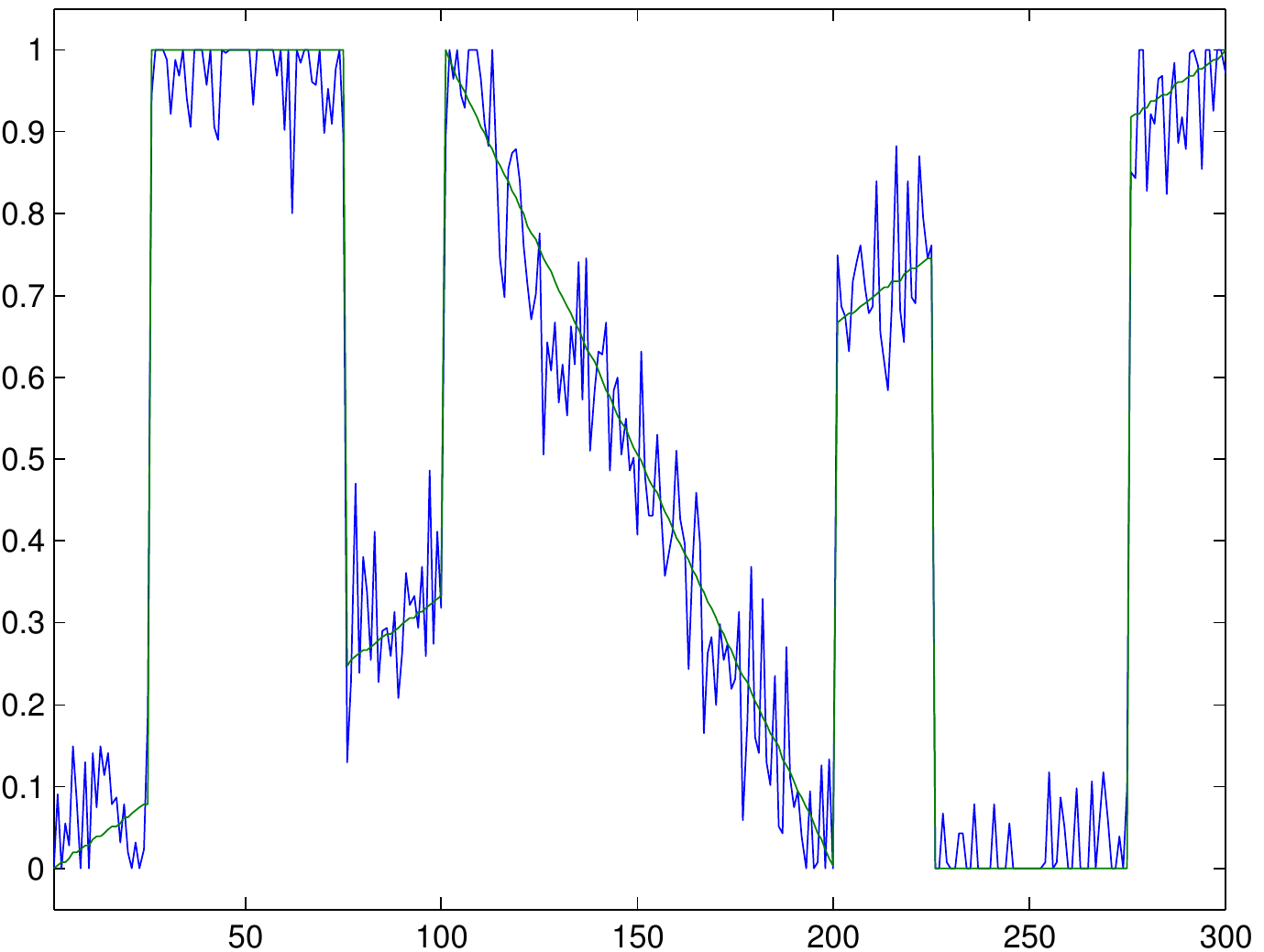} &
\includegraphics[height=.2\linewidth]{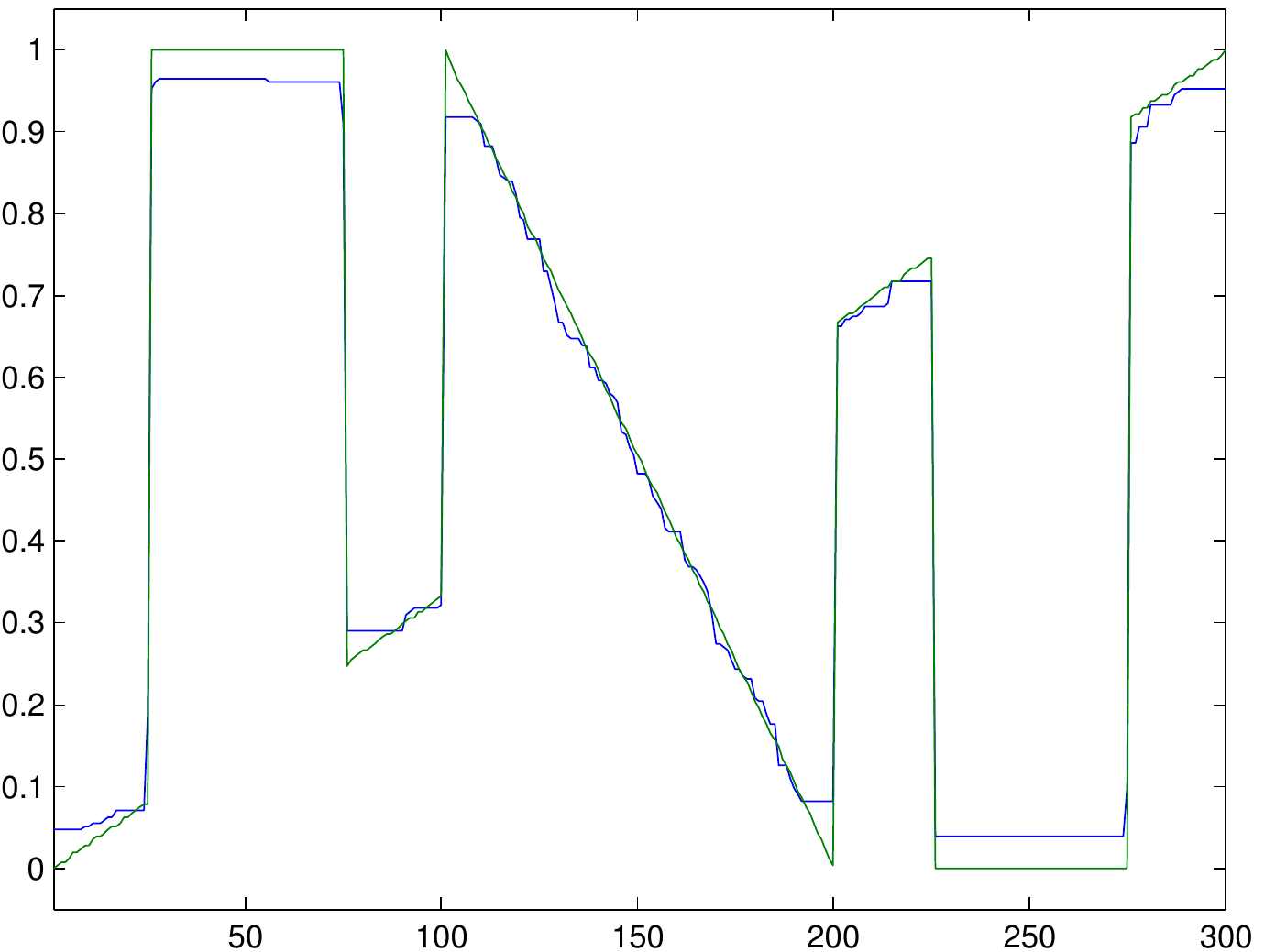} &
\includegraphics[height=.2\linewidth]{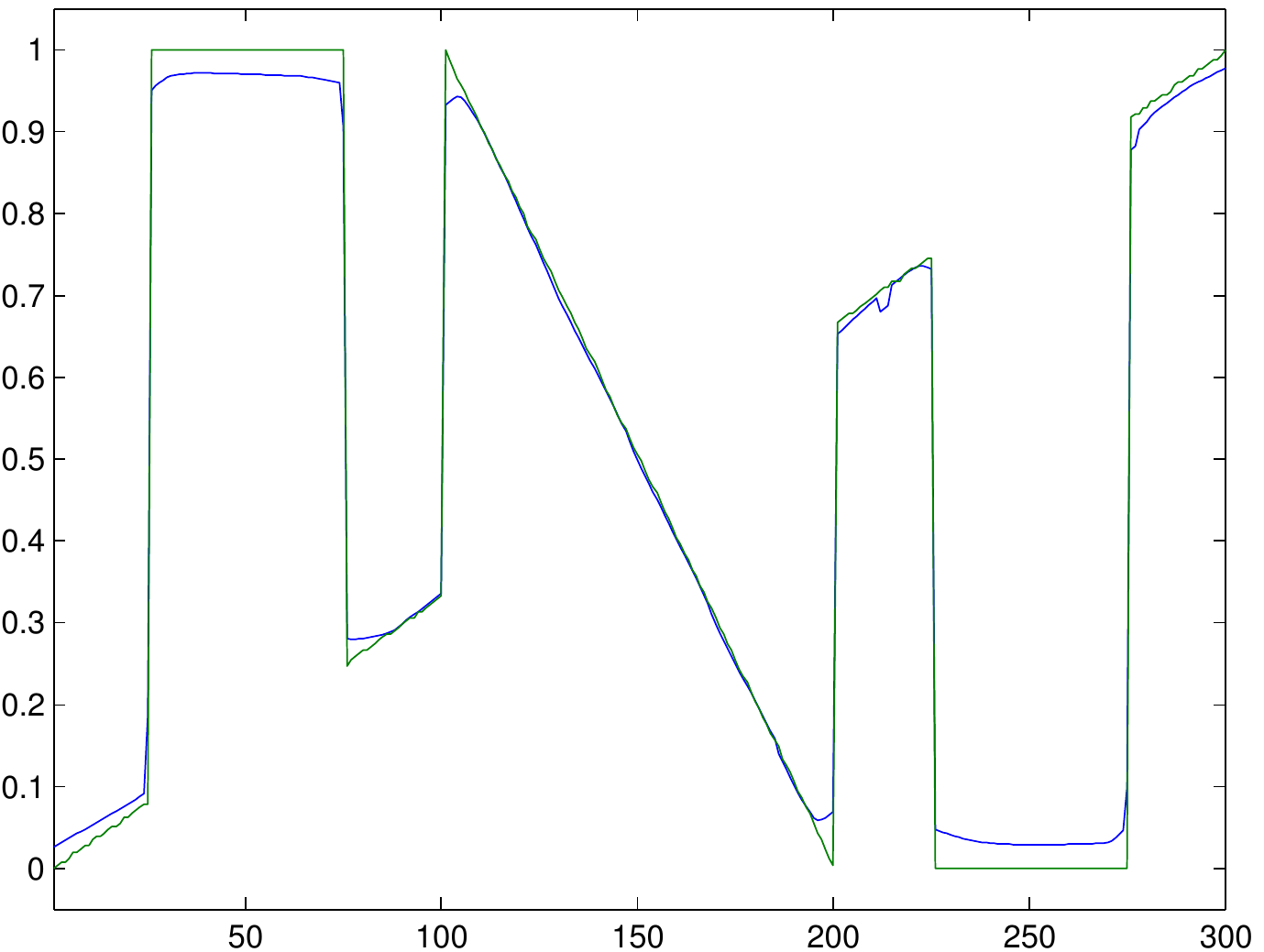} \\
input image & $L^2-\TV$ (ROF) & $L^2-\TGV$\\
\includegraphics[height=.2\linewidth]{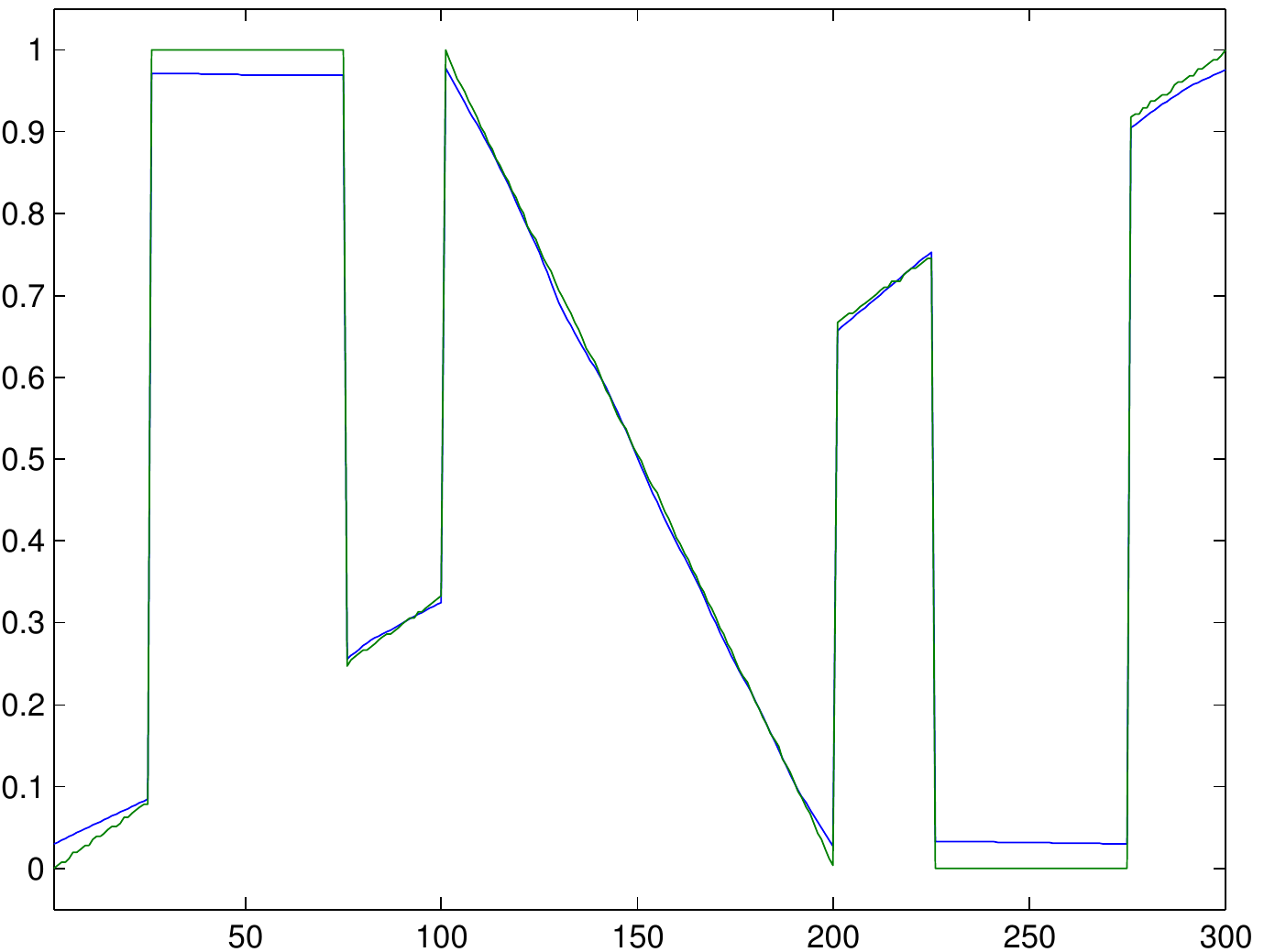} &
\includegraphics[height=.2\linewidth]{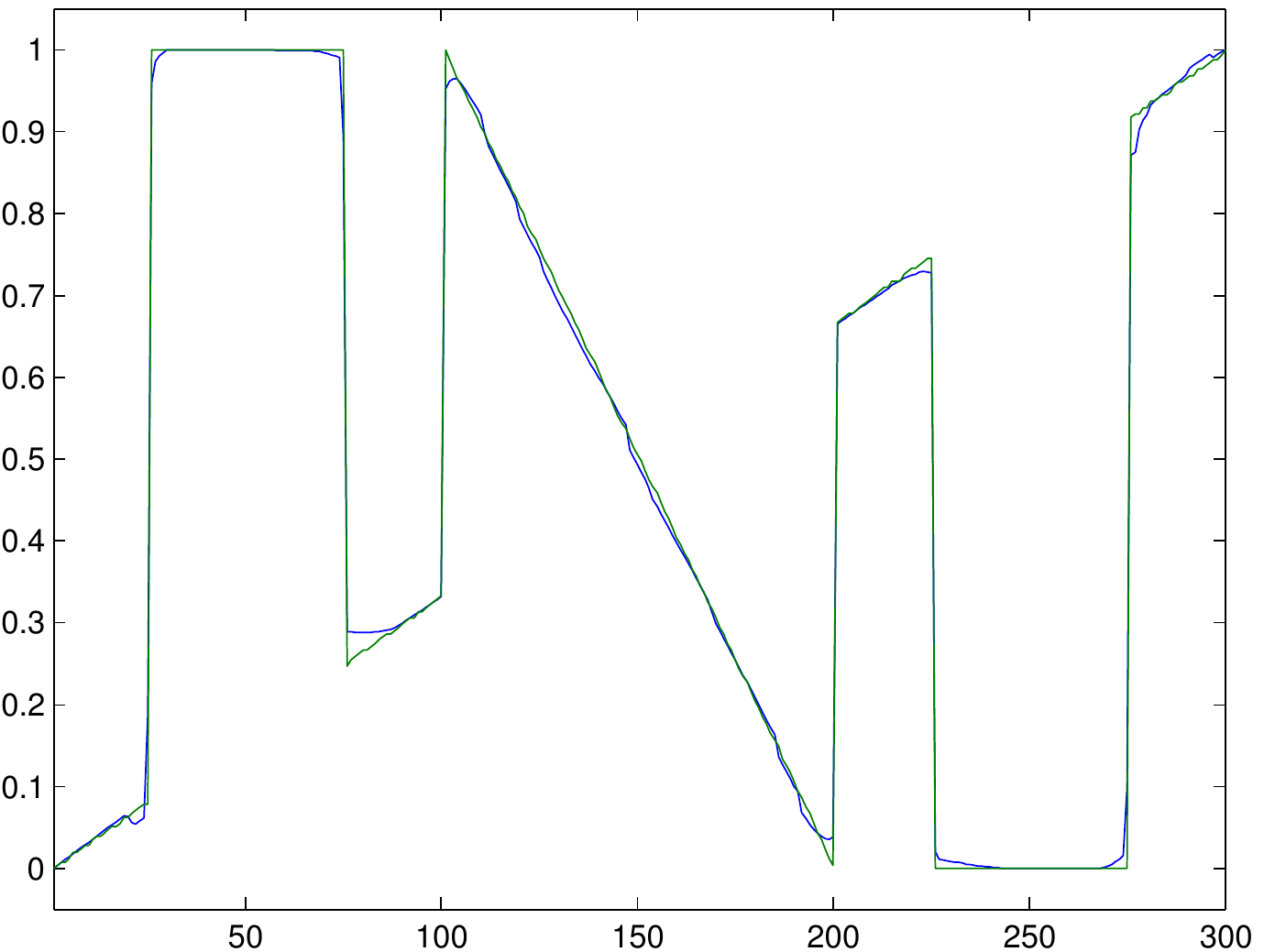} &
\includegraphics[height=.2\linewidth]{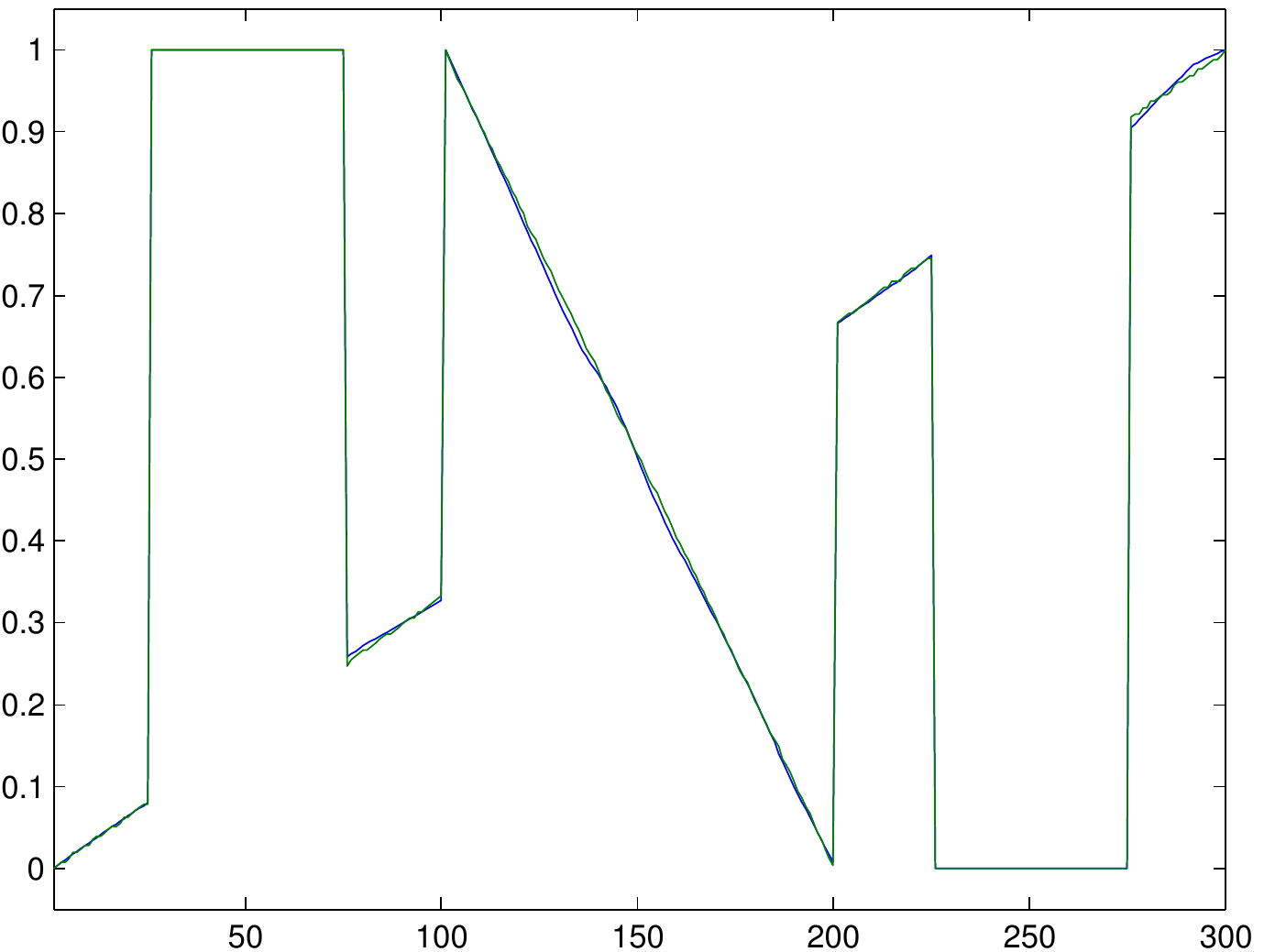}\\
$L^2-$non-local Hessian & $L^1-\TGV$ &$L^1-$ non-local Hessian \\
\end{tabular}
\ec
\caption{Slice at row $50$ through the results images in Figure~\ref{fig:shapes-img}. The $L^2-\TV$ model shows staircasing as expected.
$L^2-\TGV$ smoothes out sharp local extrema, which is preserved better by the non-local $L^2$ model.
Changing the data term to $L^1$ additionally removes the slight contrast reduction, even though the original noisy data is clipped to the interval $[0,1]$.
}\label{fig:shapes-cut}
\end{figure}

Figures~\ref{fig:shapes-img} and \ref{fig:shapes-cut} show a quantitative comparison of our non-local method with the results of the TGV approach and several classical regularizers. The parameters for each method were chosen by a grid search to yield the best PSNR. For all images we also provide a more realistic ``structural similarity index'' (SSIM) \cite{wang2004image}. The non-local approach improves on $\TGV$ with respect to PSNR (34.09 vs.~33.28). However, it is interesting to look at the spatial
distribution of the error: the parts where the non-local approach improves are exactly the local maxima, which are smoothed over by $\TGV$.
Surprisingly, this is hardly visible when looking at the images only (Figure~\ref{fig:shapes-img}), which is in line with the good results obtained with $\TGV$ for natural images.
However, this property might become a problem when the data is not a natural image, for example in the case of depth images or digital elevation maps. We refer to \cite{Lellmann2013a}
for a discussion of a problem that is more demanding in terms of the correct choice of the regularization.

A current limitation is that our choice of weights can result in a pixelated structure along slanted edges. This can happen when neighboring points are separated by a strong edge and therefore have no regularization in common. We conjecture that this effect could be overcome by ensuring that for neighbouring points $x,y$ the weight $\sigma_x(y)$ is always at least some positive constant, however we have not evaluated this direction further. Whether this is a desired behaviour or not depends on the underlying data -- for applications such as segmentation a pixel-wise decision might actually be preferable.

\begin{figure}[tp]
\bc
\begin{tabular}{ccc}
\includegraphics[height=.31\linewidth]{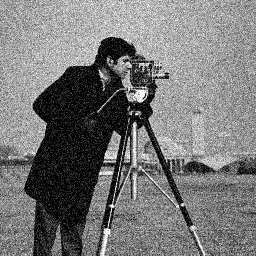} &
\includegraphics[height=.31\linewidth]{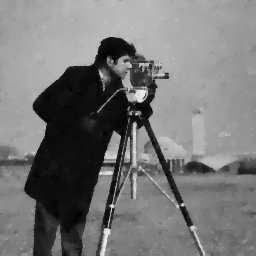}\\
input image & TGV\\
& $\alpha_0 = 0.30, \alpha_1 = 0.17$ \\
\includegraphics[height=.31\linewidth]{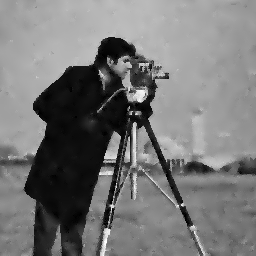} &
\includegraphics[height=.31\linewidth]{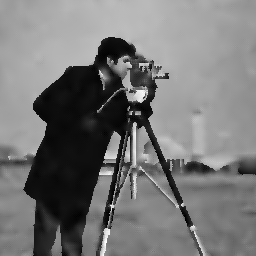} \\
$L^2-H'_u$ (1 iteration) & $L^2-H'_u$ (5 iterations)\\
$\alpha = 0.2$, $\gamma = 0.01$& $\alpha = 0.2$, $\gamma = 0.01$\\
\end{tabular}
\ec
\caption{Non-local Hessian-based denoising with $L^2$ data fidelity ($p=2$) on a real-world noisy input image (top left) with $\sigma=0.1$ Gaussian noise. Using a suitable choice of weights with a non-parametric neighborhood, small details and long, thin structures such as the camera handle can be preserved while removing noise in uniform structures (bottom left). By repeating the process several times, the result can be further improved (bottom right). Compared to TGV (upper right), edges are generally more pronounced, resulting in a slightly cartoon-like look.}
\label{fig:cameraman}
\end{figure}

Finally, Figure~\ref{fig:cameraman} shows an application to the ``cameraman''\ image with $L^2$ data term. Small details on the camera as well as long, thin structures such the camera handle and the highlights on the tripod are well preserved. In larger, more homogeneous areas the result is what would be expected from second-order smoothness.

In its basic form, our approach is data-driven, i.e., the weights are computed directly from the noisy data, while ideally they should be computed from the noise-free ground truth. We can approximate this \emph{solution-driven} approach by iterating the whole process, each time recomputing the weights from the previous solution.

\begin{figure}[tp]
\bc
\begin{tabular}{ccc}
\includegraphics[height=.2\linewidth]{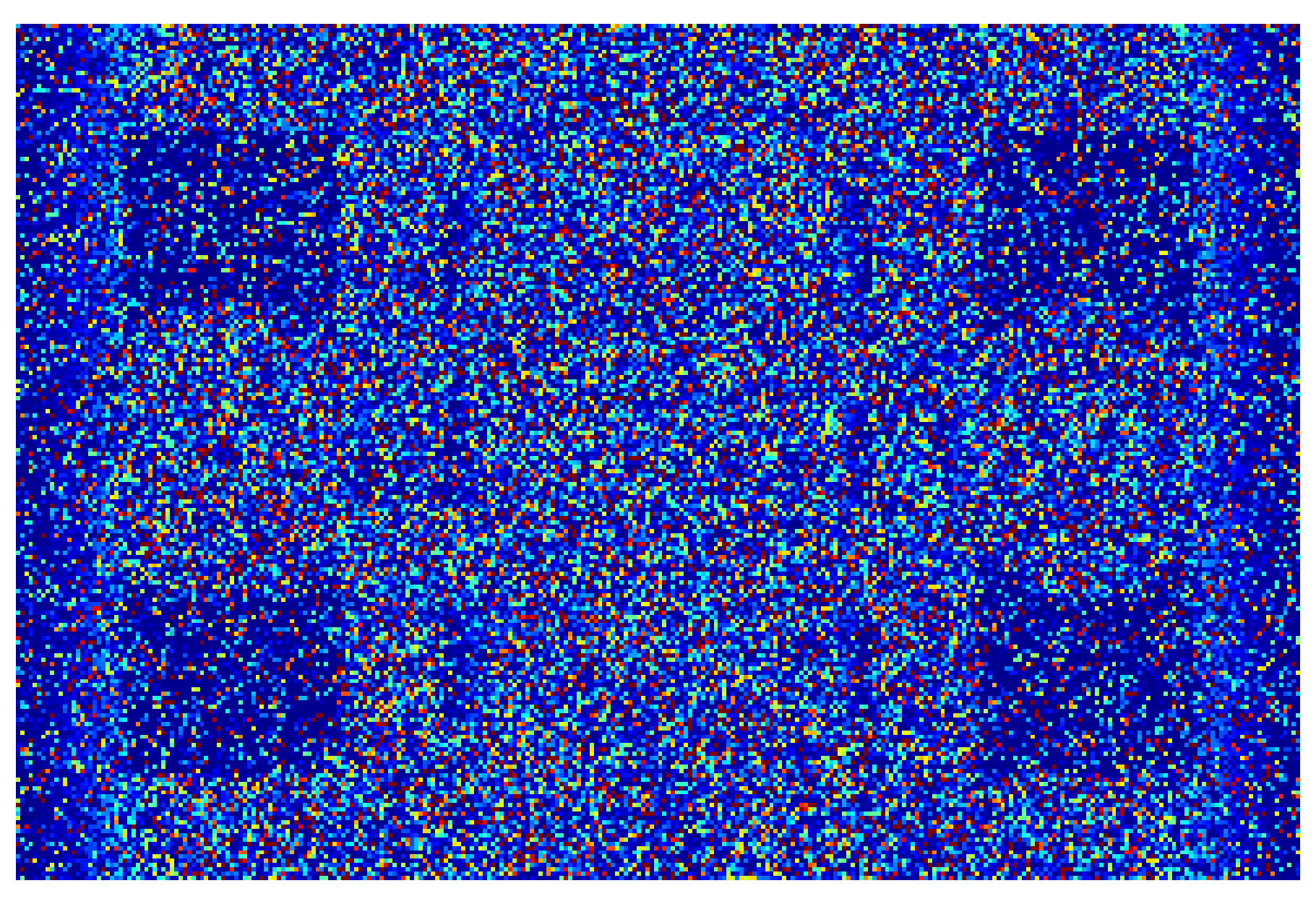} &
\includegraphics[height=.2\linewidth]{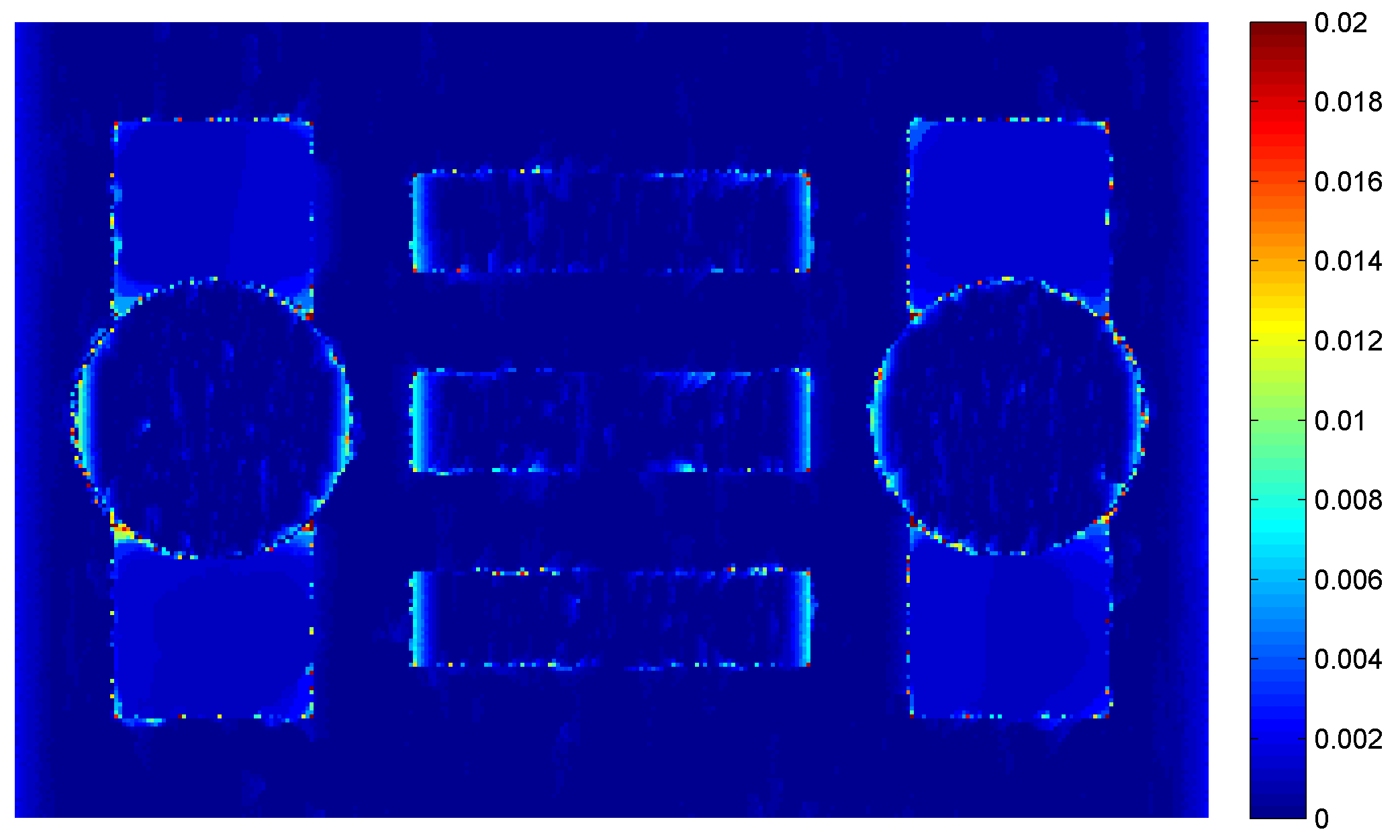} &
\includegraphics[height=.2\linewidth]{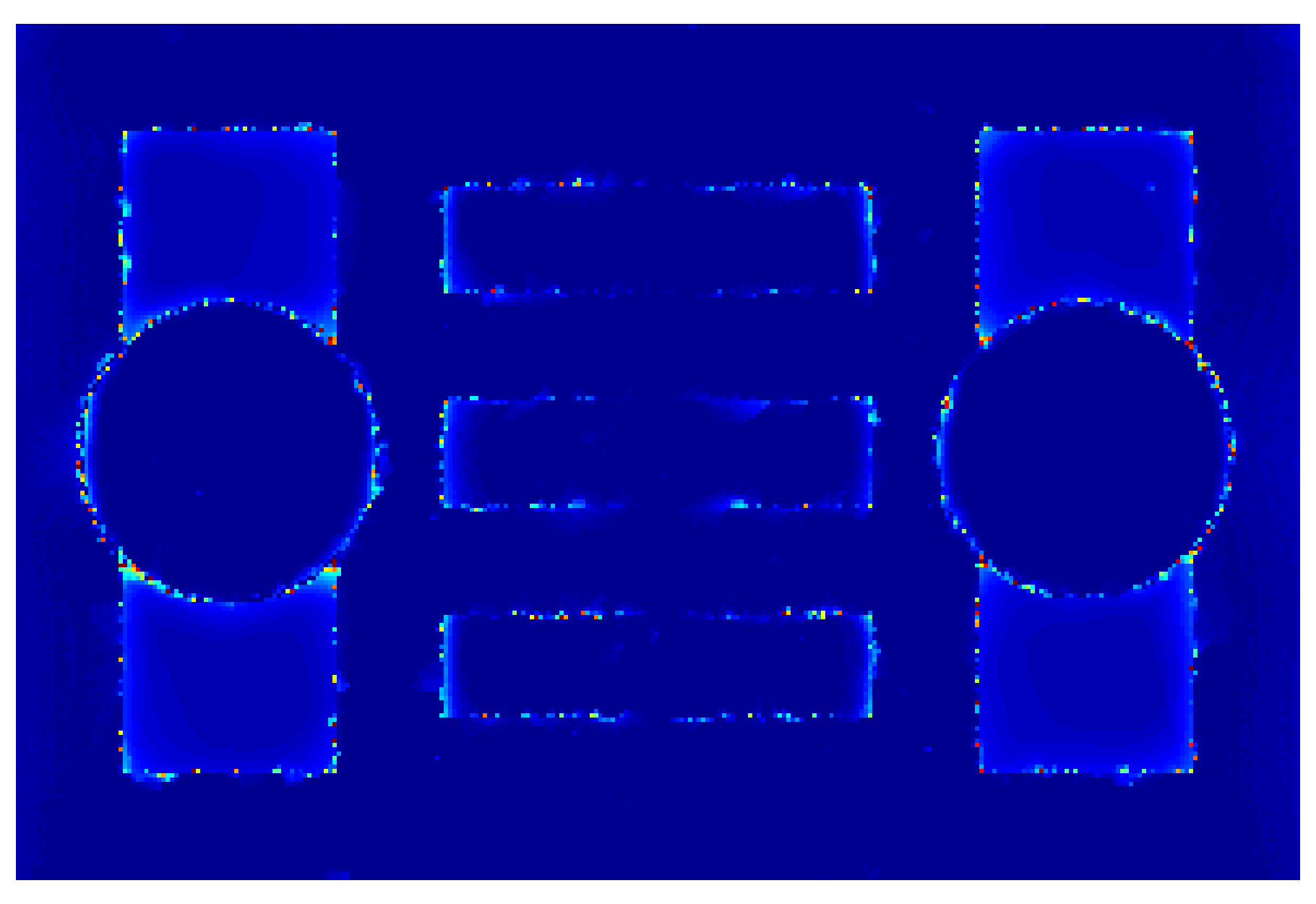} \\
input image & $L^2-\TV$ (ROF) & $L^2-\TGV$\\
\includegraphics[height=.2\linewidth]{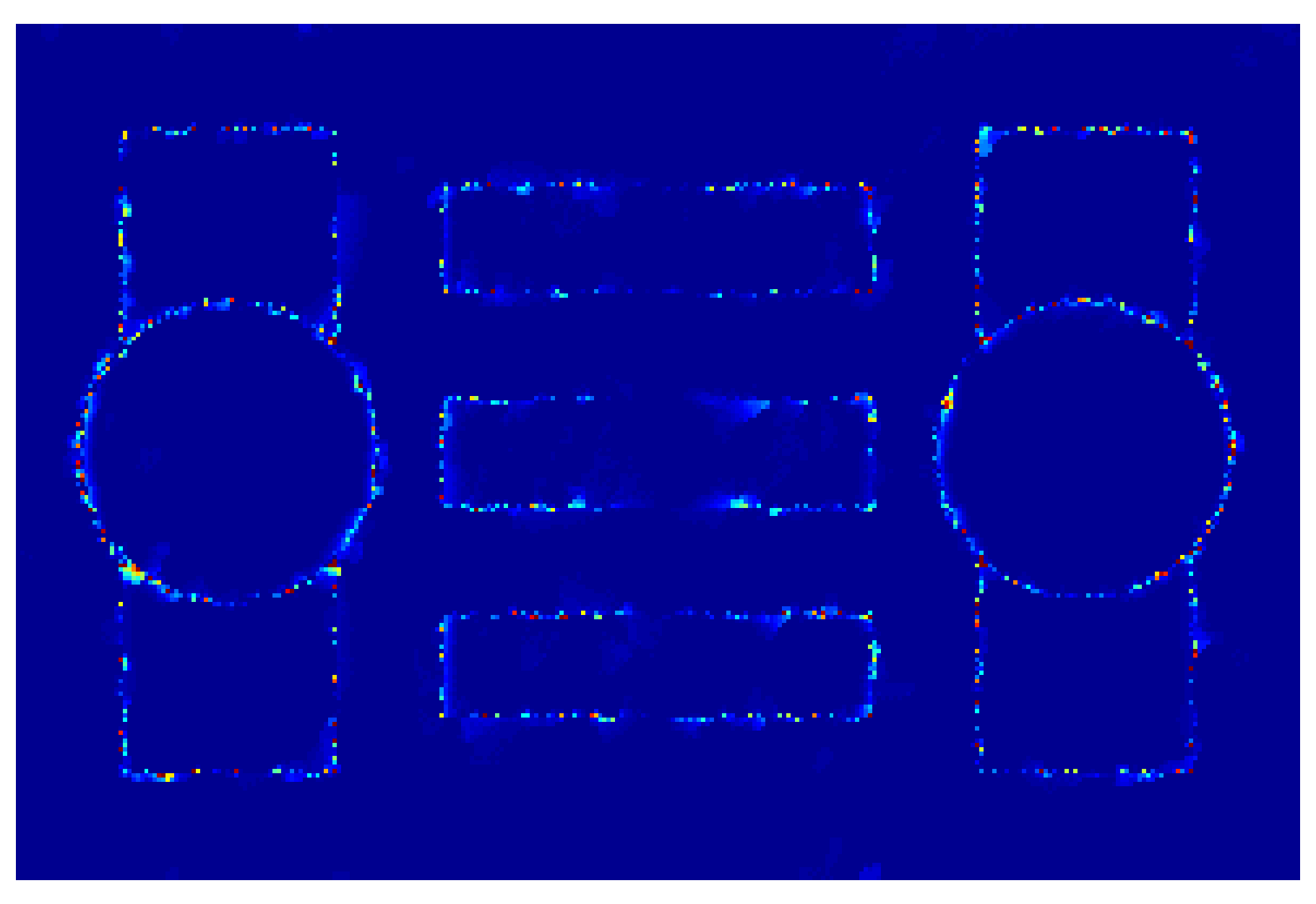} &
\includegraphics[height=.2\linewidth]{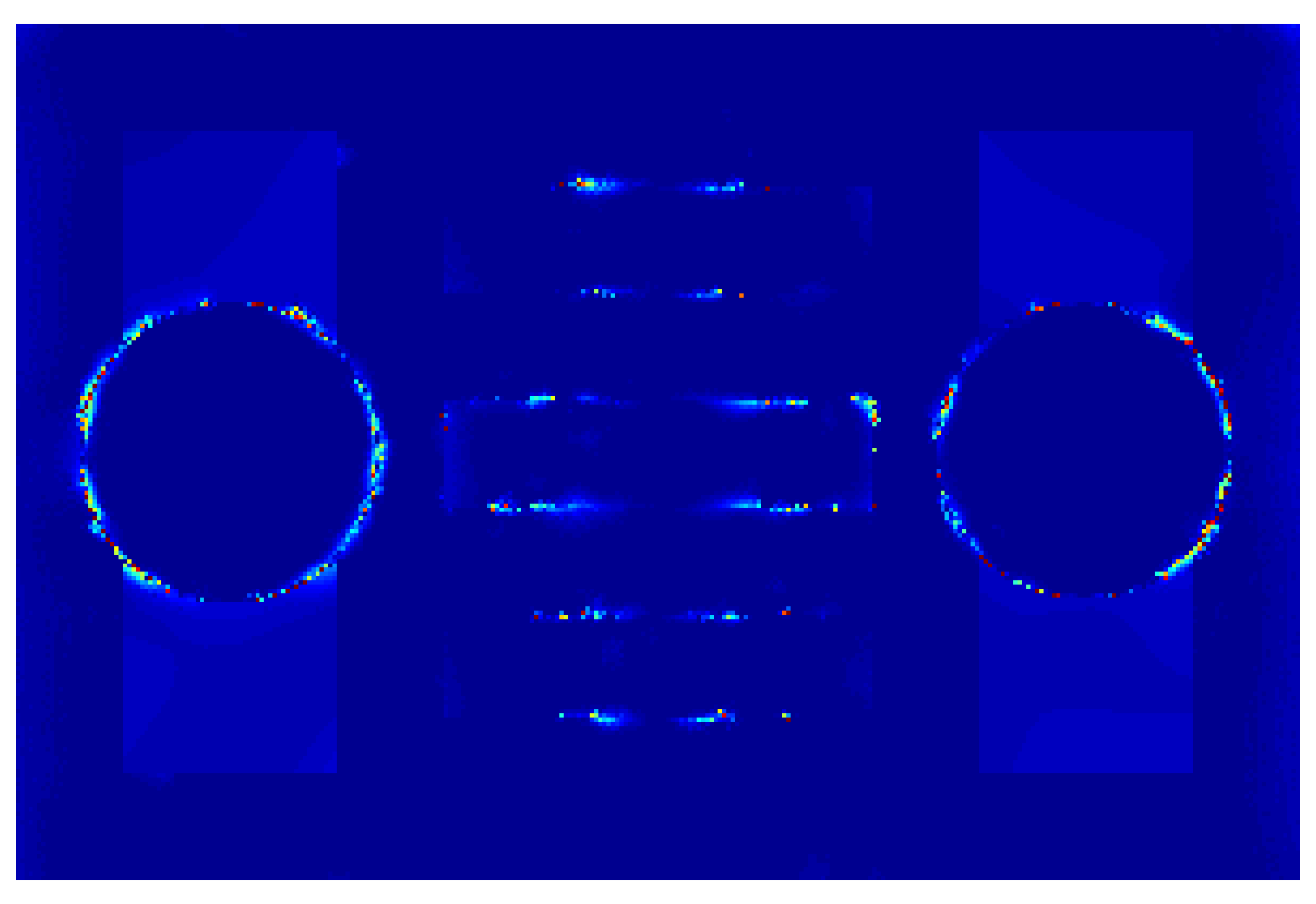} &
\includegraphics[height=.2\linewidth]{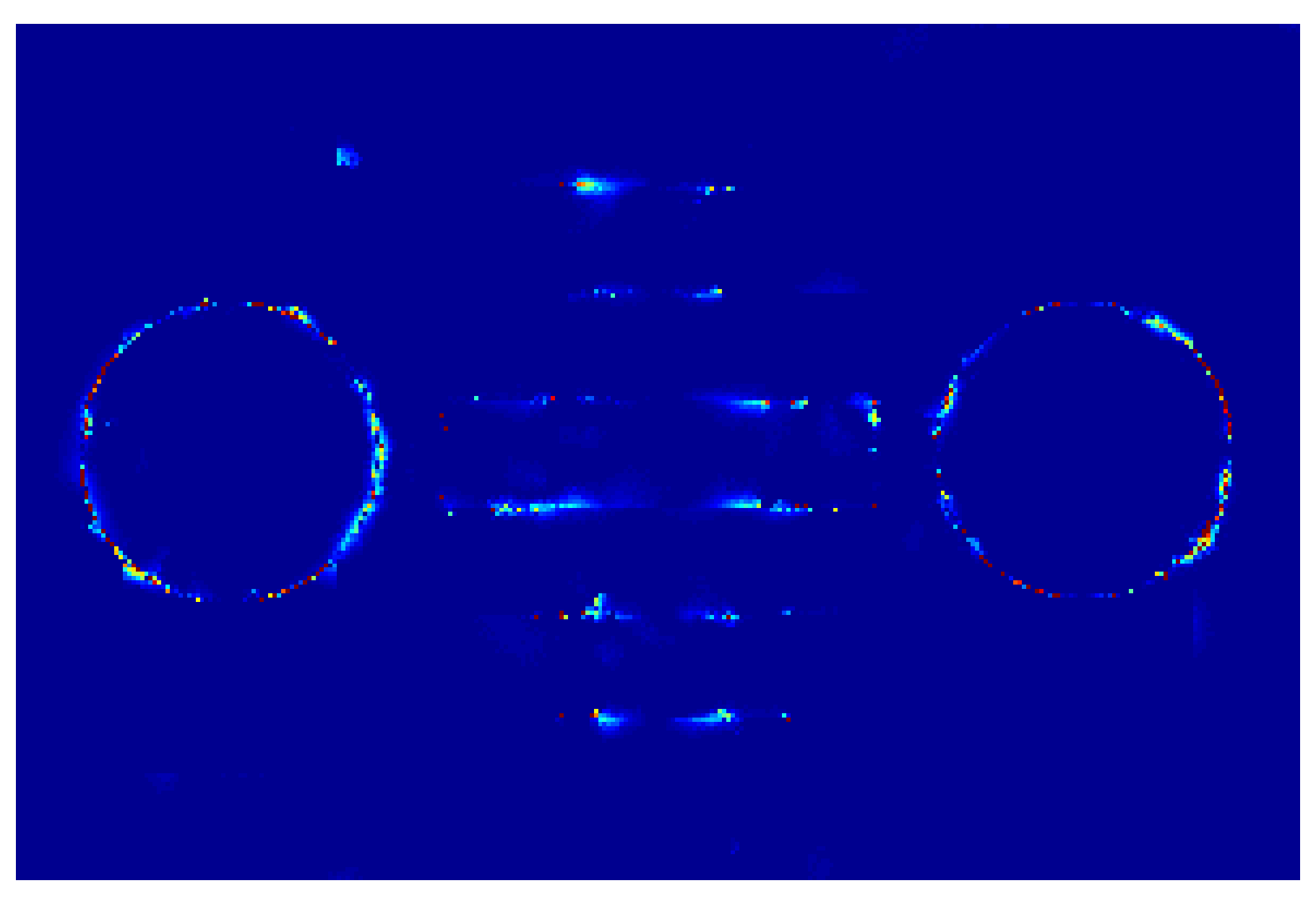} \\
$L^2-$non-local Hessian & $L^1-\TGV$ & $L^1-$ non-local Hessian \\
\end{tabular}
\ec
\caption{Distribution of the $L^2$ error to the ground truth for the experiment in Figure~\ref{fig:shapes-img}.
In the $\TGV$ case the errors spread out from the jumps. This is greatly reduced when using the non-local Hessian,
as the cost of introducing a few more local errors.
}\label{fig:shapes-error}
\end{figure}

\section{Conclusion}

From the perspective of the analysis, the study of the non-local Hessian is a natural extension of the study of non-local gradients.  While one has the straightforward observation we mention in the introduction -- that one may use non-local gradients to characterize higher-order Sobolev spaces -- the results of this paper alongside those of \cite{Spector} provide a framework for general characterizations of Sobolev spaces of arbitrary (non-negative integer) order.  This notion, that one can write down an integrated Taylor series and use this as a definition of a non-local differential object, is quite simple and yet leaves much to be explored.  Let us mention several interesting questions to this effect that would add clarity to the picture that has developed thus far.  

One of the foundations of the theory of Sobolev spaces is that of Sobolev inequalities.  While for a related class of functionals, a Poincar\'e inequality has been established by Augusto Ponce \cite{ponce2004anestimate}, the lack of monotonicity of the non-local gradient has proven difficult in adapting his argument to our setting.  This motivates the following open question.

{\bf Open Question:}
Can one find hypothesis on the approximation of the identity $\rho_n$ so that there is a $C>0$ such that for all $u$ in a suitable space one has the inequality
\[ \int_{\Omega} | u - \fint u|^p\;dx \leq C \int_{\Omega} |G_{\!\!\;n} u|^p\,dx,\]
for all $n$ sufficiently large?

In fact, beyond the multitude of questions one could explore by making a comparison of the results for non-local gradients and Hessians with known results in the Sobolev spaces, there are already several interesting questions for non-local gradients and Hessians in the regime $p=1$ that have not been satisfactorily understood from our perspective.  For instance, in the first order setting, assuming $u \in \mathrm{BV}(\Omega)$, one has the convergence of the total variations 
\[\int_{\Omega} |G_{\!\!\;n} u|\;dx  \to |Du|(\Omega).\]
While we have been able to show such a result for $H^{ij}_nu$, $i\neq j$, we have not succeeded in demonstrating this in the case $H^{ii}_nu$, which if true would settle the following conjecture.

{\bf Conjecture:}
Suppose $u \in \mathrm{BV}^2(\mathbb{R}^N)$.  Then
\[\int_{\mathbb{R}^N} |H_n u|\;dx  \to |D^2u|(\mathbb{R}^N).\]

This question is related to a larger issue needing clarification, which is that of the right assumptions for both the non-local gradient and non-local Hessian when $p=1$.  In both the paper \cite{Spector} and this paper, the assumption utilized in the characterizations has been that the object is an $L^1$ function.  The more natural assumption is that the non-local gradient or non-local Hessian exists as a measure, for which the same argument gives a characterization of $\mathrm{BV}(\Omega)$ or $\mathrm{BV}^2(\mathbb{R}^N)$.  However, the question of defining an integral functional of the non-local gradient or non-local Hessian is less clear.  More understanding is required here as to the notion of integral functionals of a measure in the local case and the right framework to place the non-local objects within.  

One can already see that from the analysis standpoint there are many interesting questions to explore in this area, while from the standpoint of applications we have seen that the relatively simple framework yields already quite interesting results. Here two interesting directions come into mind: developing strategies for choosing the weights, and moving on beyond second-order regularization.  

While in the numerical section we outlined one potential useful choice for the weights and an application for the non-local Hessian -- jump-preserving second-order regularization --, the approach is still somewhat local in the sense that we still restrict ourselves to a neighborhoods of each point. It would be interesting to find an application that allows to truly capitalize on  the fact that we can include far-reaching interactions together with an underlying higher-order model.

It would also be interesting to see whether it is useful in applications to go to even higher orders than $2$. While for natural images the use of such regularity is debatable, for other applications such as reconstructing digital elevation maps it can make the difference between success and failure \cite{Lellmann2013a}, and only requires replacing the quadratic model in the non-local Hessian \eqref{eq:implicitnl} by a cubic- or even higher-order model.

%

\section*{Appendix}

{\em Proof of Proposition \ref{prop:cij} }
  In the case $N = 1$ we always have $i_0=j_0$ and verify that in fact both the integral and the right-hand side are equal to $2$. For $N=2$ the integral becomes either
\beq
\int_0^{2 \pi} \cos^2 \alpha \sin^2 \alpha d \alpha = \frac{\pi}{4}\quad\text{or}\quad\int_0^{2 \pi} \cos^2 \alpha \cos^2 \alpha d \alpha = \frac{3 \pi}{4},
\eeq 
  which agrees with \eqref{eq:cijexplicit}. For the general case $N\geqs 3$, we apply the coarea formula in \cite[2.93]{AmbrosioBV} using $E=\Ss^{N-1}$, $M=N,$ $N=N-1$, $k=2$, $f ( x_{1} ,x_{2}
  , \ldots ) = ( x_{1} ,x_{2} )$, then
  \begin{eqnarray*}
    \int_{\Ss^{N-1}} g ( x ) C_{k} d^{E}  f_{x}  d \mathcal{H}^{N-1} ( x ) & = &
    \int_{\mathbbm{R}^{2}} \int_{\Ss^{N-1} \cap \{ y_{1} =x_{2} ,y_{2} =x_{2}
    \}} g ( y )  d \mathcal{H}^{N-3} ( y )  d  ( x_{1} ,x_{2} )\\
    & = & \int_{B_{1} ( 0 )} \int_{\sqrt{1-x_{1}^{2} -x_{2}^{2}} S^{N-3}} g (
    x_{1} ,x_{2} ,x_{r} )  d \mathcal{H}^{N-3} ( x_{r} )  d  ( x_{1} ,x_{2} )
    .
  \end{eqnarray*}
  In our case $g ( x ) =x_{1}^{2} x_{2}^{2}$, thus $g ( x_{1} ,x_{2} ,x_{r} )$
  is independent of $x_{r}$ and we obtain
  \begin{eqnarray*}
    \int_{\Ss^{N-1}} x_{1}^{2} x_{2}^{2} C_{k} d^{E}  f_{x}  d \mathcal{H}^{N-1}
    ( x ) & = & \int_{B_{1} ( 0 )} x_{1}^{2} x_{2}^{2} 
    \int_{\sqrt{1-x_{1}^{2} -x_{2}^{2}} \Ss^{N-3}} d \mathcal{H}^{N-3} ( x_{r} )
    d  ( x_{1} ,x_{2} )
  \end{eqnarray*}
  Consider $C_{k} d^{E} f_{x}$. For a given $x \in \Ss^{N-1}$, assume that $B
  \in \RN \times \R^{N-1}$ extends $x$ to an orthonormal
  basis of $\RN$. Then $x+Bt$, $t \in \R^{N-1}$,
  parametrizes the tangent plane of $E$ at $x$. The derivative of $f ( x+Bt )
  = ( x_{1} ,x_{2} )$ at $x$ in direction $t$ is
  \begin{eqnarray*}
    L \assign d^{E} f_{x} = \partial_{t} f ( x+Bt ) |_{t=0} & = & EB \nocomma
    , \hspace{1em} E \assign \left(\begin{array}{c}
      e_{1}^{\top}\\
      e_{2}^{\top}
    \end{array}\right)
  \end{eqnarray*}
  We are interested in $\sqrt{\det  L L^{\top}} = \sqrt{\det ( EBB^{\top} E
  )}$. Now
  \begin{eqnarray*}
    I_{2} & = & EI_{n} E^{\top}
     =  E \left(\begin{array}{cc}
      B & x
    \end{array}\right) \left(\begin{array}{c}
      B^{\top}\\
      x^{\top}
    \end{array}\right) E^{\top}
     =  EBB^{\top} E+Exx^{\top} E^{\top} ,
  \end{eqnarray*}
  thus
  \begin{eqnarray*}
    LL^{\top} =EBB^{\top} E & = & I_{2} -Exx^{\top} E^{\top}
     =  I_{2} - \left(\begin{array}{c}
      x_{1}\\
      x_{2}
    \end{array}\right)  \left(\begin{array}{cc}
      x_{1} & x_{2}
    \end{array}\right)
     =  \left(\begin{array}{cc}
      1-x_{1}^{2} & -x_{1} x_{2}\\
      -x_{1} x_{2} & 1-x_{2}^{2}
    \end{array}\right),
  \end{eqnarray*}
  and consequently
  \begin{eqnarray*}
    C_{k} d^{E} f_{x} & = & \sqrt{\det ( LL^{\top} )}
     =  \sqrt{( 1-x_{1}^{2} )  ( 1-x_{2}^{2} ) -x_{1}^{2} x_{2}^{2}}
     =  \sqrt{1-x_{1}^{2} -x_{2}^{2}} .
  \end{eqnarray*}
  Consequently
  \begin{eqnarray*}
    \int_{\Ss^{N-1}} x_{1}^{2} x_{2}^{2}  ( 1-x_{1}^{2} -x_{2}^{2} )  d
    \HH^{N-1} ( x ) & = & \int_{B_{1} ( 0 )} x_{1}^{2} x_{2}^{2} 
    \int_{\sqrt{1-x_{1}^{2} -x_{2}^{2}} \Ss^{N-3}} d \mathcal{H}^{N-3} ( x_{r} )
    d  ( x_{1} ,x_{2} )
  \end{eqnarray*}
  By putting the $C_{k}$ term into the denominator of $g$, we obtain the more
  useful form
  \begin{eqnarray*}
    \int_{\Ss^{N-1}} x_{1}^{2} x_{2}^{2}  d \mathcal{H}^{N-1} ( x ) & = &
    \int_{B_{1} ( 0 )} \sqrt{1-x_{1}^{2} -x_{2}^{2}}^{-1} x_{1}^{2} x_{2}^{2} 
    \int_{\sqrt{1-x_{1}^{2} -x_{2}^{2}} S^{N-3}} d \mathcal{H}^{N-3} ( x_{r} )
    d  ( x_{1} ,x_{2} ) .
  \end{eqnarray*}
  The right-hand side is
  \begin{eqnarray*}
    &  & \int_{0}^{1} \int_{0}^{2 \pi} r \sqrt{1-r^{2}}^{-1} r^{4} \cos^{2}
    \alpha \sin^{2} \alpha \left| \sqrt{1-r^{2}} \Ss^{N-3} \right|  d r d 
    \alpha\\
    & = & \int_{0}^{2 \pi} \cos^{2}   \alpha \sin^{2}   \alpha  d  \alpha
    \int_{0}^{1} \sqrt{1-r^{2}} r^{5}  \left| \sqrt{1-r^{2}} \Ss^{N-3} \right| 
    d r \\
    & = & \frac{\pi}{4}  \int_{0}^{1} \sqrt{1-r^{2}}^{-1} r^{5}  \left|
    \sqrt{1-r^{2}} \Ss^{N-3} \right|  d r\\
    & = & \frac{\pi}{4}  \int_{0}^{1} \sqrt{1-r^{2}}^{-1} r^{5} 
    \sqrt{1-r^{2}}^{N-3} | \Ss^{N-3} |  d r\\
    & = & \frac{\pi}{4}  | \Ss^{N-3} |  \int_{0}^{1} r^{5} 
    \sqrt{1-r^{2}}^{N-4}  d r\\
    & = & \frac{\pi}{4}  | \Ss^{N-3} |  \frac{8}{N (N+2 )  (N-2 )}\\
    & = & \frac{\pi}{4} \frac{(N-2 ) \pi^{(N-2 ) /2} }{\Gamma ( 1+ (N-2 )
    /2 )}  \frac{8}{N (N+2 )  (N-2 )}\\
    & = & \frac{2 \pi^{N/2} }{N (N+2 ) \Gamma (N/2 )} .
  \end{eqnarray*}
  Using the fact that $\int_{0}^{2 \pi} \cos^{4} ( \alpha )  d  \alpha =3
  \int_{0}^{2 \pi} \cos^{2} \alpha   \sin^{2} \alpha  d  \alpha$ we finally
  obtain
  \begin{eqnarray*}
    C_{i_{0}  j_{0}} = \int_{\Ss^{N-1}} \nu_{i_{0}}^{2} \nu_{j_{0}}^{2}  d
    \mathcal{H}^{N-1} ( x ) & = & \frac{ |S^{N-1}| }{N (N+2 )} \cdot \left\{ \begin{array}{ll}
      1, & i_{0} \neq j_{0} ,\\
      3, & i_{0} =j_{0} .
    \end{array} \right.
  \end{eqnarray*}
  This completes the proof.
\endproof

\bibliographystyle{amsalpha}

\bibliography{kostasbib}

\end{document}

%% file: macros.tex
\newcommand{\assign}{:=}

\newcommand{\bc}{\begin{center}}

\newcommand{\bea}{\begin{align}}
\newcommand{\beq}{\begin{equation}}
\newcommand{\beqa}{\begin{eqnarray}}
\newcommand{\bi}{\begin{itemize}}
\newcommand{\bitem}{\begin{itemize}}
\newcommand{\bpm}{\begin{pmatrix}}      
\newcommand{\BV}{\tmop{BV}}

\newcommand{\cref}[1]{ {\tiny[{#1}]}}

\newcommand{\dHHN}{d \HH^{N-1}}

\newcommand{\ec}{\end{center}}
\newcommand{\eea}{\end{align}}
\newcommand{\eeq}{\end{equation}}
\newcommand{\eeqa}{\end{eqnarray}}
\newcommand{\ei}{\end{itemize}}
\newcommand{\eitem}{\end{itemize}}

\newcommand{\epm}{\end{pmatrix}}

\newcommand{\geqs}{\geqslant}

\newcommand{\HH}{\mathcal{H}}

\newcommand{\mathbbm}{\mathbb}

\newcommand{\MM}{\mathcal{M}}

\newcommand{\N}{\mathbb{N}}

\newcommand{\nocomma}{}

\newcommand{\Om}{\Omega}

\newcommand{\R}{\mathbb{R}}

\newcommand{\Rn}{\mathbb{R}^n}
\newcommand{\RN}{\mathbb{R}^N}
\newcommand{\RNN}{\mathbb{R}^{N \times N}}

\newcommand{\Ss}{\mathcal{S}}

\newcommand{\tm}[1]{}

\newcommand{\tmop}[1]{\ensuremath{\operatorname{#1}}}

\newcommand{\TV}{\tmop{TV}}
\newcommand{\TGV}{\tmop{TGV}}

\newcommand{\half}{\frac{1}{2}}